\documentclass[11pt,leqno,a4paper,psamsfonts]{amsart}
\usepackage[english]{babel}
\usepackage{amssymb,amsmath,amscd,enumerate,amsthm}
\usepackage[psamsfonts]{amsfonts}
\usepackage{latexsym}
\usepackage{amsbsy}
\usepackage[mathscr]{eucal}
\usepackage{graphicx}
\usepackage{enumerate}
\usepackage{verbatim}
\usepackage{diagrams}
\newarrow{aux}{C}{-}{-}{-}{-}
\newarrow{Into}{C}{-}{-}{-}{>}
\usepackage[all]{xy}
\newcommand{\C}{{\mathbb{C}}}
\newcommand{\F}{{\mathcal{F}}}
\newcommand{\R}{{\mathbb{R}}}

\newcommand{\uF}{\underline{\F}}

\newcommand{\br}[1]{{\breve{#1}}}

\usepackage{color}
\newcommand{\CS}{{\mathrm{CS}}}
\newcommand{\Sat}{{\mathrm{Sat}}}
\newcommand{\tF}[1]{{\wt{\mc F}_{#1}}}
\newcommand{\tQF}[1]{{\wt{\mc Q}^{\mc F}_{#1}}}
\newcommand{\tQFpr}[1]{{\wt{\mc Q}^{\mc F'}_{#1} }}
\let\tFQ\tQF

\newcommand{\QF}{{\wt{\mc Q}^{\mc F}_{\infty} }}

\newcommand{\QFpr}{{\wt{\mc Q}^{\mc F'}_{\infty} }}

\newcommand{\Hom}{\mathrm{Hom}}

\newcommand{\proan}{{\underleftarrow{\mathrm{An}}}}
\newcommand{\protop}{{\underleftarrow{\mathrm{Top}}}}

\newcommand{\aut}{{\mathrm{Aut}}}
\newcommand{\proaut}{{\mathrm{Aut}}_{\proan}}
\let\Aut\aut

\let\monF\monf

\newcommand{\monFS}{\mf M^{\F}_S}%_{\infty}}

\newcommand{\monFSpr}{\mf M^{\F'}_{S'}}

\let\monprFS\monFSpr
\newcommand{\Sing}{{\rm Sing}}

\let\Comp\comp
\newcommand{\svee}{{\scriptscriptstyle \vee}}

\newcommand{\mf}[1]{{\mathfrak{#1}}}
\newcommand{\mr}[1]{\mathrm{#1}}
\newcommand{\mc}[1]{{\mathcal{#1}}}
\newcommand{\mb}[1]{{\mathbb{#1}}}

\newcommand{\wt}[1]{{\widetilde{#1}}}

\newcommand{\too}{\longrightarrow}
\newcommand{\limproj}{\underleftarrow\lim}
\newcommand{\limind}{\underrightarrow\lim}

\let\ssstyle\scriptscriptstyle
\let\sstyle\scriptstyle

\newcommand{\iso}{{\overset{\sim}{\longrightarrow}}}

\newcommand{\inte}[1]{\overset{\circ}{{#1}}}
\newcommand{\un}[1]{\underline{#1}}

\let\CAL=\mathcal%
\def\mathcal#1{{\CAL#1}}%
\newtheorem{teo}{Theorem}[subsection]

\newtheorem{lema}[teo]{Lemma}

\newtheorem{sublema}[teo]{Sub-Lemma}
\newtheorem{prop}[teo]{Proposition}
\newtheorem{defin}[teo]{Definition}
\newtheorem{cor}[teo]{Corollary}
\newtheorem{obs2}[teo]{Remark}
\newtheorem{convs2}[teo]{Conventions}
\newtheorem{recap2}[teo]{R\'{e}capitulation}
\newtheorem{ex2}[teo]{Example}
\newenvironment{obs}{\begin{obs2}\rm}{\end{obs2}}

\newenvironment{dem}{\begin{proof}[Proof]}{\end{proof}}
\newenvironment{dem2}[1]{\begin{proof}[Proof #1]}{\end{proof}}

\newenvironment{convs}{\begin{convs2}\rm}{\end{convs2}}
\def\bibartp#1#2#3#4#5#6#7#8
{\bibitem{#1} {\sc #2}, {\it #3}, {#4}, {\rm v. #5}, pages { #6} to
{ #7}, {({#8})}}
\def\bibart#1#2#3#4#5#6
{\bibitem{#1} {\sc #2}, {\it #3}, {#4}, {\rm #5}, {({#6})}}
\def\bibliv#1#2#3#4#5
{\bibitem{#1} {\sc #2}, {\it #3}, {#4}, {({#5})}}
\def\bibaart#1#2#3#4
{\bibitem{#1} {\sc #2}, {\it #3}, {#4}}
\begin{document}
\title[Monodromy and topological classification]{Monodromy and topological classification of germs of holomorphic foliations}
\date{\today}
\author{David Mar\'{\i}n and Jean-Fran\c{c}ois Mattei}
\thanks{The first author was partially supported by the projects
 MTM2007-65122 and  MTM2008-02294 of the  Ministerio de Educaci\'on y Ciencia de Espa\~{n}a / FEDER}
\address{Departament de Matem\`{a}tiques \\
Universitat Aut\`{o}noma de Barcelona \\
E-08193 Bellaterra (Barcelona)\\
Spain} \email{davidmp@mat.uab.es}

\address{Institut de Math\'{e}matiques de Toulouse\\
Universit\'{e} Paul Sabatier\\
118, Route de Narbonne\\
F-31062 Toulouse Cedex 9, France} \email{jean-francois.mattei@math.univ-toulouse.fr}

\noindent


\begin{abstract}
We give a complete topological classification of germs of holomorphic foliations in the plane under rather generic conditions.
The key point is the introduction of a new topological invariant
called monodromy representation. This monodromy contains all the relevant dynamical information, in particular the projective holonomy representation whose topological invariance was conjectured in the eighties by Cerveau and Sad
and proved here under mild hypotheses.
\end{abstract}

\maketitle

%\tableofcontents
\section{Introduction}
The objective of this paper is to provide a complete topological classification of germs of singular non-dicritical holomorphic foliations $\F$ at $(0,0)\in\mb C^{2}$ under very generic conditions. To do this we introduce a new topological invariant which is a representation of the fundamental group of the complement of the separatrix curve into a suitable automorphism group. We shall call this representation the \emph{monodromy of the foliation germ}.

In fact, the motivation of this work was the following conjecture of D.~Cer\-veau and P.~Sad  in 1986, cf. \cite[page 246]{CerveauSad}. Consider two germs of foliations defined by germs of differential  holomorphic  $1$-forms $\omega$ and $\omega'$ at $(0,0)\in\C^{2}$.

\medskip

\begin{description}
\item[Conjecture (Cerveau-Sad) ]  \it  If $\omega $ and $\omega'$ are topologically conjugate and if $\omega$ is a generalized curve, then their respective projective holonomies representations are conjugate.
\end{description}

It was given in two forms, each of them with natural generic hypothesis concerning the germ of the foliation $\un\F$ along the \emph{exceptional divisor}  $\mc E_{\F}:=E_{\F}^{-1}(0)$ of the reduction  $E_\F : \mc B_\F\to \mb C^2$ of the singularity of $\F$, cf. \cite{Seidenberg, MatMou}.
The weak form (named Conjecture A) assumes that the separatrix curve is the union of smooth and transverse branches. In particular, $E_{\F}$ corresponds to a single blow-up. The strong form (named Conjecture B) only asks that the reduced foliation $\un \F=E_\F^{*}(\F)$ on $\mc B_\F$ does not have any
saddle-node singularity.

\medskip

Conjecture A was established  by one of us in \cite{Mar}. We give here an affirmative answer to Conjecture B. More precisely,  Theorem I below gives a list of topological invariants containing the projective holonomy representations. In turn Theorem II gives a complete topological classification.\\

It is worthwhile to stress here that through the whole paper all the topological {conjugations} between foliations that we consider are supposed to preserve the orientations of the ambient space and also the leaves orientations.\\

As in the situation considered by D. Cerveau and P. Sad, we restrict our attention to a reasonable class of foliations that are going to be called Generic General Type.
{Let $\F$ be a non-dicritical foliation, i.e. having a finite number $n$ of irreducible analytic germ curves $S_{1},\ldots,S_{n}$ invariant by $\F$, which are called \emph{separatrices}. Here it is worth to recall the celebrated Separatrix Theorem of \cite{CamachoSad} asserting that $n>0$. In the sequel we will call $S_{\F}:=\bigcup\limits_{i=1}^{n}S_{i}$ the \emph{separatrix curve} of $\F$.}
Following \cite{MarMat} we say that
%a non-dicritical\footnote{i.e. having a finite number $n$ of irreducible analytic germ curves $S_{1},\ldots,S_{n}$ invariant by $\F$, which are called \emph{separatrices}. Here it is worth to recall the celebrated Separatrix Theorem of \cite{CamachoSad} asserting that $n>0$.}
the foliation $\F$ is of \emph{General Type} if all the singularities of $\uF$ which are not linearizable are resonant, more precisely:
\begin{enumerate}
  \item[(GT)]\it for each singularity of
 $\underline{\F}$ there are local holomorphic coordinates  $(u,v)$ such that $\un{\F}$ is locally defined by a holomorphic 1-form of one of the two following types :
\begin{enumerate}[(i)]
\item\label{singredlin} $\lambda _1udv+\lambda _2vdu$, with $\lambda_{1}\lambda_{2}\neq 0$ and  $\lambda _1/\lambda _2\notin \mb Q_{<0} $ (linearizable singularity),
\item  $(\lambda _1u+\cdots )dv+(\lambda _2v+\cdots)du$, with $\lambda _1$, $\lambda _2\in \mb N^\ast$, (resonant saddle).
\end{enumerate}
\end{enumerate}
To introduce the additional genericity condition (G) we recall first that a singularity of $\un{\F}$ is of \emph{nodal type} if it can be locally written as
$$(\lambda _1u+\cdots )dv+(\lambda _2v+\cdots)du,$$ with $\lambda_{1}\lambda_{2}\neq 0$ and $\lambda_{1}/\lambda_{2}\in\R_{<0}\setminus\mb Q_{<0}$. Such  singularities are always linearizable and consequently the only local analytical invariant of a node is its Camacho-Sad index $-\lambda_{1}/\lambda_{2}$.
The topological specificity of a nodal singularity $s$ is the existence, in any small neighborhood of $s$, of a saturated closed set whose complement is an open disconnected neighborhood of the two punctured local separatrices of the node. We call \emph{nodal separator} such a saturated closed set.
We denote by $\mr{Node}(\un{\F})$ the set of nodal singularities of $\un{\F}$.
With this notations the genericity condition can be stated as  follows:
\begin{itemize}
  \item[(G)] \it The closure of each connected component of $\mc E_{\F}\setminus(\mr{Node}(\un{\F})\cap\mr{Sing}(\mc E_{\F}))$ contains  an irreducible component of the exceptional divisor $\mc E_{\F}$
  having a non solvable holonomy group.
\end{itemize}

\noindent Notice that when $\mr{Node}(\un{\F})=\emptyset$, the genericity condition (G) only asks for a single irreducible component of $\mc E_{\F}$ having a non solvable holonomy group.
\noindent In the space of coefficients of the germ of holomorphic $1$-form defining the foliation this condition is generic in the sense of the Krull topology, cf. \cite{LeFloch}.
\\

\noindent A foliation satisfying Conditions
(G) and (GT) above will be called \emph{Generic General Type}.
For such a foliation $\F$, if $\mr{Node}(\un{\F})=\emptyset$ Theorem~I below provides a list of topological invariants. In the case that $\mr{Node}(\un{\F})\neq\emptyset$ we must restrict the class of topological {conjugations} in order to keep their invariance.
In fact, the first version \cite{MarMatMon} of this work dealt only with Generic General Type foliations $\F$ satisfying the additional requirement $\mr{Node}(\un{\F})=\emptyset$. Here this hypothesis is eliminated by modifying slightly the statements and the proofs given in \cite{MarMatMon}. In practice, this is  done by adding a prefix $\mc N$- to some notions whose new meaning is made precise when they appear for the first time.
We recommend the reader to ignore all the prefix $\mc N$- in a first reading.\\

\noindent\textbf{Definitions for the nodal case.}
\begin{itemize}
   \item A \emph{nodal separatrix of $\F$} is a separatrix whose strict transform  by the reduction map $E_{\F}$ meets the exceptional divisor  at a nodal singular point of $\un{\F}$.
        \item A \emph{$\mc N$-separator of $\un\F$} is
          the union of
         a system of nodal separators, one for each point in $\mr{Node}(\un{\F})\cap\mr{Sing}(\mc E_{\F})$ jointly with some neighborhoods of the strict transforms of the nodal separatrices of $\F$.
        A \emph{$\mc N$-separator of $\F$} is
          the image by $E_{\F}$ of a $\mc N$-separator of $\un\F$. If $\mathrm{Node}(\uF)=\emptyset$, a $\mc N$-separator is the emptyset.
   \item A \emph{$\mc N$-topological {conjugation} between two foliation germs $\F$ and $\F'$} is a germ of homeomorphism $h$ preserving the orientation of the ambient space as well as the orientation of the leaves, which is a {topological} {conjugation} between  $\F$ and $\F'$, such that for each nodal separatrix $S_{j}$ of $\F$, $h(S_{j})$ is a nodal separatrix of $\F'$ and the Camacho-Sad indices of $\un{\F}$ and $\un{\F'}$ along the strict transforms of $S_{j}$ and $h(S_{j})$ coincide.
   \item A \emph{$\mc N$-transversely holomorphic {conjugation} } between $\F$ and $\F'$ (resp. $\uF$ and $\uF'$) is a $\mc N$-topologically {conjugation} between these foliations,  which is transversely holomorphic on the complementary of some $\mc N$-separator of $\F$ (resp. $\un\F$).
\end{itemize}
Clearly the notions of $\mc N$-topological {conjugation}  and $\mc N$-transversely holomorphic {conjugation}  coincide with the usual notions of topological {conjugation}  and transversely holomorphic {conjugation}, when $\mr{Node}(\uF)=\emptyset$. In Section~\ref{invcssepartr}, cf. Remark~\ref{Ntop}, we shall prove that:
\begin{itemize}
\item{\it Any topological {conjugation} which is  transversely holomorphic in a neighborhood of each nodal separatrix minus the origin is a $\mc N$-topo\-lo\-gical {conjugation}.}
\end{itemize}
In particular, any transversely holomorphic {conjugation} is a $\mc N$-topological {conjugation}.
In order to assure the transverse holomorphy of a {conjugation} we shall use a generalized form of the following theorem of J. Rebelo \cite{Rebelo}:\\

\noindent{\bf Transverse Rigidity Theorem.} \textit{Every topological {conjugation} between two germs of non-dicritical holomorphic foliations satisfying the genericity condition  \hbox{\rm(G)} and having  singularities, after reduction, of type  $(\lambda _1u+\cdots)dv+(\lambda _2v+\cdots)du$ with  $\lambda _1\lambda _2\neq 0$, $\lambda _1/\lambda _2\notin \mb R_{<0}$, is transversely conformal.}\\

In fact the proof provided in \cite{Rebelo} shows that if we allow nodal singularities then each connected component of $\mc E_{\F}\setminus(\mr{Node}(\un{\F})\cap\mr{Sing}(\mc E_{\F}))$ possesses an open neighborhood $W$ such that the restriction of the topological {conjugation} to $E_{\F}(W)\setminus\{0\}$ is transversely conformal.
The extended version of the Transverse Rigidity Theorem asserts:
\begin{itemize}
  \item[(TRT)] \it For any orientations preserving topological  {conjugation} $\Phi$ between two germs of non-dicritical generalized curves satisfying condition (G) we have that $\Phi$ is a $\mc N$-topological {conjugation} if and only if $\Phi$ is a $\mc N$-transversely holomorphic {conjugation}.
\end{itemize}

\bigskip

\noindent{\bf Theorem I.} \textit{{For every non-dicritical Generic General Type foliation $\F$, the analytic type of the projective holonomy representation of each irreducible component of the exceptional divisor $\mc E_{\F}$ is a topological invariant when $\mr{Node}(\un\F)=\emptyset$.
More generally, the semilocal data $\mc S\mc L(\F)$}
%the family
%$$\mc S\mc L(\F) := \Big([S_{\F}]^{\mathrm{top}},\; (\CS(\un\F, D,s))_{s,\, D}\,,\; ([\un\F_{s}]^{\mathrm{hol}})_s\,, ([\mc H_{\un\F,\,D}]^{\mathrm{hol}})_D  \Big)$$
constituted by
\begin{itemize}
\item the  topological type %$[S_{\F}]^{\mathrm{top}}$
of the embedding of the total separatrix  curve %\footnote{i.e. the union $\cup_{i=1}^{n}S_{i}$ of all the separatrices of $\F$.}
    $S_{\F}$ of $\F$ into $(\C^{2},0)$,
%\item the Camacho-Sad index $\CS(\un\F, D,s)$ of the reduced foliation  $\un \F$ along  each irreducible component $D$ of the exceptional divisor $\mc E_\F:=E_\F^{-1}(0)$ and for each singularity $s$ of $\un \F$ over $D$,
  \item the {collection of} local analytic types $[\un\F_{s}]^{\mathrm{hol}}$
  {of the reduced foliation} $\underline{\F}$ at each singular point $s\in \Sing(\un\F)$, {codifying in particular  the Camacho-Sad index  $\CS(\un\F, D,s)$ of $\un\F$ at every singular point $s$ along each irreducible component $D$ of $\mc E_{\F}$ containing $s$,}
  \item the analytic type  %$[\mc H_{\un\F,\,D}]^{\mathrm{hol}}$
  of the holonomy representation $\mc H_{\un\F,\,D}$ of each {irreducible} component  $D$ of $\mc E_\F$,
\end{itemize}
is a $\mc{N}$-topological invariant of the germ of $\F$ at $0\in\mb C^2$.}
\\

{Notice that the Camacho-Sad index $\CS(\un\F,D,s)$ determine the analytic type of $\un\F$ at $s$ when $s$ is not a resonant singularity after the assumption (GT). On the other hand, }
% Notice that
the genericity condition (G) is strictly necessary in Theorem~I. Indeed, inside the family of homeomorphisms $\Psi(x,y)=(x|x|^{a},y|y|^{b})$ there is a topological {conjugation} between any pair of  linear hyperbolic singularities having different Camacho-Sad indices.\\

Theorem I asserts that
$\mc S\mc L(\F)$  is a topological invariant for the class of Generic General Type foliations with $\mr{Node}(\un{\F})=\emptyset$.
\begin{comment}
In order to give a precise sense to the above statement, we need to explain the meaning of the term \emph{topological invariant}. Two  germs $\F$ and $\F'$ are \emph{topologically conjugate} if there exists a homeomorphism $\Psi: U \to U'$, $\Psi (0,0)=(0,0)$, between two neighborhoods of the origin in $\C^{2}$, sending leaves of a representative $\F_{U}$ of the germ $\F$ over $U$ into leaves of a representative $\F'_{U'}$ of $\F'$ over $U'$. As we have already mentioned, we shall also suppose that $\Psi$ preserve the orientation of the ambient space as well as the orientation of the leaves. This assumption is important to have the implication:
$\Psi$ transversely conformal $\Rightarrow$ $\Psi$ transversely holomorphic.
We note that if $\F$ is given by a holomorphic differential form with real coefficients, then the homeomorphism $\Psi(x,y)=(\bar x,\bar y)$, conjugating $\F$ to itself, preserves the orientation of the ambient space but reverse the orientation of the leaves.\\
\end{comment}
{In fact,}
the equality $\mc S\mc L(\F)=\mc S\mc L(\F')$ need to be precised because the index sets of the families are different for $\F$ and $\F'$.
In order to do this, {we recall} that a topological {conjugation}  between $\F$ and $\F'$ as above transforms $S_{\F}$ into $S_{\F'}$ and induces a unique homeomorphism  $$\Psi ^\sharp : \mc E_\F\to\mc E_{\F'}\,,\quad \Psi ^\sharp(\Sing(\un\F))=\Sing(\un\F')$$ between the exceptional divisors up to isotopy. This is {a} consequence of the following result proved in a previous work \cite{MarMatMarq}.

\vspace{1em}

\noindent{\bf Marking Theorem.} \textit{Let $S$ and $S'$ be two germs of analytic curves at the origin in $\C^{2}$ and let  $h : (\mb C^2,0)\iso  (\mb C^2,0)$ be a germ of homeomorphism such that $h(S)=S'$. If $E_{S}$ and $E_{S'}$ denote the minimal reduction of singularities of $S$ and $S'$, then there is a germ of homeomorphism $h_1 : (\mb C^2,0)\iso  (\mb C^2,0)$
such that:
\begin{enumerate}[(i)]
  \item\label{condhomo}
  $h_1(S)=S'$ and the restrictions of $h$ and  $h_1$ to the complements of $S$ and $S'$ are homotopic,
  \item\label{conexcel}  $E_{S'}^{-1}\circ h_1\circ E_{S}$ extends to a homeomorphism
     from a neighborhood of $\mc D_S:=E_S^{-1}(S)$ onto a neighborhood of $E_{S'}^{-1}(S')$, being holomorphic over an open neighborhood of $\Sing(\mc D_{S})$ and compatible with the Hopf fibrations outside another open neighborhood of  $\mathrm{Sing}(\mc D_{S})$.
\end{enumerate}}

\begin{comment}
\medskip

\noindent The above mentioned  homotopic condition (i)  is made precise  in \cite{MarMatMarq}, and it is shown to be equivalent to the fact that the action of $h$ and $h_{1}$ on the fundamental group of the complement of $S$ and $S'$, in a small Milnor ball, coincide up to inner automorphisms, see Section~\ref{subsecmarquage}.
\end{comment}

\medskip

\emph{The topological invariance of $\mc S\mc L(\F)$ induced by the {conjugation} $\Psi $ between $\F$ and $\F'$} means that for each irreducible component $D$
of $\mc E_\F$ and each singular point   $s\in \Sing(\F)$, the following conditions hold:
\begin{enumerate}[a)]
  \item  $
\CS(\un\F', \Psi ^\sharp(D), \Psi ^{\sharp}(s))=\CS(\un\F, D,s))$ and  $[\un\F'_{\Psi ^\sharp(s)}]^{\mr{hol}}=[\un\F_s]^{\mr{hol}}$,
\item there exists  a germ of biholomorphism  $\psi $ between two germs $(\Delta ,m)$ and  $(\Delta ',m')$  of transverse analytic curves to $\F$ and $\F'$ through points $m\in D\setminus\Sing(\un\F)$ and $m':=\Psi ^\sharp(m)$ respectively, such that the following diagram is commutative:
\begin{equation}\label{diagdholonomie}
\begin{diagram}
 \pi _1(D\setminus\mr{Sing}(\un\F), \, m)& \rTo^{\phantom{aaa}\mc H^{\F}_D\phantom{aaa}} & \mathrm{Diff}(\Delta , m)\\
 \dTo{\Psi ^{\sharp}_\ast}& {\sstyle \circlearrowleft}& \dTo{\psi_{\ast}}\\
  \pi _1(D'\setminus\mr{Sing}(\un\F'), \, m')& \rTo^{\phantom{aaa}\mc H_{D'}^{\F'}\phantom{aaa}} & \mathrm{Diff}(\Delta ', m')\
\end{diagram}
\end{equation}
with $\psi _\ast(\varphi ):=\psi \circ\varphi \circ\psi ^{-1}$
 and $\Psi ^{\sharp}_\ast(\dot{\gamma }):=\Psi ^\sharp\circ\dot{\gamma }$.
\end{enumerate}

\medskip

Notice that  $\mc S\mc L(\F)$ is a ``semi-local'' invariant in the sense that it only contains information along the irreducible components of  $\mc E_\F$, but it does not provide any information about the combinatorial gluing of these data. Thus,  it cannot be reasonably a complete invariant of $\F$.
To remedy for this situation,
the idea is to consider the separatrix curve $S_{\F}$ as the ``organization center'' of the topology of $\F$ as was conjectured by Ren\'e Thom in the seventies. The incompressibility of the leaves inside the complement of $S_{\F}$ proved in \cite{MarMat} plays a major role here and it indicates that the fundamental group of the complement of $S_{\F}$ ``controls'' the topology of the leaves of $\F$. It also suggest the possibility to replace the usual notion of ``holonomy'' by that of ``monodromy''. The ``holonomy'' consists of the pseudo-group of local automorphisms of the ambient space coming from the ambiguity of the (multivalued) first integrals of $\F$. The ``monodromy'' reports the automorphisms of the set of (multivalued) first integrals coming from the ambiguity of the ambient space. We precise this notion in a general setting:

\medskip

\noindent\textbf{Definition.}\textit{
Let  $\mc G$ be a differentiable foliation on a manifold $M$ and consider the universal covering $q : \wt M\to M$ of $M$. We denote by
$\wt {\mc G}$ the lift of $\mc G$ in $\wt M$
and by  $\wt M/\wt {\mc G}$  the space of leaves of $\wt{\mc G}$.
Then the \emph{monodromy of $\mc G$} is the morphism
$$
\mf M^{\mc G}_{M} : \Aut(\wt M, q)\to \Aut(\wt M/\wt {\mc G})\,,
$$
sending an element  $\varphi $ of the group of deck transformations of $q$ to the automorphism of   $\wt M/\wt {\mc G}$  obtained by factorizing $\varphi $, i.e. $\mf M^{\mc G}_{M}(\varphi )\circ \tau =\tau \circ\varphi $, where $\tau:\wt M\to \wt M/\wt {\mc G}$ denotes the natural quotient map.
}\\

For instance, if $(M,\mc G)$ is a foliated bundle over $B$ with simply connected fibre $F$ then $\wt M=\wt B\times F$, $\wt M/\wt{\mc G}= F$ and the monodromy $\mf M_{M}^{\mc G}$ can be identified with the global holonomy representation $\pi_{1}(M)\to\mr{Aut}(F)$.\\

In \cite{MarMat} we have constructed a fundamental system  $(U_\alpha )_\alpha $ of neighborhoods of  $S_\F$ in a Milnor ball such that the space of leaves of the universal covering  $\wt U^\ast_\alpha \to U^\ast_\alpha $ of $U^\ast_\alpha := U_\alpha \setminus S_\F$ is a (in general non Hausdorff) holomorphic manifold. Thus, the monodromy of the global foliation
$\F_{|U_\alpha^\ast}$ is a representation of
$\pi _1(U_{\alpha }^\ast)$ to the group of holomorphic automorphisms of the leaf space {$\wt{\mc Q}_{U_{\alpha}}^{\F}$} of $\wt U^\ast_\alpha $. However,  this monodromy representation  depends on the choice of the open set where is defined the chosen representative of
 the germ $\F$.
 This notion of monodromy admits a reformulation in term of germs by considering the \emph{category of pro-objects}, cf. Section~\ref{subsec.germification}. This allows us to introduce the notion of \emph{monodromy of a germ of foliation} in Definition~\ref{mongermefeui}.
{Roughly speaking, the image of a deck transformation $g$ of the universal covering $\wt{ U}_{\alpha}^{*}\to U_{\alpha}^{*}$ by the monodromy representation of $\F$ is given by taking the germ over all the open sets $U_{\alpha}$ of the mappping from $\wt{\mc Q}_{U_{\alpha}}^{\F}$ to itself defined by $\wt{L}_{\alpha}\mapsto g(\wt L_{\alpha})$.}
The interest of this notion lies in the fact that it takes into account simultaneously the ``transverse structure'' of the foliation and the topology of the complement of their separatrix curves.

We {must} also consider the quite technical but highly rellevant
% the key
notions of \emph{geometric {conjugation} of monodromies} (Definition~\ref{conjugtopol})  {\emph{preserving the Camacho-Sad indices} (Definition~\ref{CS-preserving})
%
%and Remark~\ref{CS-obs})
%
and} \emph{realizable over transversals} (Definition~\ref{defconj}),
which allows us to compare {in a precise way} the monodromies of two germs of foliations.
{We refer the reader to Section~\ref{sec.monodromie} to have precise definitions of these technical notions. Using them, the statement of}
the main result of this paper is the following:\\

\it
\noindent{\bf Theorem II.} If $\F$ and $\F'$ are Generic General Type foliations, then the following properties are equivalent:
\begin{enumerate}
 \item\label{equloc} there exists a $\mc N$-topological {conjugation} between {the germs} $\F$ and $\F'$,
 \item\label{equecltrhol}
  there exists a $\mc N$-transversely holomorphic {conjugation} between {the foliations} $\un\F$ and $\un\F'$, {which is defined on open neighborhoods of the exceptional divisors}%inducing a  $\mc N$-excellent {conjugation} between $\F$ and $\F'$, see Definition~\ref{excellent},
 \item\label{conjmon} there exist a geometric $\mc N$-{conjugation} of the monodromies of the germs  $\F$ and $\F'$, preserving the  Camacho-Sad indices, which is realizable over $\mc N$-collections {of transversals} $\Sigma$ and $\Sigma'$ %of germs  of holomorphic regular curves transverse to the separatrix sets
 of $\F$ and $\F'$.
 \end{enumerate}
\rm

If $\mr{Node}(\un{\F})=\emptyset$ then \emph{a $\mc N$-collection of transversals of $\F$}  consists in
a single germ of regular holomorphic curve transverse to any separatrix of $\F$, at a regular point.
When $\mr{Node}(\un{\F})\neq\emptyset$ we must precise the location of the connected components of a $\mc N$-collection of transversals of $\F$ by using the following theorem of \cite[Corollary 4.1]{OBRGV10} generalizing the main result of \cite{CamachoSad}:\\

\noindent{\bf Strong Camacho-Sad Separatrix Theorem.} {\it Each connected component of $\mc E_{\F}\setminus(\mr{Node}(\un{\F})\cap\mr{Sing}(\mc E_{\F}))$
contains  a singular point of $\un{\F}$ lying on the strict transform of a  separatrix of $\F$ whose Camacho-Sad index has positive real part.}\\

{In particular, for each connected component $C$ of $\mc E_{\F}\setminus(\mr{Node}(\un{\F})\cap\mr{Sing}(\mc E_{\F}))$ exists a (not necessarily unique) non-nodal separatrix of $\F$ whose strict transform meets $C$. We define a \emph{$\mc N$-collection  of transversals of $\F$}
as a collection $\Sigma=\{(\Sigma_{1},p_{1}),\ldots,(\Sigma_{m},p_{m})\}$ where each $(\Sigma_{i},p_{i})$ is the germ at $p_{i}\in S_{\F}\setminus\{0\}$ of a regular curve $\Sigma_{i}$ transverse to $\F$ and the whole collection fulfils the following property:
\begin{enumerate}
\item[($\sigma $)] \it for each connected component $C$ of  $\mc E_{\F}\setminus(\mr{Node}(\un{\F})\cap\mr{Sing}(\mc E_{\F}))$ there exists a germ $(\Sigma_{i},p_{i})$ of $\Sigma$ with $p_{i}$ belonging to a non-nodal separatrix of $\F$ whose strict transform meets $C$.
\end{enumerate}
}

\medskip

In fact, Theorems I and II are easy consequences of Theorems~\ref{invariance} and~\ref{classification} that are going to be proved in Sections~\ref{7} and~\ref{8}. They are more general than Theorems I and II because  the genericity hypothesis (G) is replaced by weaker but more  technical conditions.
Notice also that Corollary~\ref{coroltrriggeneral} is of interest even in the case $\F=\F'$. If, in addition, we assume that $\mr{Node}(\un\F)=\emptyset$  then it implies the following result concerning the automorphism group
 $\Aut_{0}(\F)$ of orientations preserving homeomorphisms germs conjugating $\F$ to itself.\\

\noindent\textbf{Corollary.} \textit{If $\F$ is a Generic General Type foliation with $\mr{Node}(\un\F)=\emptyset$, then for each  $h\in \Aut_{0}(\F)$ there exists a homeomorphism $h_{1}\in \Aut_{0}(\F)$ satisfying for $S=S'=S_{\F}$
the properties  (i) and (ii) of Marking Theorem.}

\section{Preliminary notions}\label{notpreliminaires}
Through all the paper we will use the following notations:
$$\mb B_{r}=\{(x,y)\in\mb C^{2},\ |x|^{2}+|y|^{2}\le r\},\quad\mb D_{r}=\{z\in\mb C,\ |z|\le r\}\,,$$ and if {$B\subset A$, $B'\subset A'$,
$f:(A,B)\to(A',B')$} will denote the germ of a map $f_{1}$ defined {on a neighbourhood of $B$} in $A$ into $A'$, such that ${f_{1}(B)\subset B'}$.\\

In this section $S\subset \mb C^2$ denotes a  holomorphic curve with an isolated singularity at the origin $0=(0,0)$ of $\mb C^2$ and ${\mb B}:={\mb B}_{r_0}$ is
a closed \emph{Milnor ball} for $S$,  i.e.
each sphere $\partial {\mb B}_{r}$, $0<r\leq r_0$,  is transverse to $S$,  cf. \cite{Milnor}. We denote by $E_S : \mc B_S \to {\mb B}$ the \emph{minimal desingularization map} of $S$ such that the \emph{total divisor} $\mc D_S:= E_S^{-1}(S)$ has normal crossings. We denote by $\mc E_S:=E_S^{-1}(0)$ the \emph{exceptional divisor} and  by $\mc S:= \overline{\mc D_S\setminus \mc E_S}$ the \emph{strict transform} of $S$.
We will also use the following conventions along all the paper: for $A\subset\mb B$ and for $\mc A\subset\mc B_{S}$ we put
\begin{equation}\label{notast}
    A^\ast := A\setminus S\quad \hbox{ and }\quad \mc A^\ast := \mc A\setminus \mc D_S\,.
\end{equation}

\subsection{Incompressibility of the leaves}\label{subsec.bons.vois}
Let $\F$ be a singular non-dicritical holomorphic foliation defined in a neighborhood of $\mb B$, having $0$ as the unique singularity and $S$ as the \emph{separatrix curve in $\mb B$}, i.e. $S$ is invariant by $\F$ and every analytic invariant  curve passing through $0$ is contained in $S$. Denote by  $E_\F : \mc B_\F\to {\mb B}$ the minimal reduction of singularities map of $\F$ and by  $\underline{\F}= E_{\F}^*\F$ the reduced foliation over $\mc B_\F$.
The hypothesis (GT) on $\F$ implies that $\un{\F}$ does not have any saddle-node singularity, so that $\F$ is a \emph{generalized curve}
and consequently $\mc B_\F=\mc B_S$ and  $E_S=E_\F$, cf. \cite{CSL}.\\

Now we fix a  $\mc N$-collection of holomorphic transversal curves  $\Sigma\subset {\mb B}$ (i.e. satisfying Condition ($\sigma$) given in the Introduction). %It is easy to see that
{The saturation of a small transversal to one of the separatrices of a non-nodal (GT)-singularity $s$ jointly with the two local separatrices of $s$, is a neighborhood of $s$. Consequently the following property holds:}
\begin{enumerate}
  \item[\it a)]  \textit{there exist a $\mc N$-separator $N$ of $\F$, such that for any open neighborhood $W$ of $S$ in $\mb B$, the closure of $$\mr{Sat}_{\un{\F}}(E_{S}^{-1}({\Sigma}),E_{S}^{-1}(W)\setminus N)$$ is a
      neighborhood of  $\mc D_S\setminus \mr{Node}(\underline{\F})$, in $E_{S}^{-1}(W\setminus N)$,}
\end{enumerate}
where $\mathrm{Sat}_{\uF}(A,B)$ denotes the union of all the leaves of $\uF_{|B}$ meeting $A$ and it is called \emph{the saturation of $A$ in $B$ by $\uF$}.\\

We say that an open neighborhood $U$ of $S$ in $\mb B$ is  \emph{$(\F,\Sigma )$-admissible}, if for each leaf $L$ of the regular foliation $\F_{|U^\ast}$, the following properties hold:
\begin{enumerate}\it
  \item [b)] $L$ is \emph{incompressible} in $U^\ast$, i.e. the inclusion $L\subset U^\ast$ induces a mono\-morphism $\pi _1(L,p)\hookrightarrow \pi _1(U^\ast, p)$, $p\in L$; in addition the map $\pi _1(U^\ast, p)\to \pi _1({\mb B}^\ast, p)$ induced by the inclusion $U^\ast\subset {\mb B}^\ast$ is an isomorphism;
  \item [c)] every path in  $L$  whose ends
  are in  ${\Sigma}^\ast$ and which is homotopic in  $U^\ast$ to a path contained in $\Sigma^\ast$, is a  null homotopic loop in $L$.
\end{enumerate}
{Notice that an admissible open set $U\subset\mb B$ is not necessarily saturated in $\mb B$.}
We denote by
$\mf U_{\F,\Sigma }$ the collection of connected open neighborhoods of $S$ which are  $(\F,\Sigma )$-admissible.
Property {\it c)} above {will play a key role in this work. In fact, it} is equivalent to the \emph{foliated $1$-connexity} of $\Sigma \cap U ^\ast$ in $U^\ast$.
{This notion introduced in \cite{MarMat} plays a major role in  the proof of
the main result  of \cite[Theorem 6.1.1]{MarMat} which can be stated as follows:}
\begin{teo}\label{incfeuilles}
If $\F$ is a foliation of General Type, then $\mf U_{\F,\Sigma }$ is a fundamental system of neighborhoods of $S$ in the closed Milnor ball $\mb B$.
\end{teo}

\subsection{Leaf spaces}\label{subsec.esp.des.feuilles}
We fix once for all a universal covering  $q : \wt{{\mb B}}^\ast \to {\mb B}^\ast$ of ${\mb B}^\ast$ and for every subset $A\subset{\mb B}$ we will denote
\begin{equation}\label{notrevuniv}
\wt A^\ast := q
^{-1}(A^\ast)\quad\hbox{and}\quad  q_{\ssstyle A}:=q_{\ssstyle |\wt A^\ast} : \wt A^\ast\too A^\ast\,.
\end{equation}
If $U\in\mf U_{\F,\Sigma }$ then $q_{\ssstyle U}$ is a universal covering of  $U^*$. The group  $\Gamma
:= \aut_q
(\wt {\mb B}^\ast )$
of deck transformations of the covering $q$ can be identified with the group $\Gamma_{U}$ of deck transformations of the covering $q_{\ssstyle U}$ by the restriction map  $g\mapsto g_{|\wt U^\ast}$. Hence, we can also identify $\Gamma$ with
\begin{equation}\label{idgamma}
\Gamma _\infty:= \limproj_{U\in \mc U_{\F,\Sigma }}\Gamma _U\,.
\end{equation}

\indent
On $\wt{\mb B}^{\ast}$ we consider the regular foliation $\wt \F$, pull-back of $\F$ by $q$. For  $U\in \mf U_{\F,\Sigma }$,
we denote by  $\tF{U}$ its restriction to $\wt U^\ast$ and for  % a saturated
{an arbitrary}
subset $W$ of $U$, we denote by
\begin{equation}\label{espquot}
    \tQF{W}  := \left(\wt{W}^\ast\left/ \tF{W}\right)\right.\,,\quad \varrho _{\ssstyle  W} : \wt W^\ast \rightarrow \tQF{W}\subset \tQF{U}\,,
\end{equation}
the leaf space of the restriction $\tF{W}$ of $\tF{U}$ to $\wt W^\ast$, endowed with the quotient topology and the quotient map $\varrho_{\ssstyle W}$. {It turns out that if $W$ is saturated in $U$ then} the natural map $\tQF{W}\hookrightarrow \tQF{U}$ is a \emph{topological embedding}, i.e. a homeomorphism  onto its image, and we may  consider $\tQF{W}$ as a subset of $\tQF{U}$.

%
%On $\wt{\mb B}^{\ast}$ we consider the regular foliation $\wt \F$ pull-back of $\F$ by $q$. For  $U\in \mf U_{\F,\Sigma }$,
%we denote by  $\tF{U}$ its restriction to $\wt U^\ast$ and by
%$$\tQF{U}  := \left(\wt U^\ast\left/ \tF{U}\right)\right.\,,\quad \varrho _{\ssstyle  U} : \wt U^\ast \rightarrow \tQF{U}\,,$$
%the leaf space of $\tF U$ endowed with the quotient topology and the quotient map $\varrho_{\ssstyle U}$.

Properties  {\it a)}, {\it b)} and {\it c)} satisfied by the open sets of  $\mf U_{\F,\Sigma }$
can be understood as geometric properties of the foliation  $\wt\F_{U}$ with respect to the transverse section
$\wt{\Sigma}_U^*:=\wt{\Sigma}^\ast\cap\wt U^\ast$, cf. \cite[\S6.2]{MarMat}:
\begin{itemize}\it
  \item every leaf of  $\wt\F_{U}$ is simply connected;
  \item the intersection of every leaf of $\wt\F_U$ with each connected component of  $\wt{\Sigma}^\ast _U$
       is either empty or consists in a single point;
 \item  the  restriction of $\varrho _{\ssstyle U}$ to each connected component of  $\wt{\Sigma}_ U^{*}$, is a topological embedding.
\end{itemize}
Let $N$ be a $\mc N$-separator  of $\F$ satisfying Property \emph{a}) of the previous section.
The inverse maps of the $\varrho _{\ssstyle U}$'s restrictions to the connected components of  $\wt{\Sigma}_ U^{*}$ form a holomorphic atlas on $\tQF{U\setminus N}$,
defining a  structure of one dimensional complex manifold (non necessarily Hausdorff) over it.  It is easy to {check} that this  structure  extends to an unique structure of complex manifold over $\tQF{U}$, such that:
\begin{itemize}\it
\item for every holomorphic map $g :\mb D_{1}\rightarrow \wt U^\ast$, the composition
$$ \varrho _{\ssstyle U}\circ g  : \mb D_{1} \rightarrow
\tQF{U}$$ is also holomorphic.
\end{itemize}

It is clear that each element  $g$ of $\Gamma_U $ preserves the foliation  $\wt\F_U$ and factorizes by an element $g_U^{\flat}$ of  the analytic automorphism group $\mathrm{Aut}_{\mathrm{An}}(\tQF{U})$.
In  \cite{MarMat} we have defined the \emph{monodromy} of $\F_{U}$ as the morphism
\begin{equation}\label{monGlobale}
 \monF_U :   \Gamma _U \longrightarrow\mathrm{Aut}_{\mathrm{An}}(\tQF{U})\,,\qquad g\mapsto g_U^{\flat}\,.
\end{equation}
This representation of $\Gamma_{U}$ is clearly an {analytic} invariant of the  foliation $\F_{U}$ {by biholomorphisms preserving  the open set $U$.} In order to obtain an {analytic} invariant of the germ of $\F$ at $0$, or along $S$, we need to ``germify'' this notion. {This will be done in the following section.}

\section{Monodromy of a germ of foliation}\label{sec.monodromie}
\subsection{Germification}\label{subsec.germification}
The set $\mf
U_{\F, \Sigma }$  is cofiltered by the partial order $$U\preceq V :
\Longleftrightarrow U\supset V.$$
The maps
\begin{equation}\label{applrestr}
\rho_{\ssstyle  U V} :  \tFQ V\too \tFQ U\,,\quad V\subset U, \quad U,V\in \mc U_{\F,\,\Sigma }\;,
\end{equation}
sending each leaf $L$ of $\tF V$ into the unique leaf of $\tF U$ containing $L$, are open and holomorphic. They form a projective or inverse system of complex manifolds
   $$ \QF := \left(\left(\tFQ U\right)_{U\in \mf U_{\F,\Sigma }}, \;\left(\rho
_{\ssstyle  U V}\right)_{\ssstyle U, V \in \mf U_{\F,\Sigma},\, U \preceq V}\right)\,.$$
called the  \emph{leaf pro-space of  $\wt\F$}. It is an object in the category
$\proan$ of pro-objects associated to the category $\mathrm{An}$ of (non necessarily Hausdorff) complex manifolds and holomorphic maps.
We recall that the objets of $\proan$
are the projective families of complex manifolds; on the other hand, if $\mf A$ and $\mf B$ are cofiltered sets and
$$M = \left( (M_\alpha
)_{\alpha \in \mf A},\, (\zeta _{\alpha \alpha ' })_{\alpha \geq \alpha '}\right)\quad {\rm and}\quad
M' = \left( (M'_\beta )_{\beta \in \mf B},\, (\zeta' _{\beta \beta ' })_{\beta  \geq \beta' }\right)$$
are two objects of $\proan$ then
the set of
\emph{$\proan$-morphisms} of $M$ into $M'$
is by definition
\begin{equation}\label{promorphismes}
    \Hom_{\proan}( M, \, M')
    :=  \limproj_{\beta \in \mf B}\;\limind_{\alpha \in \mf A}
    \;\,\mc O
    (M_\alpha ,\, M'_\beta )\,,
\end{equation}
where $\mc O(M_\alpha ,\, M'_\beta )$ denotes the set of holomorphic maps of
$M_\alpha$ into $M'_\beta $.
For further details, see \cite{Douady}.

\subsection{Pro-germs at infinity} Let  $M$ be a complex submanifold of  $\wt{\mb B}^\ast$ such that  $\overline{q(M)}\cap S\neq\emptyset$. The projective system
$$(M, \infty ) := \left( (M\cap \wt U^\ast)_{U\in \mf U_{\F,\Sigma }} ,\,
(\iota_{UV})_{U,V\in \mf U_{\F, \Sigma },\,  V\subset U}\right)\,,
$$
formed by the inclusion maps $\iota_{UV}$ of $M\cap \wt V^\ast$ into $M\cap \wt U^\ast$ is a pro-object in $\proan$.
{Let $T$
be a complex manifold.}
Every element $g$ of the set
$${\mc O}((M,\infty),\,T) := \limind_{U\in \mf U_{\F,\Sigma }} \mc O(M\cap \wt U^\ast, \,
T)\,
$$
will be called \emph{germ at infinity} of $M$ into $T$ and denoted by $g : (M, \infty)\to T$. Identifying $T$ to the constant projective system, ${\mc O}((M,\infty),\,T) $ can be naturally identified with  $\Hom_\proan((M,\infty),\,T)$.

\begin{obs}\label{attention}
Endowing $M$ with the induced topology of  $\wt{\mb B} ^\ast$, the pre-sheaf
 $W\mapsto \mc O((W,\infty),\,T)$ is not a sheaf and two different germs at infinity
 $f,g\in \mc O((M,\infty), T)$  can coincide as elements of $\mc O((M\cap V_j,\infty)$, over the intersection of  $M$ with each open set of a covering $(V_j)_{j\in J}$.
In particular, if $M$ has infinitely many connected components  $M^{\alpha }$, each of them satisfying  $\overline{q(M^{\alpha })}\cap S\neq\emptyset$, then the restriction map $\mc O(M,\infty)\to\prod_\alpha \mc O(M^\alpha ,\infty)$ is never surjective.
\end{obs}
Assume now that $T$ is contained in  $\wt{\mb B}^\ast$.
Each element $f$ of the set
$$%\mc O((M,\infty), \,(T,\infty)):=
\Hom_{\proan}((M,\infty), (T,\infty) )$$
is called \emph{pro-germ at infinity} of $M$ into $T$ and it shall be denoted by  $f :(M,\infty)\to (T,\infty)$. Thus, $f$ is a family of germs at infinity $$f=(f_V)_{V\in \mc U_{\ssstyle \F,\,\Sigma }}\in\prod_{V\in\mf U_{\F,\Sigma }}\mc O((M,\infty),\,T\cap \wt V^{*})\,,$$
such that $\varsigma_{VW}\circ f_W=f_V$,  $W\subset V$,  where $\varsigma_{WV}$ denotes the inclusion map of  $T\cap \wt W^{*}$ into $T\cap \wt V^{* }$.
The same notions in the category
$\mathrm{Top}$ of topological spaces and continuous maps define the set of continuous pro-germs at infinity
 $$%\mc \mc C^0((M,\infty), \,(T,\infty)):=
 \Hom_{\underleftarrow{\mathrm{Top}}}((M,\infty), (T,\infty) )\,.$$

Notice that the group of pro-germs at infinity of deck transformations of the
covering can be canonically identified with the group  $\Gamma _\infty$ defined by (\ref{idgamma}):
$$\Gamma _\infty \simeq \left\{\varphi \in \Aut_{\proan}(\wt {\mb B} ^\ast,\infty)
\;\left|\; q_{\ssstyle \infty }\circ \varphi = q_{\ssstyle \infty}\right.\right\}\,,
$$
where
$q_{\ssstyle \infty} : (\wt{\mb B}^\ast ,\infty)\to {\mb B}$ denotes the germ at infinity of the covering map $q$.

\subsection{Canonical pro-germs}\label{progermcan} {For every $U\in \mf U_{\F, \Sigma}$}
we denote by  $$\tau _{\ssstyle M,U} : (M,\infty)\to \tQF{U}$$ the germ at infinity of the quotient map $\wt U\to \tQF{U}$ restricted to {$M\cap\wt U^{*}$}. The element
$$\tau_M :=(\tau_{\ssstyle M,U})_{U\in \mf U_{\F, \,\Sigma }}\in \Hom_\proan((M,\infty),\, \QF)\subset \prod_{U\in\mf U_{\F,\Sigma }}\mc O((M,\infty),\,\tQF{U})
$$
will be called the \emph{canonical pro-morphism} associated to $M$.
Next proposition follows easily from the geometric properties of the foliation $\wt\F_{U}$ relatively to $\wt \Sigma {}^\ast_U$ stated in Section~\ref{subsec.bons.vois}.
\begin{prop}\label{monoprogermcan}
If $M$ is a connected component of  $\wt\Sigma^\ast $ then $\tau_{\ssstyle M}$ is a monomorphism in the category $\proan$.
\end{prop}

\subsection{Monodromy of a germ}
Let $g$ be an element of $\Gamma_{\infty}$ and consider $U,V\in\mc U_{\F,\,\Sigma }$ with $V\subset U$. With the notations  (\ref{monGlobale}) and  (\ref{applrestr}) we have the following commutation relations:
$$
g^\flat_U\circ \rho _{\ssstyle UV}= \rho _{\ssstyle UV}\circ g_V^\flat \in \mc O(\tFQ{V}\,,\; \tFQ{U})\,.
$$
Hence, by denoting  $\mc O(\tFQ{\infty},\,\tFQ{U}):= \limind_{V\in \mf U_{\F,\,\Sigma }}\mc O(\tFQ{V},\,\tFQ{U})$,
the elements  $$ g_{U\infty}^\flat:=\limind_V (g^\flat_U\circ \rho _{\ssstyle UV})\in \mc O(\tQF{\infty},\tQF{U})$$
form a projective family. The  $\proan$-endomorphism
 $$  g^\flat := (g^\flat_{U\infty})_{U\in \mf U_{\F,\,\Sigma }} \in  \mathrm{End}_{\proan}(\QF)\subset \prod_{U\in \mf U_{\F,\,\Sigma }}\mc O(\tFQ{\infty},\,\tFQ{U})\,,
$$
is invertible and its inverse is  $(g^{-1})^\flat$. More generally, we have the following covariance relations:
$$(g\circ h)^\flat = g^\flat\circ h^\flat\,,\quad g,\; h\in \Gamma_\infty\,.$$
\begin{defin}\label{mongermefeui}
The morphism of groups
$$\monFS : \Gamma_\infty \;\longrightarrow\;\Aut_{\proan} ( \QF)\,,\qquad
g \mapsto \monFS(g):= g^\flat \,.$$
is called the \emph{monodromy of the germ of $\F$ along $S$.}
\end{defin}

We now fix once for all a second curve $S'\subset \mb C^2$ with an isolated singularity at the origin  as well as a closed Milnor ball ${\mb B}'$ for $S'$.
We denote by  $E_{S'} : \mc B_{S'}'\to {\mb B}'$, $\mc D_{S'}$, $\mc E_{S'}$, $\mc S'$,
the minimal desingularization map, the total transform, the exceptional divisor and the strict transform of $S'$
 respectively. In order to avoid any ambiguity with the notation (\ref{notast}), we denote
\begin{equation}\label{notationstar}
A^\star := A\setminus S'\,,\quad \mc A^{\star} := \mc A\setminus \mc D_{S'}\,,\quad \hbox{\rm for}\quad A\subset {\mb B}'\quad \hbox{\rm and}\quad \mc A\subset \mc B_{S'}'\,.
\end{equation}
We also fix a singular non-dicritical holomorphic foliation $\F'$ of General Type defined in a neighborhood of  ${\mb B}'$ having $S'$ as total separatrix curve.
Let $\Sigma '\subset {\mb B}'$ be a $\mc N$-collection of transversals of $\F'$, we denote by  $\mf U_{\F', \Sigma '}$ the set of open neighborhoods of $S'$ in $\mb B'$ which are  $(\F',\Sigma ')$-admissible
and we fix a universal covering $q' : \wt{\mb B} '{}^\star\to {\mb B} '{}^\star$.
For  $A\subset {\mb B} '$, we put
\begin{equation}\label{notationqpr}
\wt A^\star :=q'^{-1}(A^\star)\,,\quad \hbox{and}\quad q'_{A}:=q'_{|\wt A^\star} : \wt A^\star \to A^\star\,.
\end{equation}
As in (\ref{idgamma}), we identify the projective limit of the groups $\Gamma'_{U}$ of deck transformations of the covering $q'_{\ssstyle U'}$ with the group of pro-automorphisms at infinity of $\mb B'$ preserving the germ at infinity $q'_\infty $ of $q'$:
$$
\Gamma '_{\infty}:=
\limproj_{U\in \mf U_{\F',\Sigma '}}\Gamma '_U\simeq
\left\{\varphi\in \proaut({\mb B} ',\infty)
\;\left|\; q'_{\ssstyle \infty }\circ \varphi = q'_{\ssstyle \infty}\right.\right\}\,.
$$
When $\mathrm{Node}(\uF)\neq \emptyset$, we also need to consider a mixed class of pro-germs, called \emph{$\mc N$-analytic pro-germs}:
\begin{defin}\label{NAnprogerms} A  germ
$f : (\mb B,S)\to (\mb B',S')$ will be called \emph{$\mc N$-analytic} (with respect to $\F$ and $\F'$) if it can be represented by a continuous map $\underline{f} : U\to U' $, $U\in \mf U_{\F,\Sigma }$, $U'\in \mf U_{\F',\Sigma '}$ for which there are $\mc N$-separators $N$ and $N'$ of $\F$ and $\F'$, such that
$\underline{f}$ is holomorphic in $U\setminus N$ and $\underline{f}(U\cap N)\subset U'\cap N'$.
An element of
$$\Hom_{\mc N-\proan}(\wt{\mc Q}^{\F}_{\infty},\wt{\mc Q}^{\F'}_{\infty})\subset\Hom_{\protop}(\wt{\mc Q}^{\F}_{\infty},\wt{\mc Q}^{\F'}_{\infty})$$ consists in a collection $f=(f_{U'})_U'$ of germs which  can be represented by continuous maps $\underline{f}_{U'} : \wt{\mc Q}^{\F}_{U}\to \wt{\mc Q}^{\F'}_{U'} $, $U\in \mf U_{\F,\Sigma }$, $U'\in \mf U_{\F',\Sigma '}$, for which there are $\mc N$-separators $N$ and $N'$ of $\F$ and $\F'$, such that
$\underline{f}_{U'}(\wt{\mc Q}^{\F}_{U\cap N})\subset \wt{\mc Q}^{\F}_{U'\cap N'}$ and $\underline{f}_{U'}$ is holomorphic in $\wt{\mc Q}^{\F}_{U\setminus N}$, cf. (\ref{espquot}).
\end{defin}
%
%\begin{defin}\label{NAnprogerms}
%An element of
%$$\Hom_{\mc N-\proan}(\mb B,\infty), (\mb B',\infty) )\subset \Hom_{\protop}((\mb B,\infty), \,(\mb B',\infty))$$
%consists in a collection $f=(f_V)_V$ of germs which  can be represented by continuous maps $\underline{f}_V : U\to V $, $U\in \mf U_{\F,\Sigma }$, $V\in \mf U_{\F',\Sigma '}$, for which there are $\mc N$-separators $N$ and $N'$ of $\F$ and $\F'$ such that $\underline{f}_{V}(U\cap N)\subset V\cap N'$ and $\underline{f}_{V}$ is holomorphic in $U\setminus N$. An element of
%$$\Hom_{\mc N-\proan}(\wt{\mc Q}^{\F}_{\infty},\wt{\mc Q}^{\F'}_{\infty})\subset\Hom_{\protop}(\wt{\mc Q}^{\F}_{\infty},\wt{\mc Q}^{\F'}_{\infty})$$ consists in a collection $f=(f_V)_V$ of germs which  can be represented by continuous maps $\underline{f}_V : \wt{\mc Q}^{\F}_{U}\to \wt{\mc Q}^{\F'}_{V} $, $U\in \mf U_{\F,\Sigma }$, $V\in \mf U_{\F',\Sigma '}$, for which there are $\mc N$-separators $N$ and $N'$ of $\F$ and $\F'$ such that $\underline{f}_{V}(\wt{\mc Q}^{\F}_{U\cap N})\subset \wt{\mc Q}^{\F}_{V\cap N'}$ and $\underline{f}_{V}$ is holomorphic in $\wt{\mc Q}^{\F}_{U\setminus N}$, \JF{cf. Section~\ref{subsec.bons.vois}}.
%\end{defin}

{The usual notion of conjugation of group representations induce the notion of  \emph{conjugation} (resp. $\mc N$-\emph{conjugation}) between the mondromies $\monFS$ and $\monFSpr$ as a pair $(\mf g,h)$ where $\mf g: \Gamma _\infty \iso \Gamma '_{\infty}$ is an isomorphism of groups and $h\in\mathrm{Isom}_{\proan}(\QF,\,\QFpr)$ (resp. $h\in\mathrm{Isom}_{\mc N-\proan}(\QF,\,\QFpr)$) satisfying the commutation relation $h_{*}\circ \monFS=\monFSpr\circ \mf g$ where
\begin{eqnarray}
h_{*}\ :\ \Aut_{\proan}
( \QF)&\iso& \Aut_{\proan} ( \QFpr)\label{hhast}\\
 \varphi &\mapsto & h\circ \varphi\circ h^{-1}\, .\nonumber
\end{eqnarray}
This notion of conjugation is algebraic and does not take in account the many  topological informations contained in the monodromy morphisms. Then we introduce the following more specific notion:}

\begin{comment}
We say that the monodromies  $\monFS$ and $\monFSpr$ are   \emph{algebraically equivalent} if there are group isomorphisms  $$\mf g: \Gamma _\infty \iso \Gamma '_{\infty}\quad \hbox{ and }\quad  {h} : \Aut_{\proan}
( \QF)\iso \Aut_{\proan} ( \QFpr)$$ satisfying the  commutation relation ${h}\circ \monFS=\monFSpr\circ \mf g$.
%
The pair   $(\mf g,{h})$ will be called  \emph{a {conjugation}} (resp. \emph{a $\mc N$-{conjugation}})  \emph{between  $\monFS$ and $\monFSpr$}, if ${h}$ is the conjugation morphism associated to an  an element  $h\in\mathrm{Isom}_{\proan}(\QF,\,\QFpr)$ (resp. $h\in\mathrm{Isom}_{\mc N-\proan}(\QF,\,\QFpr)$):
\begin{equation}\label{hhast}
{h}=h_{\ast}: \varphi\mapsto h\circ \varphi\circ h^{-1}\,.
\end{equation}
%If, in addition, there exists an element  $h\in\mathrm{Isom}_{\proan}(\QF,\,\QFpr)$ (resp. $h\in\mathrm{Isom}_{\mc N-\proan}(\QF,\,\QFpr)$), such that ${h}$ is the \emph{conjugation morphism} $h_\ast$,
%\begin{equation}\label{hhast}
%{h}=h_{\ast}: \varphi\mapsto h\circ \varphi\circ h^{-1}\,,
%\end{equation}
%we will say that the pair  $(\mf g,{h})$ is an \emph{algebraic {conjugation}} (resp. algebraic $\mc N$-{conjugation}) between  $\monFS$ and $\monFSpr$.
%
\end{comment}
\begin{defin}\label{conjugtopol}
A \emph{geometric {conjugation}} (resp. \emph{geometric $\mc N$-{conjugation}}) between the monodromies  $\monFS$ and $\monFSpr$, is a {conjugation} (resp. $\mc N$-{conju\-ga\-tion}), $(\mf g, {h})$ such that there are a  homeomorphism germ $g$ {not necessarily foliated)} from $({\mb B} ,S)$ into $({\mb B}',S')$
preserving the orientations of $\mb B$, $\mb B'$ and $S$, $S'$ and  a pro-germ at infinity  $\wt g$ from $(\wt{\mb B} ^\ast, \infty)$ into $(\wt {\mb B} '{}^\star,\infty)$ lifting  $g$, i.e. $q'_\infty\circ \wt g = g \circ q_\infty$, such that
$\mf g$ equals the conjugation morphism defined by  $\wt g$:
\begin{equation}\label{defconjaut}
\mf g = \wt g _{\ast } :\Gamma _\infty\to\Gamma '_\infty\,,\quad
 \varphi\mapsto
{\wt g}\circ \varphi\circ{\wt g}^{-1}\,.
\end{equation}
\end{defin}
\noindent We {then} have the following commutative diagram:
$$\begin{array}{ccc}
\Gamma _\infty&\hspace{-2em} \stackrel{\monFS} {\longrightarrow}&
\hspace{-0,5em}\mathrm{Aut}_{\proan}(\QF)
\\
\scriptstyle\mf g= \wt{ g}_{*}\; \downarrow\phantom{\wt\scriptstyle\wt{g}_{*}\,=\,\mf f\;}
&\hspace{-1em}\circlearrowright &\hspace{-2em}
\phantom{{\scriptstyle {h_{*}}}}\downarrow {\scriptstyle {h_{*}}}\\
\Gamma' _\infty& \hspace{-2em}\stackrel{\monFSpr}{\longrightarrow}&
\hspace{-0,5em}\mathrm{Aut}_{\proan}(\QFpr)
\end{array}\,,
$$
\noindent and we will say that the triple  $(g, \wt g, h )$  \emph{represents geometrically} the\linebreak ($\mc N$-){conjugation}  $(\mf g, {h})$.
%
\begin{comment}
Notice that if $(\mf f,\mf  h)$ is a geometric ($\mc N$-){conjugation} between  $\monFS$ and $\monFSpr$, then for each $\varphi\in \Gamma _\infty$ and for each  $\varphi '\in\Gamma '_\infty$, the pairs  $$(\mf f\circ \varsigma _{\varphi},\, {h}\circ \monFS(\varphi ) )\quad\hbox{ and }\quad (\varsigma _{\varphi'} \circ\mf f,\, \monFSpr (\varphi ')\circ{h})$$
are also geometric ($\mc N$-){conjugations} of these monodromies,
where
$\varsigma _\varphi $ and $\varsigma _{\varphi '}$ denote the inner automorphisms of  $\Gamma _\infty$ and $\Gamma '_\infty$ defined by  $\varphi $ and $\varphi '$, respectively.
\end{comment}
\begin{obs}\label{invparisotopie}
Let $\Theta _t : U\iso\Theta _t(U)$, $U\subset{\mb B}$, be a \emph{$S$-isotopy}, i.e. a continuous family of homeomorphisms depending on a parameter $t\in[0,1]$ such that  for each $t\in[0,1]$, $\Theta_{t}(U)$ is an open neighborhood of $S$ in $\mb B$,  $\Theta_{t}(S)=S$, $\Theta_{t}(p)=p$ for all $p\in U\cap\partial\mb B$ and $\Theta_{0}=\mr{id}_{U}$. We denote by $\wt{\Theta }_t$ the open embedding of  $\wt U^\ast$ into $\wt{{\mb B} }^\ast$ lifting  $\Theta _t$,
such that $\wt\Theta_t (\wt p)=\wt p$ if $q(\wt p)\in \partial {\mb B}$. It depends continuously on $t$, and    $\wt\Theta _0=id_{\wt U^{\ast}}$. {Since the induced maps  $(\wt g\circ \wt\Theta _t)_\ast : \Gamma _\infty\to \Gamma'_{\infty}$ depend continuously on $t$, they are constant and consequently} if $(g, \wt g, h)$ is a geometric representation of a ($\mc N$-){conjugation}  $(\mf g, {h})$, then  $(g\circ\Theta _1, \wt g\circ\wt\Theta _1, h)$ is also a geometric representation of $(\mf g, {h})$.  In the same way, if   $\Theta '_t: U'\iso \Theta '_t(U')$, $U'\in{\mb B} '$, $t\in [0,1]$, is a  $S'$-isotopy, then  $(\Theta '_1\circ g, \wt\Theta '_1\circ\wt g, h)$ is also a geometric representation of $(\mf g,{h})$.
\end{obs}

\subsection{Marking a germ of curve}
\label{subsecmarquage}
In  \cite{MarMatMarq} we have introduced the notion of marking of a germ of curve $(S',0)$ by another germ of curve $(S,0)$ as
a fundamental equivalence class of germs of homeomorphisms of $(\C^{2},0)$ conjugating $(S,0)$ to $(S',0)$.  From Proposition 2.8 of \cite{MarMatMarq} follows that two homeomorphic germs $\phi_{0}$ and $\phi_{1}$ are fundamentally equivalent if and only if
one of the following equivalent properties are satisfied:
\begin{enumerate}
\item \it there exists  $\varepsilon>0$ and a homotopy $\Phi\in \mc C^0({\mb B}^\ast_{\varepsilon}\times[0,1], {\mb B}'{}^\star)$ such that  $\Phi(\cdot, 0)$ and $\Phi(\cdot, 1)$ are representatives of the
 restrictions to the complement of $S$ of the germs $\phi_{0}$ and $\phi_{0}$ respectively;
\rm  \item \it there exist representatives
$\underline{\phi}_0$, $\underline{\phi}_1$  of the germs $\phi_0$ and $\phi_1$ on a small ball  ${\mb B}_{\varepsilon}$ such that for all $p\in {\mb B}_{\varepsilon}^\ast$ there is a path $\alpha $ contained in ${\mb B}'{}^\star$ with endpoints
 $\underline{\phi}_0(p)$ and $\underline{\phi}_1(p)$ such that the morphism  $\alpha _\ast : \dot{\gamma }\mapsto \dot{\alpha}^{-1} \svee\dot{\gamma} \svee \dot{\alpha}$ makes commutative the following diagram:

\hspace{5mm}
\xymatrix{   \pi _1({\mb B}_{\varepsilon}^\ast,\,p) \ar[rr]^{\underline{\phi}_{0\ast}}\ar[rdd]_{\underline{\phi}_{1\ast}} & & \pi _1({\mb B}'{}^\star, \underline{\phi}_0(p))\ar[ldd]^{\alpha_{*}}\\
& & \\
& \pi _1({\mb B}'{}^\star, \underline{\phi}_1(p)) &}
\end{enumerate}
\begin{obs}\label{3.5.1}
The morphism  $\mf g$ of a geometric ($\mc N$-){conjugation}  $(\mf g, {h})$ between the monodromies $\monFS$ and $\monFSpr$ determines a marking of $(S',0)$ by $(S,0)$, that we will denote again by $\mf g$.
\end{obs}

In \cite{MarMatMarq} we have seen that every marking can be represented by a homeomorphism having good regularity properties.
In order to precise these properties we need to consider some auxiliary geometric data. For each irreducible components $D$ of $\mc D_{S}$ and $D'$ of $\mc D_{S'}$ we fix germs of submersions
$$\pi _D : (\mc B_{S},D)\to D\quad\hbox{and}\quad \pi _{D'} : (\mc B_{S'},D)\to D'
$$
whose respective restrictions to $D$ and $D'$ are the identity. They will be called the \emph{Hopf fibrations of $D$ and $D'$}.
For each singular point   $s\in \Sing(\mc D_{S})$ and  $s'\in\Sing(\mc D_{S'})$, we fix holomorphic charts  $(x_s,y_s) : \mc W_s\iso \mb D_{1}^2$,  $(x_{s'},y_{s'}):\mc W_{s'}\iso \mb D_{1}^2$, with disjoint domains and such that the local equations of $\mc D_S$ and  $\mc D_{S'}$ are monomial in these coordinates.
The collections \begin{equation}\label{systloc}
\mc L:=((\pi _D)_{D},\,(x_s,y_s)_{s})\,,\quad \mc L':=((\pi _{D'})_{D'},\,(x_{s'},y_{s'})_{s'})
\end{equation}
will be called \emph{local data} for $S$ and $S'$.

\begin{defin}\label{excellent}
A germ of homeomorphism $g$ from $({\mb B},S)$ into $({\mb B}',S')$ will be called  \emph{excellent} (resp. \emph{$\mc N$-excellent}) with respect to the local data $\mc L$ and $\mc L'$ (resp. and a foliation $\F$ with separatrix set $S_{\F}=S$), if it admits a homeomorphic lifting $G$ from  $(\mc B_S,\mc D_S)$ into $(\mc B'_{S'}, \mc D_{S'})$ satisfying the following conditions:
\begin{enumerate}
  \item\label{consdiv} $G(\mc D_S)=\mc D_{S'}$ and $G(\mc D_S\cap \mc W_s)=\mc D_{S'}\cap \mc W'_{G(s)}$,  $s\in \Sing(\mc D_S)$,
  \item\label{Gholo} $G$ is holomorphic in a neighborhood of
 $\Sing(\mc D_{S})$\\ (resp. $\Sing(\mc D_{S})\setminus(\mr{Node}(\uF)\cap\mr{Sing}(\mc E_{S}))$),
  \item\label{comhopf}  the restriction of $G$ to a neighborhood of {the adherence of}
  $$\mc D_S\setminus \bigcup_{s\in\Sing(\mc D_S)}{\mc W}_s$$ commute with the Hopf fibrations,
 i.e. $\pi _{G(D)}\circ G=G\circ \pi _D$,
\end{enumerate}
\end{defin}

\noindent Once we fix the local data  $\mc L$ and $\mc L'$ we can precise the Marking Theorem stated in the introduction, in the following way.

\begin{teo}\label{marquageprecis} \cite{MarMatMarq} Every marking of  $S'$ by $S$ possesses an excellent representative with respect the local data  $\mc L$ and $\mc L'$.
\end{teo}

\begin{cor}\label{reprfibree}
Every geometric ($\mc N$-){conjugation}  between the monodromies $\monFS$ and $\monFSpr$ can be geometrically represented by a triple $(g, \wt g, h)$, where  $g$ is excellent  with respect to the local data  $\mc L$ and $\mc L'$.
\end{cor}

\begin{defin}\label{CS-preserving}
{We say that a geometric ($\mc N$-)conjugation $(\mf g,{h})$ between the monodromies of two foliation germs $\F$ and $\F'$ \emph{preserves the Camacho-Sad indices} if
once we  represent geometrically it by  a triple $(g,\wt g,h)$ with $g$ an excellent homeomorphism germ, then its lifting  $G:(\mc B_{S},\mc D_{S})\to(\mc B_{S'}',\mc D_{S'})$ satisfies
$\CS(\un\F,D,s)=\CS(\un\F',G(D),G(s))$
for every irreducible component $D$ of $\mc D_{S}$ and every singular point $s\in D$ of $\un\F$.}
\end{defin}

%\begin{obs}\label{CS-obs}
{We will see in Section~\ref{invtousCSpreuve} that a ($\mc N$-)conjugation $(\mf g,h)$ between the monodromies of $\F$ and $\F'$ preserves the Camacho-Sad indices if and only if there is a representative $g$ of the marking determined by $\mf g$ (cf. Remark~\ref{3.5.1}) such that for each irreducible component $\br S$ of $S$ we have the equality
$\CS(\un\F,\br{\mc S},\br{s})=\CS(\un\F',\br{\mc S'},\br{s'})$ where $\br{\mc S}$ and $\br{\mc S'}$ denote respectively the strict transforms of $\br S$ and $\br S'$ and $\br s$ and $\br s'$ are their corresponding attaching points in the exceptional divisors.}
%\end{obs}

\subsection{Realizations of {conjugations}}\label{subsectrealconj}
Consider subsets $V\subset {\mb B} $ and $V'\subset{\mb B}'$ such that
 $\overline{V^\ast}\cap S$ and $\overline{V'{}^\star}\cap S'$ are non-empty. We denote by
$\Gamma _{\wt V^\ast,\infty}$
(resp. $\Gamma '_{\wt V'{}^\star,\,\infty}$), the group of germs at infinity $\varphi$ of deck transformations of the (possibly non-connected) covering
 $(\wt V^\ast,\infty)$ (resp.  $(\wt V'{}^\star,\infty)$), i.e.
satisfying $q_\infty\circ \varphi=q_\infty$ (resp. $q'\circ \varphi=q'_\infty$). Clearly, the restriction maps define group monomorphisms $\iota :\Gamma _\infty\hookrightarrow\Gamma _{\wt V^\ast,\,\infty}$ and $\iota' :\Gamma '_\infty\hookrightarrow\Gamma '_{\wt V'{}^\star,\,\infty}$.

\begin{defin}\label{defconj}
A {geometric conjugation} (resp. geometric $\mc N$-{conjugation}),  $(\mf g,{h})$ between the monodromies  $\monFS$ and $\monFSpr$ is called \emph{realizable over the germs $(V,S)$ and $(V',S')$} if there exists
%a triple $(\psi,\wt\psi, h)$ constituted by an $\proan$-isomorphism (resp. $\mc N$-$\proan$-isomorphism), $h$ from $\QF$ onto $\QFpr$  such that  ${h}=h_\ast$,  cf. (\ref{hhast}),
a germ of homeomorphism  $\psi$ from $(V,S)$ into $(V',S')$, and a continuous pro-germ at infinity $$\wt\psi\in \Hom_{\underleftarrow{\mathrm{Top}}}((\wt V^\ast,\infty),\,(\wt V'{}^\star,\infty))$$ lifting  $\psi$, i.e. $q'_{\infty}\circ \wt\psi=\psi\circ q_{\ssstyle\infty}$, such that the following diagrams commute:
$${%\sstyle
(\star)}%_{\wt{\psi}}
\hspace{1em}
 \begin{array}{ccc}
(\wt V^\ast, \infty) & \stackrel{\wt \psi}{\longrightarrow} & (\wt V'{}^\star, \infty)\\
{\sstyle\tau _{\ssstyle \wt V^\ast}}\downarrow\phantom {\tau
_{\ssstyle \wt V^ast}} & {\sstyle\circlearrowright}&\phantom{\scriptstyle\tau _{\ssstyle \wt V'{}^\star}}\downarrow
{\scriptstyle\tau _{\ssstyle \wt V'{}^\star}} \\
\QF & \stackrel{h}{\longrightarrow} & \QFpr
\end{array}\,,\qquad
{%\sstyle
(\star\star)}%_{\wt{\psi}}
\hspace{1em}
\begin{array}{ccc}
\Gamma_\infty & \stackrel{\iota}{\hookrightarrow} & \Gamma _{\wt V^\ast,\,\infty }\\
{\scriptstyle\mf g}\downarrow\phantom{\scriptstyle\phi_{*}} &
\hspace{-0,2cm}{\sstyle\circlearrowright}
 & \phantom{\scriptscriptstyle \wt\psi_{\ast}}\downarrow {\scriptscriptstyle \wt\psi_{*}}\\
\Gamma_\infty' &
\hspace{-0,2cm}
\stackrel{\iota'}{\hookrightarrow} & \Gamma'_{\wt V'{}^\star, \,\infty }\end{array}$$
where
$\tau _{\ssstyle\wt V^\ast}$ and $\tau _{\ssstyle\wt V'{}^\star}$ are the canonical pro-germs defined in  Section~\ref{progermcan} and  $\wt \psi_*$ is the {conjugation} morphism  $\varphi\mapsto
\wt\psi\circ \varphi\circ \wt\psi^{-1}$.
We will say then that the triple
$(\psi,\wt\psi, h)$ is a \emph{realization} of $(\mf g,{h})$ over the germs  of $V$ and $V'$.
\end{defin}
\begin{obs}\label{conjfeuilconjmon} If $g : ({\mb B} ,S)\to ({\mb B} ',S')$ is a ($\mc N$-)transversely holomorphic {conjugation} between
  $\F$ and $\F'$, then every lifting $\wt g : (\wt
{\mb B}^\ast ,\infty) \rightarrow (\wt {\mb B} '{}^\star, \infty)$ of $g$ determines a ($\mc N$-)$\proan$-isomorphism $h$ from $\QF$ onto $\QFpr$. The pair  $(\wt g_\ast, {h})$  constituted by the conjugation isomorphism  (\ref{defconjaut}) % and  (\ref{hhast})
is a geometric ($\mc N$-){conjugation} between the monodromies of these germs of foliations and $(g, \wt g, h)$ is a geometric realization of this {conjugation} {over $(\mb B,S)$ and $(\mb B',S')$.}
\end{obs}
\begin{obs}\label{conjgeom}
If $W$ is a submanifold of $V$ such that $\overline{W^\ast}\cap S\neq\emptyset$, then the  restriction $(\psi _{|W}, \wt\psi _{|\wt{W}^\ast}, h)$ of a realization $(\psi , \wt \psi , h)$ of $(\mf g,{h})$ over $V$ and $V'$, is a realization of $(\mf g,{h})$ over $W$ and $\psi (W)$.
\end{obs}

A $S$-isotopy $\Theta _t : U\to \Theta _t(U)$, $U\subset {\mb B}$, in the sense of  Remark~\ref{invparisotopie} will be called \emph{$\F$-isotopy}, if   for all $p\in U$ the path
 $[0,1]\ni t\mapsto \Theta _t(p)$ is contained in a leaf of $\F$.
We have the following invariance property:
\begin{prop}\label{lemmeisotopie}
Let $\Theta _t : U\to \Theta _t(U)$
(resp. $\Theta '_t : U'\to \Theta '_t(U')$),  $t\in[0,1]$, be a $\F$-isotopy (resp. $\F'$-isotopy), defined in an open neighborhood  $U\supset S$ of ${\mb B} $ (resp. $U'\supset S'$ of ${\mb B} '$) and let   $(g, \wt g, h)$ be a realization of a {geometric} ($\mc N$-){conjugation} $(\mf g, {h})$ between the  monodromies of $\F$ and $\F'$, over subsets $V$ and $V'$ of ${\mb B} ^\ast$ and ${\mb B} '{}^\star$. Then $(\Theta '_1\circ g\circ \Theta _1^{-1}, \wt\Theta '_1\circ\wt g\circ \wt\Theta _1^{-1}, h)$ is a realization of  $(\mf g, {h})$ over $\Theta _1(V)$ and $\Theta '^{-1}_1(V')$, where the lifting  $\wt \Theta _t$ and $\wt \Theta '_t$ are defined as in Remark~\ref{invparisotopie}.
\end{prop}
\begin{dem}
{The idea  is  quite simple: on one hand}  $\tilde\Theta_{t}$ preserves the leaves of $\wt\F$ and consequently the {corresponding} diagram $(\star)$
%_{\tilde\Theta_{1}'\circ\tilde g\circ\tilde \Theta_{1}^{-1}}$
is commutative; on the other hand, Remark~\ref{invparisotopie} implies the equality
$(\tilde\Theta_{1}'\circ\tilde g\circ\tilde \Theta_{1}^{-1})_{*}=\tilde g_{*}$ and therefore the diagram  $(\star\star)$ %_{\tilde\Theta_{1}'\circ\tilde g\circ\tilde \Theta_{1}^{-1}}$
is also commutative. We leave the details of the proof to the reader.
\end{dem}

\section{Monodromy and projective holonomy}
\subsection{Holonomy representation of a JSJ block}\label{holblocJSS}
First of all, we recall some classical notions which will be used in the sequel. For a curve $\mc D$, the \emph{valence} of an irreducible component $D$ of $\mc D$ is the number $v(D)$ of irreducible components of $\mc D$ other than $D$ meeting $D$. We call \emph{dead branch of {the exceptional divisor} $\mc E_S$ any} maximal connected union of  components of $\mc E_{S}$ having  valence two in $\mc D_{S}$, except for one of them which has valence one in $\mc D_{S}$. Every dead branch $\mc M$ has a single \emph{attaching point} belonging to a unique component of $\mc E_{S}$ having valence at least three.\\

Let  $F :{\mb B}\to \mb C$ (resp. $F' : {\mb B}'\to \mb C$), be a reduced equation of {the separatrix curve} $S$ (resp. $S'$), and consider $\mc L$, $\mc L'$ two local data as in (\ref{systloc}). For $\varepsilon>0$ small enough, {the composition of these equations by the reduction morphisms $E_S$ and $E_{S'}$  define} the real smooth hypersurfaces $\{|F\circ E_S|=\varepsilon\}$ and $\{|F'\circ E_{S'}|=\varepsilon\}$. {These hypersurfaces bound} the \emph{Milnor tubes}  of $\mc D_S$ and $\mc D_{S'}$,
  \begin{equation}\label{tubemilnorB}
  \mc T_\varepsilon:=\{|F\circ E_S|\leq\varepsilon\}\subset \mc B_S\,,\quad \mc T'_\varepsilon:=\{|F'\circ E_{S'}|\leq\varepsilon\}\subset \mc B_{S'}\,,
\end{equation}
{and they }are transverse to the hypersurfaces  $\{|x_s|=1\}$, $\{|y_s|=1\}$ and $\{|x_{s'}|=1\}$, $\{|y_{s'}|=1\}$, for all  $s\in \Sing(\mc D_S)$, $s'\in \Sing(\mc D_{S'})$, as well as to the spheres  $E_S^{-1}(\partial {\mb B})$ and  $E_{S'}^{-1}(\partial {\mb B}')$. {We extend for subset of these Milnor tubes, the convenient notation introduced in (\ref{notast})}.
\begin{convs}\label{convtubemiln} For  ${\mc A}\subset \mc T_{\varepsilon }$ and  ${\mc A'}\subset \mc T_{\varepsilon '}$, we denote
$$
\wt{{\mc A}}^\ast:=q^{-1}(E_S({\mc A}))\,,\qquad
\wt{{\mc A'}}{}^\star:=q'{}^{-1}(E_{S'}({\mc A'}))\,,
$$
and we consider the following universal coverings:
$$\un q:=E^{-1}_{S}\circ {q }_{|\wt {\mc T}^\ast_\varepsilon }: \wt {\mc T}^\ast_\varepsilon \too \mc T_\varepsilon ^{\ast}\,,\quad
\un q':=E^{-1}_{S'}\circ {q'}_{|\wt {\mc T}'{}^\star_{\varepsilon'} } : \wt {\mc T}'{}^\star_{\varepsilon'} \too \mc T'{} ^{\star}_{\varepsilon '}\,.$$
Given a germ of homeomorphism $\phi _S:({\mc A},\mc D_S)\to ({\mc A'},\mc D_{S'})$  and a pro-germ at infinity $\wt \phi  _\infty: (\wt{{\mc A}}^\ast,\infty)\to(\wt{{\mc A'}}{}^\star,\infty)$ lifting $\phi  _S$,
we will say that the triple  $(\phi_S , \wt\phi _\infty, h)$ \emph{carries out over  ${\mc A}$ and ${\mc A'}$} {a geometric} {conjugation} (resp. {a geometric} $\mc N$-{conjugation}), $(\mf g, {h})$ , if the triple $(\phi_{S}^\flat, \wt\phi _{\infty}, h)$, with $\phi _S^\flat :=E_{S'}\circ\phi _S\circ {E_S}_{|{\mc A}^\ast}^{-1}$, is a realization of the {a geometric} {conjugation} (resp. {a geometric} $\mc N$-{conjugation}), $(\mf g,{h})$ over  $E_S({\mc A})$ and $E_{S'}({\mc A'})$.
\end{convs}

Let  $D$ be an irreducible component of $\mc D_S$ having valence ${v=v(D)}\geq 3$. If $D$ has  ${v-r}$ attaching points belonging to dead branches, we numerate the points $s_1,\ldots,s_{{v}}$ of $\Sing(\mc D_S)\cap D$ in such a way that   $\{s_j\,|\;j>{r}\}$ is the set of attaching points of the adjacent dead branches to $D$. We also denote
 \begin{equation}\label{comptrouee}
    D^\sharp:= D\setminus\cup_{j=1}^{{r}}\{|x_{s_j}| < 1\}\,,\quad
     D^\circ:= D\setminus\cup_{j=1}^{{v}}\{|x_{s_j}| < 1\}\,.
 \end{equation}
{We assume} that for $j=1,\ldots, {v}$, $y_{s_j}=0$ is a local reduced equation of  $D$.

 \begin{defin} We define the
 \emph{Jaco-Shalen-Johannson block} (JSJ for short) $B_{D}(\varepsilon)$ of $\mc T_{\varepsilon}$ associated to $D$ as the adherence of the connected component of  $\mc T_\varepsilon\setminus\cup_{j=1}^{{r}}\{|x_{s_j}|=1\}$ which contains  $D^\sharp$.
%In the same way  $B_{D}(\varepsilon)^\ast$ will be called the  JSJ block of $\mc T_\varepsilon^\ast$ associated to $D$.
 \end{defin}

\noindent We fix  $\varepsilon>0$ small enough and we denote $B_{D}(\varepsilon)$ and $\mc T_{\varepsilon}$ simply by $B_{D}$ and $\mc T$. {Then  the following properties hold, cf.   \cite{MarMat}:}
\begin{itemize}
  \item   $\mc T$ is a deformation retract of  $\mc B_S$, and for  $0<\varepsilon'\leq\varepsilon$, $B_D(\varepsilon')$ and $\mc T_{\varepsilon'}$  are deformation retract of  $B_D$ and  $\mc T$ respectively;
      \item $B_D^\ast$ is incompressible in $\mc T^\ast$, hence it is so in $E^{-1}_S({\mb B} ^\ast)$;
\item   a presentation of the fundamental group of $B_{D}^{\ast}$ by generators and relations can be obtained in the following way.
We consider loops $\gamma _1,\ldots, \gamma _{{v}}$ in
$B_D\cap\pi_D^{-1}(D^\circ)$ having the same origin $m$ and such that the projections $\pi _D\circ\gamma _j$ are the boundaries of closed conformal disks  ${\mc V}_j\subset D$ satisfying  $\inte{\mc V}_j\cap \Sing(\mc D_S)={\mc V}_j\cap \Sing(\mc D_S)=\{s_j\}$; then we consider a loop $c$ in the fiber
$\Delta  :=\pi _D^{-1}(m_0)$, $m_0:=\pi _D(m)$, with the same origin $m$ and having rotation index one with respect to the point $m_{0}$. Then we  have
               $$
               \pi _1(B_D^\ast, m) =\left< \dot{c},\;\dot{\gamma} _1,\ldots \dot{\gamma} _{{v}}\;\left|\; [\dot{\gamma} _j,\,\dot{c} ]=1,\dot{\gamma} _k^{p_k}= \dot{c}^{q_k}\right. \right>_{\stackrel{j=1,\ldots, {v}} {\ssstyle k={r}+1,\ldots,{v}}}\,,$$
where $\mathrm{gcd}(p_k, q_k)=1$ and $-\frac{q_k}{p_k}$ is the Camacho-Sad index of $\underline{\F}$ along $D$ at the point $s_k$;
\item the germ of $\underline{\F}$ at each  point $s=s_{{r}+1},\ldots,s_{{v}}$, possesses a holomorphic first integral that can be written as $x_{s}^{p_k}y_{s}^{q_k}A(x_{s},y_{s})$,
$A(0,0)\neq 0$.
\end{itemize}

\noindent

{Since $\mc H_{D}(\dot{\gamma })$ is the map
sending a point $p$ of $\Delta$ to the end of the path having origin $p$ and lifting  $\gamma ^{-1}$ in the leaf of $\underline{\F}$ passing through  $p$,}
the kernel of the holonomy representation of $\underline{\F}$ along $D^\circ$,
\begin{equation}\label{holo}
   \mc  H_{D} : \pi _1(D^\circ,\, m_0)=\mb Z\dot{\gamma } _{{r}+1}\ast\cdots\ast\mb Z\dot{\gamma} _{{v}}\too  \mathrm{Diff}(\Delta ,\, m_0)\,,
\end{equation}
contains the normal subgroup generated by the elements  $\dot{\gamma} _k^{p_k}$, for $k= {r}+1,\ldots , {v}$.
The morphism $\mc H_D$ factorizes through a morphism $\mc H_D^{\mathrm{orb}}$ defined on the {quotient group}
%\emph{$\pi _1$-orbifold of $D^\sharp$},
$$\pi _1^{\mathrm{orb}}(D^\sharp,\, m_0):={\pi _1(D^\circ,\, m_0)}/\ll\dot{\gamma} _{r+1}^{p_{r+1}},\ldots,\dot{\gamma}^{p_{{v}}} _{{v}}\gg\,.$$
The morphism $\pi _{D\ast}$ from $\pi _1(B_D\cap\pi_D^{-1}(D^\circ),\, m)$ onto $\pi _1(D^\circ,\,m_0)$ induced by the fibration  $\pi _D$, determines a morphism $\pi _{D\ast}^{\mathrm{orb}}$ which enters in the following exact sequence
\begin{equation}\label{piunorb}
    1\too \pi _1(\Delta ^\ast, m)=\mb Z\dot{ c}\too\pi _1(B_D^\ast,\,m)\stackrel{\pi _{D\!\ast}^{\mathrm{orb}}}{\longrightarrow}\pi_1^{\mathrm{orb}}(D^\sharp,\, m_0)\too 1\,.
\end{equation}
\begin{defin}\label{reprholon}
 We call the morphism
 $ \mc H_{B_D}:=\mc H_D^{\mathrm{orb}}\circ \pi _{D\ast}^{\mathrm{orb}}$,
$$
    \mc H_{B_D}  : \pi _1(B_D^\ast, \, m) \too \mathrm{Diff}(\Delta  , \, {m_0})\,,\quad \dot{\gamma} \mapsto \mc H_{B_D}(\dot{\gamma} ) = \mc H_D(\pi _D\circ\dot{\gamma} ),
$$ the
   \emph{holonomy representation of $\F$ along $B_{D}$} realized on the transverse section~$\Delta$.
\end{defin}

\subsection{Extended holonomy and monodromy}\label{ssectionholetendue}
With Conventions~\ref{convtubemiln} and the precedent notations, we denote by  $(\wt\Delta ^{\ast\alpha })_{\alpha \in \pi _0(\wt \Delta ^\ast)} $  and  $(\wt B_D^{\ast{\beta  }})_{\beta\in \pi _0(\wt B_D^{\ast})} $
the collection of connected components of $\wt\Delta ^\ast$ (resp. of  $\wt B_D^\ast$).
Thanks to the incompressibility of  $B_D^\ast$ in $E^{-1}_S({\mb B} ^\ast)$, the restriction of $\un q$ to each connected component  $\wt B_D^{\ast\beta }$ is a universal covering of $B_D^\ast$.
Thus, once we fix a point $m\in B_D^\ast$ and  $\wt m\in\un q^{-1}(m)\cap \wt B_D^{\ast\beta }$, the group  $\Gamma _\infty$ can be canonically identified with $\pi _1(\mc T^\ast, m)\simeq\pi _1({\mb B} ^\ast,  m)$ and the subgroup  $\Gamma _\infty(\beta )$ consisting of the elements  $\varphi$ in $\Gamma _\infty$ preserving $B_D^{\ast\beta  }$, can be identified with  $\pi _1(B_D^\ast,  m)$.
If $\wt m$ belongs to $\wt\Delta ^{\ast\alpha
}\subset \wt B_{D}^{\ast\beta}$ then we have the following exact sequence of groups:
$$1\too \Gamma _\infty(\beta , \alpha)
\too \Gamma_\infty(\beta)\stackrel{\sigma}{\too} \pi _1^{\mathrm{orb}}(D^\sharp,\, m_0)\too 1 \,,
$$
with $\Gamma _\infty(\beta, \alpha )\simeq\pi _1(\Delta ^\ast, m)$ denoting the subgroup consisting of those $\varphi\in \Gamma _\infty(\beta )$ which preserve  $\wt\Delta ^{\ast \alpha } $ and $\sigma $ is the well-defined morphism determined by
$$\sigma (\varphi):=\pi _D^{\mathrm{orb}}(\underline{q}\circ\mu_\varphi )\,,\qquad \varphi\in \Gamma _\infty^{\beta }\,,$$ where $\mu _\varphi$ is a path in $\wt B_D^{\ast\beta }$ whose endpoints are  $\wt m$ and $\varphi(\wt m)$. {Notice that  $\Gamma _\infty(\beta ,\alpha )$ is a normal subgroup in $\Gamma_{\infty}(\beta)$, because $\pi _1(\Delta ^\ast, m)$ is the center of $\pi _1(B_D^{\ast\beta })$; this means that each  $\varphi \in\Gamma _\infty(\beta, \alpha  )$ preserves every component  $\wt\Delta ^{\ast\alpha' }\subset B_D^{\ast\beta }$.}

\begin{prop}\label{mseudohol} If  $\wt\Delta ^{\ast\alpha }$ and $\wt\Delta ^{\ast\alpha'}$
are contained in the same connected component $\wt B_D^{\ast\beta }$, then there is a unique  pro-germ $h_{\alpha '\alpha } : (\wt\Delta ^{\ast\alpha },\infty)\to (\wt\Delta ^{\ast\alpha '}, \infty)$ commuting with the canonical pro-germs, {i.e.} %defined in  Section~\ref{progermcan}:
$\tau _{\wt\Delta ^{\ast\alpha '}}\circ h_{\alpha '\alpha }=\tau _{\wt\Delta ^{\ast\alpha }}$.
\end{prop}
\begin{dem} For all $\br U\in \mf U_{\F, \Sigma }$  put $U:=E^{-1}_S(\br U)$ and  denote by $W^{\alpha '}_U$ the $\wt\F$-saturation of  $\wt\Delta ^{\ast\alpha '}\cap\wt U^\ast$ inside $\wt B^{\ast\beta }_D\cap \wt U^\ast$. The map  $$h_U : {\mf O}_{U}^{\alpha '\alpha }:=W^{\alpha '}_U\cap\wt\Delta ^{\ast\alpha }\to\wt\Delta ^{\ast\alpha '}$$ obtained by following the leaves of  $\wt {\underline{\F}}_{|\mc W_U^{\alpha '}}$  is defined
without ambiguity because every leaf meets each transversal $\wt \Delta ^{\ast \alpha }$ and $\wt \Delta ^{\ast \alpha' }$ in at most one point, cf. Section~\ref{subsec.esp.des.feuilles}. We will see that ${\mf O}_{U}^{\alpha '\alpha }$
always contains a non-empty open set of type  $\wt V^\ast\cap \wt \Delta ^{\ast\alpha }$, with $V:=E_S^{-1}(\br V)$, $\br V\in \mf U_{\F,\Sigma }$.  In order to conclude, it suffices to put
$$
h_{\alpha '\alpha }:=( \,\limind_{\br V} \;h_{VU} \,)_{\br U\in \mf U_{\F, \Sigma }}\,,\quad \hbox{with} \quad h_{VU} :=h_{U}{}_{ |_{\wt V^{\ast}\cap \wt\Delta ^{\ast\alpha }}} : \wt V^{\ast}\cap \wt\Delta ^{\ast\alpha } \to \wt\Delta ^{\ast\alpha '}\,.
$$
\indent Consider  $\varphi\in\Gamma _\infty(\beta,\,\alpha )$ and choose $\br V\in \mf U_{\F,\Sigma }$, $\br V\subset\br U$, small enough so that all $p\in V\cap \Delta ^\ast$
be the origin of a path $\gamma  _p$ ending in $\Delta ^\ast$, which is contained in a leaf of the restriction of $\un\F$ to  $B_D^\ast\cap U\cap\pi_D^{-1}(D^\circ)$, and such that the homotopy  class of $\pi _D\circ \gamma _p$ in  $\pi _1^{\mathrm{orb}}(D^\sharp, m_0)$ coincides with  $\sigma (\varphi)$.
The lift $\mu _{\wt p}$ of $\gamma _p$ in $\wt B_D^{\ast\beta }$, passing through an arbitrary point $\wt p$ of $\un q^{-1}(p)\cap\wt\Delta ^{\ast\alpha }$, is contained in a leaf of $\wt {\un\F}$. We will see that its endpoint, which coincides with $h_U(\wt p)$, always belongs to  $\wt\Delta ^{\ast\alpha '}$; then the inclusion ${\mf O}_{U}^{\alpha '\alpha }\supset \wt V^\ast\cap \wt\Delta ^{\ast\alpha }$ will follow.

\indent We consider a path $\xi$ in $\wt\Delta^\ast $ having origin $\mu _{\wt p}(1)$ and endpoint in $\un q^{-1}(m)$, as well as a path $\delta$ having origin in $\un q^{-1}(m)$ and endpoint  $\wt p$. The homotopy class  $\dot{\zeta}\in \pi _1(B_D^\ast, m)$ of the loop  $\zeta:= \un q \circ (\delta\vee\mu_{\wt p}\vee \xi)$, satisfies:
$$
\pi _{D\ast}^{\mathrm{orb}}(\dot{\zeta})= \overline{\pi _{D}\circ \gamma  _p  } =\sigma (\varphi)\,,
$$
where  $\overline{\pi _D\circ \gamma _p  }$ denotes  the class of $\pi _D\circ \gamma  _p$ in $\pi _1^{\mathrm{orb}}(D^\sharp, m_0)$.
Let $\mu_\varphi $ be a path in  $\wt B^{\ast\beta }$ joining $\wt  m$ to $\varphi(\wt m)$. Thanks to (\ref{piunorb}), the homotopy class of the loop $\un q\circ \mu_\varphi$ in $\pi _1(B_D^{\ast},\,m)$  differs from $\dot{\zeta}$ in an element of $\pi _1(\Delta ^\ast, m)$; therefore the paths   $\delta\vee\mu_{\wt p}\vee \xi$ and $\mu _\varphi$ have their endpoints on the same connected component of  $\wt \Delta ^\ast$; the same property holds for  $ \mu _{\wt p}$ and $\mu_\varphi $.
\end{dem}

%\begin{obs}
For three components $\wt\Delta ^{\ast\alpha }$, $\wt\Delta ^{\ast\alpha' }$, $\wt\Delta ^{\ast\alpha '' }$ contained in $\wt B_D^{\ast\beta }$, we clearly have the relation
$$h_{\alpha ''\alpha '}\circ h_{\alpha '\alpha }=h_{\alpha ''\alpha }\,.$$
On the other hand, the above constructions are ``compatible'' with
 the action of
 $\Gamma _\infty(\beta)$, because $\Gamma_{\infty}$ preserves  $\wt\F$.
More precisely, with the  above notations, if $\varphi\in \Gamma _\infty(\beta)$ then the paths $\varphi\circ \mu _{\wt p}$ and $\mu _{\varphi(\wt p)}$ coincide. Thus, {denoting also by}
%if we denote by
$\varphi%_{\bullet} 
 : \pi _0(\wt \Delta ^{\ast})\iso \pi _0(\wt \Delta ^{\ast})$ the bijection induced by {deck transformation} $\varphi$, we have  that
$$\varphi\circ h_{\alpha' \alpha }\circ \varphi^{-1}= h_{\varphi%_{\bullet}
(\alpha ' )\varphi%_{\bullet}
(\alpha )}\,.$$ We easily deduce that the map
$$
\wt{\mc H}_D^{\alpha  } : \Gamma _\infty(\beta) \too \Aut (\wt\Delta ^{\ast\alpha },\,\infty)\,,\quad \varphi\mapsto  h_{\alpha \varphi%_{\ssstyle \bullet}
(\alpha )}\circ \varphi=h^{-1}_{\varphi%_{\ssstyle \bullet}
(\alpha )\alpha }\circ \varphi\,,
$$
is a morphism of groups. {In the case that $E_{\F}$ consists in a single blow-up this morphism was considered already in \cite{OBRGV} where it was called \emph{extended holonomy}. We adopt their definition in the general context.}
%\end{obs}
\begin{defin}
We will call $\wt{\mc H}_D^{\alpha }$, the  \emph{extended holonomy morphism of $D$ over  $\wt\Delta ^{\ast\alpha }$}.%, cf. \cite{OBRGV}.
\end{defin}
\noindent In order to justify this definition we note that if $\dot{\mu  }\in \pi _1( B_D^{\ast}, m)$ and $\varphi\in \Gamma _\infty(\beta )$ satisfy $\pi _{D\,\ast}^{\mathrm{orb}}(\dot{\mu  })= \sigma (\varphi)$, then $\wt{\mc H}_D^{\alpha  }(\varphi)$
is the lifting of the holonomy diffeomorphism $\mc H_{B_D}(\dot{\mu  })$ on the connected component  $\wt\Delta ^{\ast\alpha }$, considered as a universal covering of $\Delta ^\ast$. We finally obtain the following commutative diagram:
$$ (\diamond)_\varphi\qquad\hbox{
\begin{diagram}
(\Delta ,\,m_0) & \lOnto^{\phantom{aaa}\un q_{\ssstyle\infty}\phantom{aaaa}} & (\wt\Delta^{\ast\alpha },\infty) & \rInto^{\phantom{aaa}\tau_{\wt\Delta ^{\ast\alpha }}\phantom{a}} & \QF\\
\dTo{\mc H_{B_D}(\dot{\mu } )}&   &\dTo^{\wt{\mc H}_D^{\alpha }(\varphi)} & &  \dTo_{\monFS(\varphi )}\\
(\Delta ,\,m_0) & \lOnto^{\phantom{aaa}\un q_{\ssstyle\infty}\phantom{aaaa}} & (\wt\Delta ^{\ast\alpha },\infty) & \rInto^{\phantom{aaa}\tau_{\wt\Delta ^{\ast\alpha }}\phantom{a}} & \QF\\
\end{diagram}\,.
}\phantom{(\diamond)_\varphi\qquad}
$$

\medskip

\subsection{Relationship between holonomy and monodromy {conjugations}}
We fix a geometric ($\mc N$-){conjugation} $(\mf g, {h})$ between the monodromies $\monFS$ and $\monFSpr$ of the General Type foliations $\F$ and $\F'$.
Thanks to Corollary~\ref{reprfibree}, there exists a geometric representation  $(g, \wt g, h)$ of $(\mf g, {h})$, cf. Definition~\ref{conjugtopol},
where  $g : ({\mb B} , S)\to ({\mb B} ', S')$ is the germ of  an excellent homeomorphism and $G : (\mc B_S, \mc D_S)\to (\mc B_{S'}', \mc D_{S'})$ is its lifting to the total spaces of their reductions of singularities. Let $(\Delta  , m_0)$, $m_0\notin \Sing(\mc D_{S})$, be a fibre of the Hopf fibration of an irreducible component $D$ of $\mc D_S$, such that $(\Delta ', m'_0):=(G(\Delta ), G(m_0))$
is also a fibre of the Hopf fibration associated to the component $D':=G(D)$ of $\mc D_{S'}$. Finally, we will identify  $\Delta $ and $\Delta ' $ with their images $E_S(\Delta )$ and $E_{S'}(\Delta ')$ by the reduction of singularities maps.
\begin{teo}\label{monodimplholo}
If there is a realization $(\psi, \wt\psi, h)$ of $(\mf g, {h})$ over  $\Delta $ and  $\Delta '$, then $\psi$ and the restriction $G_{|D} :D\iso D'$ conjugate the holonomy representations associated to the irreducible components  $D$ and $D'$, i.e. the following diagram is commutative:
\begin{diagram}
 \pi _1(D^\circ, \, m_0)& \rTo^{\phantom{aaa}\mc H_D\phantom{aaa}} & \mathrm{Diff}(\Delta , m_0)\\
 \dTo{G_\ast}& {\sstyle \circlearrowleft}& \dTo{\psi_{\ast}}\\
  \pi _1(D'{}^\circ, \, m_0)& \rTo^{\phantom{aaa}\mc H_{D'}\phantom{aaa}} & \mathrm{Diff}(\Delta ', m_0')\
\end{diagram}
\noindent where $\psi_{\ast}(\varphi):=\psi\circ \varphi\circ\psi^{-1}$ and $G_\ast$ is induced by the restriction of  $G$ to $D^\circ$.%, cf. (\ref{comptrouee}).
\end{teo}

\begin{figure}[htbp]
\begin{center}
\begin{diagram}
(\Delta ^\ast,m_0) & & \rTo^{\mc H_D(\dot{\gamma})} & & (\Delta ^\ast , m_0) & & \\
& \rdTo_{\psi} & & & \uOnto^{\un q_{\ssstyle\infty}} & \rdTo_{\psi} & \\
\uOnto^{\un q_{\ssstyle\infty}} & & (\Delta '{}^\star,m'_0) & \rTo^{\mc H_{D'}(G_\ast(\dot{\gamma}))} & \HonV & & (\Delta '{}^\star, m'_0) \\
& & \uOnto^{\un q_{\ssstyle\infty}'} & & \dLine & & \\
(\wt\Delta ^{\ast\alpha },\infty) & \hLine & \VonH & \rTo^{\wt{\mc  H}_D^{\alpha }(\varphi)} & (\wt\Delta ^{\ast\alpha },\infty) & & \uOnto_{\un q_{\ssstyle\infty}'} \\
& \rdTo_{\wt\psi} & & & \daux & \rdTo_{\wt\psi} & \\
& & (\wt\Delta '{}^{\star\alpha '},\infty) & & \rTo^{\wt{\mc H}_{D'}^{\alpha '}(\mf g(\varphi))} & & (\wt\Delta '{}^{\star\alpha '},\infty) \\
\dInto^{\tau_{\wt\Delta ^{\ast\alpha }}} & & \dInto^{\tau_{\wt\Delta '{}^{\star\alpha '}}} & & \dTo^{\tau_{\wt\Delta ^{\ast\alpha }}} & & \\
\QF & \hLine & \VonH & \rTo^{\monFS(\varphi)} & \QF & & \dInto_{\tau_{\wt\Delta '{^{\star\alpha '}}}} \\
& \rdTo_{h} & & & & \rdTo_{h} & \\
& & \QFpr & & \rTo^{\monprFS(\mf g(\varphi))} & & \QFpr \\
\end{diagram}
\caption{Diagram concerning the {conjugation} between the holonomies and monodromies.}\label{diagramacubico}
\end{center}
\end{figure}

\begin{dem} Consider $\dot{\gamma }\in \pi_1(D^\circ, m_0)$, $\dot{\mu  }\in \pi _1(\wt B_D^{\ast}, m)$ and $\varphi\in \Gamma _\infty(\beta )$  such that $\pi _{D\,\ast}^{\mathrm{orb}}(\dot{\mu  })= \sigma (\varphi)=\dot{\gamma }$ and therefore $\mc H_{B_D}(\dot{\mu } ) = \mc H_D(\dot{\gamma })$. Consider also the diagram of  Figure~\ref{diagramacubico}.
Both frontal sides (behind and ahead) are constituted by the commutative diagrams  $(\diamond)_\varphi$ and $(\diamond)_{\mf g(\varphi)}$;
both lateral sides are constituted by the commutative diagram  $(\star)$ %_{\wt{\psi}}$
and that one expressing  that
  $\wt \psi$ lifts $\psi$;
the commutativity of the bottom horizontal diagram, $h\circ \monFS(\varphi)=\monprFS(\mf g(\varphi))\circ h$, follows from the
fact that $(\mf g, h_{*})$ is a {conjugation} between $\monFS$ and $\monprFS$.
Since the canonical pro-germs $\tau _{\wt \Delta ^{\ast\alpha }}$ and $\tau _{\wt \Delta '{}^{\ast\alpha '}}$
are  monomorphisms, cf.  Proposition~\ref{monoprogermcan}, the
median horizontal diagram is also commutative, i.e. $\wt\psi\circ\wt {\mc H}_D^{\alpha }(\varphi)=\wt {\mc H}_{D'}^{\alpha '}\circ\wt\psi$. Finally, the commutativity of the top horizontal diagram, $$\psi\circ \mc H_D(\dot{\gamma })= \mc H_{D'}(G_\ast(\dot{\gamma }))\circ\psi\,,$$ follows from the fact that the pro-germs
 $\un q_{\ssstyle\infty}$ and $\un q'_{\ssstyle\infty}$ are epimorphisms.
\end{dem}

\section{Statements and proofs of Theorems I and II}\label{sect.th.classification}
We keep the notations introduced in the precedent sections about the germs of foliations  $\F$ et $\F'$, in particular  (\ref{notast}), (\ref{notrevuniv}), (\ref{notationstar}) and (\ref{notationqpr}).
We shall introduce here 	{two statements that imply}
%the more general statements of the main theorems of this paper, from which will immediately follow
Theorems I and II in the Introduction.

\setcounter{teo}{0}
\begin{teo}[\textbf{of invariance}]\label{invariance}
Assume that $ \F $ and $ \F '$ are of General Type
%,i.e. they satisfy the condition (GT) given in the introduction,
and they are
conjugate  by a germ of transversely holomorphic (resp. $\mc N$-transversely holomorphic) homeomorphism $ \Psi: (\mb C^2, 0) \iso (\mb C^2,0)$. Consider  a germ $ g: (\mb C^2, 0) \iso (\mb C ^2,0)$ of  excellent homeomorphism  fundamentally equivalent to $ \Psi $, cf. Section~\ref{subsecmarquage}, and  denote by $G : (\mc B_S,\mc D_S)\to(\mc B'_{S'},\mc D_{S'})$ its lifting. %, i.e. $E_{S'}\circ G=g\circ E_S$.
Then,
\begin{enumerate}
  \item \label{invCS} for each irreducible component  $D$ of $\mc D_S$ and for each singular point $s\in \Sing(\underline{\F})\cap D$, we have equality of Camacho-Sad indices:
  \begin{equation}\label{asts}
     \CS(\underline{\F}, D, s)= \CS(\underline{\F}', G(D), G(s))
  \end{equation}
  \item\label{conjreal} there is a {geometric} {conjugation} (resp. a {geometric} $\mc N$-{conjugation}), $(\mf g, {h})$ between the monodromies $\monFS$ and $\monFSpr$ realizable over some $\mc N$-collections of transversals (i.e. satisfying Condition ($\sigma$) in the Introduction) of $\F$ and $\F'$,
such that $\Psi$ is a representative of the marking determined by $\mf g$.
\end{enumerate}
\end{teo}

\noindent The following theorem can be considered as a sort of converse of the precedent result.
\begin{teo}[\textbf{of classification}]\label{classification} Let  $(\mf g, {h})$ be a geometric $\mc N$-{conjugation} between the monodromies  $\monFS$ and  $\monFSpr$ of two General Type foliations, which possesses a realization  $(\psi , \wt \psi , h)$  over $\mc N$-collections of transversals   $\Sigma $ and $\Sigma '$ of $\F$ and $\F'$.
Denoting by $(g,\wt g, h)$ a geometric representation of $(\mf g, {h})$, assume that the following conditions
%\footnote{Properties (\ref{memesep}) and (\ref{egCSseparatrices}) do not depend on the choice of the geometric representation of $(\mf g,{h})$.}
 hold:
\begin{enumerate}
\item\label{memesep}
for each connected component $\br \Sigma$ of $\Sigma$
intersecting an irreducible component  $\br S$ of $S$ we have $\psi(\br\Sigma)\cap g(\br S)\neq\emptyset$,
\item\label{egCSseparatrices}
along each irreducible component $\br S$ of $S$, we have coincidence of Camacho-Sad indices: $ \CS(\underline{\F}, \br{\mc S}, \br s)= \CS(\underline{\F}', \br{\mc S'}, \br s')$, where $\br{\mc S}$ (resp. $\br{\mc S}'$), denote the strict transform of $\br S$ (resp. $g(\br S)$), and $\br s\in \Sing(\underline{\F})$ (resp. $\br s'\in \Sing(\underline{\F}')$), are their attaching points in the exceptional divisors.
\end{enumerate}
Then there exists a homeomorphism $\Psi$ defined on an open neighborhood  $U$ of $S$ onto an open neighborhood $U'$
of $S'$ and there is a lift
 $\wt \Psi : \wt U ^\ast\iso\wt U '{}^{\star}$ of $\Psi$, such that:
\begin{enumerate}
\item[(a)] $\Psi (\Sigma \cap U)\subset \Sigma '\cap U'$, the germ of $\Psi _{|\Sigma \cap U}$ at  the finite set $\Sigma\cap S$
equals $\psi $  and the germ at infinity of $\wt \Psi _{|\wt \Sigma \cap \wt U}$ coincides with  $\wt \psi $,
  \item[(b)] $\Psi$ is $\mc N$-excellent, conjugates $\F_{|U }$ to $\F_{|U'}'$ and it is $\mc N$-transversely holomorphic,
  \item[(c)] denoting by  $\Psi _S$ the  germ of $\Psi $ along  $S$ and by $\wt \Psi _\infty$ the pro-germ at infinity of $\wt\Psi $, then $(\Psi_S,\wt\Psi_\infty,h)$ is a realization of $(\mf g,{h})$ over
the complement of some $\mc N$-separators of $\F$ and $\F'$.
\end{enumerate}
\end{teo}

\noindent In fact we will prove in Section~\ref{invtousCSpreuve} that Condition~(\ref{egCSseparatrices}) in  Classification Theorem~\ref{classification} is equivalent to Assertion~(\ref{invCS}) in  Invariance Theorem~\ref{invariance}.
The relationship between these two results and Theorems I and II in the introduction come from the following {property of transversal rigidity.}

\begin{defin}\label{deftrasvrigide}
We say that the germ of $\F$ at $(0,0)$ is \emph{transversely rigid} (resp. \emph{$\mc N$-transversely rigid}), if every germ of homeomorphism preserving the orientations of $(\mb C^2, 0)$ and those of the leaves, and conjugating $\F$ to a General Type foliation, is necessarily transversely  (resp.  $\mc N$-transversely) holomorphic.
\end{defin}

After the extended version (TRT) of the Transverse Rigidity Theorem of J. Rebelo \cite{Rebelo} stated in the introduction, the hypothesis (G) on $\F$ implies the $\mc N$-transverse rigidity of $\F$. As we have already point out in the introduction, the genericity of this property in the sense of the Krull topology was proved in \cite{LeFloch}. We remark that there are other interesting, although more particular, situations inducing also the transverse rigidity of the foliation, cf. \cite{Loray}.

\begin{cor}\label{coroltrriggeneral}
Let  $\F$ and $\F'$ be two germs of General Type $\mc N$-transversely rigid foliations which are conjugated by a germ of orientations preserving homeomorphism $\Psi_0$. Then there is a germ of $\mc N$-excellent homeomorphism  $\Psi$  conjugating $\F$ and $\F'$. In particular, denoting by $\Psi^\sharp : \mc D_S\iso \mc D_{S'}$ the restriction to the total divisors of the lifting of $\Psi$ over the reduction of singularities, the following properties hold:
 \begin{enumerate}[(a)]
  \item for each $s\in \Sing(\F)$, the germs of $\underline{\F}$ at  $s$ and that of $\underline{\F}'$ at  $\Psi^\sharp(s)$ are holomorphically equivalent;
  \item
  for each irreducible component  $D$ of $\mc D_S$, the holonomy representations  $\mc H_D$ and $\mc H_{\Psi ^\sharp(D)}$  defined in (\ref{holo}) are holomorphically conjugate via $\Psi^\sharp$.
      \end{enumerate}
In addition, if $\mr{Node}(\un\F)=\emptyset$ then $\Psi$  is fundamentally equivalent to $\Psi_0$, cf. Section~\ref{subsecmarquage}.
\end{cor}
\noindent Assertion (b) in Corollary~\ref{coroltrriggeneral} means that if $\Delta $ and $\Delta':=\Psi(\Delta)$ are holomorphic curves transverse to  $D$ and $D':=\Psi ^\sharp(D)$, at the points $m\in D\setminus \Sing(\underline{\F})$ and $m':=\Psi^\sharp(m)$, then the diagram (\ref{diagdholonomie}) in the introduction commutes taking $\psi=\Psi_{|\Delta}$.

\begin{dem} We can apply first  Theorem of Invariance~\ref {invariance}, because $\Psi $ is $\mc N$-transver\-sely holomorphic; then we apply Classification Theorem~\ref{classification}.
\end{dem}

Notice that this corollary is a more precise statement of Theorem~I. On the other hand, thanks to the transverse rigidity theorem, Theorem~II follows from the invariance and classification
theorems below, using the additional generic hypothesis (G) considered in the introduction.

\section{Peripheral structure of a germ of curve}\label{structuper}
{Before proving Theorems~\ref{invariance} and \ref{classification},  we shall examine some auxiliary topological notions that we will need in the sequel. In fact, this section deals uniquely with curves and there is no foliation there.}
\subsection{Peripheral groups}\label{subsectstrper}
Following the previously introduced notations
(\ref{notast}), (\ref{notationstar}), (\ref{tubemilnorB}) and Conventions~\ref{convtubemiln}, {let $\br S$ be an irreducible component of $S\subset\mb B$. We consider a tubular neighborhood $W_{\br S}$ of   $\br S^\circ:=\br S\setminus\{0\}$ in $\mb B\setminus\{0\}$. The pair $(W_{\br S},\br S^\circ)$ is homeomorphic to  $(\br S^\circ \times \mb D_{1}, \br S^\circ \times \{0\})$.} Let  $s\in \Sing(\mc D_S)$ be the attaching point of the strict transform $\br{\mc S}$ of $\br S$. Up to permutation of the coordinates  $(x_s, y_s)$ of the local datum  fixed in (\ref{systloc}), we assume that $x_s=0$ is a local reduced equation of  $\br{\mc S}$. We choose  $\varepsilon>0$ small enough so that  $\mc W_{\br{\mc S}}^\ast:=E_S^{-1}(W_{\br S}^\ast)$ retracts over the 2-torus $\{|x_s|=\varepsilon, |y_s|=1\}$.
\begin{prop}\label{incompbordmiln} $W_{\br S}^\ast$ is incompressible in ${\mb B}^\ast$.
\end{prop}
\begin{dem}
It suffices to show the incompressibility of
the torus $$\{|x_s|=\varepsilon, |y_s|=1\}$$ inside $\mc B_S^\ast$. This can be done by using Van Kampen's Theorem, see for instance the construction of an open neighborhood of  $\mc D_S$ by ``boundary assembly'' made in
\cite{MarMat}.
\end{dem}
Consider the loops  $m$ (resp. $p$), in  $\mc  W_s^\ast$,
having the same origin, defined by $(x_s,y_s)\circ m(t)= (\varepsilon e^{2i\pi t},1)$ (resp. $(x_s,y_s)\circ p(t)=(\varepsilon,e^{2i\pi t})
$).
At the point $\br c:=E_S(m(0))$, the homotopy classes $\mf m_{\br c}$ (resp.  $\mf p_{\br c}\in \pi _1(W_{\br S}, \br c)$), of the loops $E_S\circ m$ (resp. $E_S\circ p$), allow to decompose $\pi _1(W_{\br S}^\ast, \br{c})=\mb Z \mf m_{\br c}\oplus\mb Z \mf p_{\br c}$. The abelianity of this group implies that the isomorphism from $\pi _1(W_{\br S}^\ast, c_1)$ onto $\pi _1(W_{\br S}^\ast, c_2)$, induced by a path joining the points $c_{1}$ and $c_{2}$ inside $W_{\br S}^\ast$,   does not depend on the particular choice of this path.
Thus, the direct sum decomposition of $\pi _1(W_{\br s},\br c)$ is canonical, i.e. it can be unambiguously defined  for every base point in  $W_{\br S}^\ast$:
$$\mc P_{\br S,\,c}:=\pi _1(W_{\br S}^\ast, c)=\mb Z \mf m_c\oplus\mb Z \mf p_c\subset \pi _1({\mb B}^\ast, c)\,,\quad c\in W_{\br S}^\ast\,.$$
\begin{defin}\label{periph}
We will call  $\mf m_c$ the \emph{meridian}, $\mf p_c$  the \emph{parallel} and $\mc P_{\br S,\,c}$ the \emph{peripheral subgroup} associated to the component  $\br S$ at the point $c$.
\end{defin}

The following geometric property states that this decomposition is ``intrinsic'', see \cite{Heil,MarMatMarq}.

\begin{prop}\label{normalisateur} The subgroup  $\mc P_{\br S,\,c}$ coincides with its normalizer inside $\pi _1({\mb B}^\ast, c)$, i.e. $\Big(\zeta \in \pi _1({\mb B}^\ast,c)\;\;\hbox{\rm and}\;\;\zeta\,\mc P_{\br S,\,c}\, \zeta ^{-1}\subset\pi _1(W_{\br S}^\ast, c)\Big)\Rightarrow\zeta\in\mc P_{\br S,\,c}$.
\end{prop}
\noindent We immediately deduce:
\begin{cor}\label{per}  The direct sum decomposition  $$P=\mb Z\mf m_{\ssstyle P}\oplus\mb Z\mf p_{\ssstyle P}\,,\qquad \mf m_{\ssstyle P}:= \zeta \,\mf m_c\,\zeta ^{-1}\,,\quad \mf p_{\ssstyle P}:= \zeta \,\mf p_c\,\zeta ^{-1}\,,$$ of every subgroup $P=\zeta \,\mc P_{\br S, c}\,\zeta ^{-1} $, $\zeta \in \pi _1({\mb B}^\ast, c)$ conjugated to  $\mc P_{\br S, c}$ is intrinsic,  i.e. it does not depend on $\zeta \in \pi _1({\mb B}^\ast, c)$.
\end{cor}

\subsection{Conjugation of peripheral structures}
We will see that the canonical meridians and parallels introduced in Definition~\ref{periph}, associated to the irreducible components of $S$, are topological invariants.
\begin{teo}\label{invstruperi}
Let $U$ be an open neighborhood of $S$ in ${\mb B}$ and let $\Phi$ be a homeomorphism from $U$  onto a neighborhood  $U'$ of $0$ in ${\mb B}'$,  such that   $\Phi(S)=S'\cap U'$. Then for each irreducible component  $\br S$ of $S$ and for all point  $c$ in a tubular neighborhood of  $\br S\setminus \{0\}$ inside ${\mb B}$, the isomorphism  $\Phi_\ast$ from $\pi _1({\mb B}^\ast, c)$ onto $\pi _1({\mb B}'{}^{\!\ast},c')$ induced by $\Phi$, sends respectively the meridian $\mf m_{c}$ and the parallel $\mf p_{c}$ associated to $\br S$ to the meridian $\mf m'_{c'}$ and the parallel $\mf p'_{c'}$ associated to the component $\Phi(\br S)$, at the point $c':=\Phi(c)$.
\end{teo}
\begin{dem} First we note that  $\Phi_\ast$ induces an isomorphism from the peripheral group $\mc P_{\br S,\,c}$  of $\br S$,  onto the peripheral group  $\mc P'_{\br S',\,c'}$ of $\br S':=\Phi(\br S)$. Indeed, we can consider a tubular neighborhood $W$ of $\br S$ and two tubular neighborhoods $W'$ and $W''$  of $\br S'\setminus\{0\}$ in ${\mb B}'$, as well as a ball   ${\mb B}''\subset {\mb B}'$ centered at the origin, such that  $W''\cap{\mb B}''\subset \Phi(W)\subset W'$. These inclusions induce two $\mb Z$-linear morphisms at the fundamental group level,
 $$\mb Z^2\simeq \mc P_{\br S',\, c'}'\to\Phi_\ast(\mc P_{\br S,\,c})\to \mc P'_{\br S',\,c'}\simeq\mb Z^2\,,$$
whose composition is an isomorphism. Hence
 \begin{equation}\label{isper}
 \Phi_\ast(\mc P_{\br S,\,c})=\mc P'_{\br S',\,c'}.
\end{equation}
\noindent The Marking Theorem~\ref{marquageprecis} provides the existence of an excellent homeomorphism $g$ fundamentally equivalent to  $\Phi$. We can assume that  $W''\cap{\mb B}''\subset g(W)\subset W'$ and $g$ induces an isomorphism $g_\ast$ from $\mc P_{\br S, c}$ onto $\mc P_{\br S', g(c)}$. Clearly   $g_\ast(\mf m_c)=\mf m'_{g(c)}$ and $g_\ast(\mf p_c)=\mf p'_{g(c)}$, because the lifting  $E'_{S'}{}^{\!-1}\circ g\circ E_S$ extends to the exceptional divisor.
The fundamental equivalence between $g$ and $\Psi$ implies the existence of an element  $\zeta $ of $\pi _1({\mb B}'{}{\!^\ast}, c')$ such that
\begin{equation}\label{imper}
I _{\zeta }\circ\Phi_\ast = \kappa \circ g_{\ast} : \pi _1({\mb B}^\ast, c')\longrightarrow\pi _1({\mb B}^\ast, c')\,,
\end{equation}
where  $I_{\zeta }$ denotes the interior automorphism of $\pi _1({\mb B}^\ast, c')$ determined by  $\zeta $ and $\kappa $ is the canonical isomorphism from $\pi _1({\mb B}^\ast, g(c))$ onto $\pi _1({\mb B}^\ast, c')$, determined by an arbitrary path in $W'$ joining  $g(c)$ to  $c'$. The relations (\ref{isper}) and (\ref{imper}) give the equality
$
\zeta \mc P'_{\br S',\,c'}\zeta ^{-1}= \mc P'_{\br S',\,c'}
$. By applying Proposition~\ref{normalisateur} we obtain that $\zeta$ belongs to  $\mc P_{\br S',\,c'}$. The restriction of  $I_{\zeta }$ to $\mc P_{\br S',\,c'}$ is the identity because this group is abelian. By restricting   (\ref{imper}) to $\mc P_{\br S,\,c}$ we obtain the relation $\Phi_\ast=\kappa \circ g_\ast$.  Hence  $$\Phi_\ast(\mf m_c)=\kappa (g_\ast(\mf m_c))=\kappa (\mf m'_{g(c)})=\mf m'_{c'}\,;$$
and analogously $\Phi_\ast(\mf p_c)=\mf p'_{c'}$.
\end{dem}

\section{Proof of Invariance Theorem \ref{invariance}}\label{7}

\subsection{Proof of Assertion (\ref{conjreal})}\label{demasrtdeux} By  Proposition~\ref{lemmeisotopie}, we can compose $\Psi$ on the left by a germ of homeomophism $\Theta _1 : ({\mb B}',S')\iso({\mb B}',S')$ which is $\F'$-isotopic to the identity.
Let $\Sigma$ and $\Sigma'$ be $\mc N$-collections of transversals of $\F$ and $\F'$ respectively, such that for each connected component $\br\Sigma$ of $\Sigma$ intersecting an irreducible component $\br S$ of $S$ we have that $\Psi(\br\Sigma)$ and $g(\br\Sigma)$ meet the same irreducible component $\br S'$ of $S'$. {We can} construct a $\F'$-isotopy $\Theta_{1}$ such that $\Theta_{1}(\Psi(\br\Sigma))$ and $\br\Sigma'$ define the same germ at $\br\Sigma'\cap \br S'$.
Assertion~(\ref{conjreal}) of Theorem~\ref{invariance} follows directly from Remarks~\ref{conjfeuilconjmon} and~\ref{conjgeom}.

\subsection{Invariance of the Camacho-Sad indices associated to the separatrix curve}\label{invcssepartr}
We will prove the equality (\ref{asts}) %in  Assertion (\ref{invCS}) of Theorem~\ref{invariance},
when  $s$ is the attaching point of the strict transform
 $\br{\mc{S}}$ of an irreducible component $\br S$ of $S$. Following the notations of Section~\ref{subsectstrper} we denote by  $\rho : \mc W_s\to \br{\mc S}\cap \mc W_s$ the  disk fibration such that $y_s\circ\rho =y_s$, by $\gamma _n$ the loop contained in  $x_s=0$, such that $y_s\circ\gamma _n(t) :=e^{2i\pi nt}$,  $0\leq t\leq 1$, by $q$  the point having coordinates  $(0,1/2)$ and finally by  $T$ the transverse section $\rho ^{-1}(q)$. \\

Consider a sequence $(q_n)_{n\in\mb N}$ of points in  $T$  tending to $q$, such that the loop $\gamma _n$ lifts, via  $\rho $, to a path  $\Gamma _n$ contained in a leaf of
$\underline{\F}$. {We can see} that such a sequence always exists and that
\begin{equation}\label{camsa}
    \CS(\underline{\F},\br{\mc S}, s)= \lim_{n\to\infty}\frac{1}{2i\pi n}\int_{\Gamma _n}\frac{dx_s}{x_s}\,.
\end{equation}
Fix a real number $\theta _n\in]-\pi ,\pi ]$ different from the  arguments of $q_n$ and that of $\Gamma _n(1)$ and choose a path $\xi_{n}$ in
$$T\cap\left\{\arg(x_s)\neq \theta _n\,,\;0<|x_n|<1/n\right\}\,,$$ having endpoints  $q_n$ and  $\Gamma _n(1)$. Since the real part of  $\frac{1}{2i\pi n}\int_{\xi _n}\frac{dx_s}{x_s}$ is bounded, the real part of the Camacho-Sad index  (\ref{camsa}) is given by
 $$\mathrm{Re}(\CS(\underline{\F},\br{\mc S}, s)) = \lim_{n\to\infty}\frac{I_n}{n}\,,\quad\hbox{ where }\quad  I_n:=\frac{1}{2i\pi }\int_{\Gamma_n{\ssstyle\vee}\,\xi_n}\frac{dx_s}{x_s}.$$
Using the peripheral structure of %the peripheral group
{$\mc P_{\br S,\,c_n}$ given in Corollary~\ref{per}},
%$\subset \pi _1({\mb B}^\ast, c_n)$,
the homotopy class of the loop
 $\Gamma _n{\sstyle\vee}\xi_n$
 in
$\pi _1({\mb B}^\ast, c_n)$, $c_n:=E_S(q_n)$, can be decomposed as
%in the following way:
 $$\overline{\Gamma _n{\sstyle\vee}\xi_n}=I_n\mf m_{c_n}\,+\,n\mf p_{c_n}\,.$$

If $\alpha _n$ is an arbitrary path in $T\setminus\{q\}$ with endpoints $q_0$ and $q_n$ then the homotopy class of the loop
$\lambda _n:=\alpha _n{\sstyle\vee}\Gamma _n{\sstyle\vee}\xi_n{\sstyle\vee}\alpha _n^{-1}$ in  $\pi _1({\mb B}^\ast, c_0)$  is
 $$   \overline{\lambda _n}=I_n\mf m_{c_0}\,+\,n\mf p_{c_0}\,.
$$
Fix now the same data at the attaching point
 $s':=G(s)$ of the strict transform $\br{\mc S}'$ of $\br S':=\Psi(S)$. We denote by $(x_{s'},y_{s'}) : \mc W'_{s'}\iso\mb D_{1}^2$ the local coordinates at $s'$ determined by the local datum
 $\mc L'$, by $\rho ': \mc W'_{s'}\to \br{\mc S'}\cap \mc W'_{s'}$ the  disk fibration defined by $y_{s'}\circ \rho '=y_{s'}$, and by  $q'\in \br{\mc S}$ the point having coordinates $(0,1/2)$. We also denote   $T':=\rho' {}^{-1}(q')$. It is easy to see that, after composing it by a homeomorphism $\F'$-isotopic to the identity, $\Psi$ satisfy the following properties:
\begin{itemize}
  \item $\Psi(V)\subset V'$,
  where  $V$ is the image by  $E_S$ of a tubular neighborhood of the circle  $\{x_s=0, |y_s|=1\}$ inside $\{|x_s|\leq 1, |y_s|=1\}$ and $V'$ is the image by  $E_{S'}$ of the torus  $\{|x'_{s'}|\leq 1, |y_{s'}|=1\}$;
  \item $\Psi_{|V}$ conjugates the fibrations, i.e. $E'_{S'}\circ\rho '\circ E_{S'}^{-1}\circ\Psi_{|V}=\Psi\circ E_S \circ\rho\circ E_S^{-1} $.
\end{itemize}
As in (\ref{camsa}) we have the  following equality:
$$
 \CS(\underline{\F}',\br{\mc S}', s')= \lim_{n\to\infty}\frac{1}{2i\pi n}\int_{\Psi\circ\Gamma _n}\frac{dx'_{s'}}{x'_{s'}}\,.
$$
The variation of the argument of  $x'_{s'}\circ\xi_n$ is bounded because the restriction of $\Psi$ to $T$ is holomorphic (recall that $\Psi$ is transversely holomorphic on $V$). Thus,
\begin{equation}\label{recspr}
    \mathrm{Re}(\CS(\underline{\F'},\br{\mc S}', s')) = \lim_{n\to\infty} \frac{J_n}n\,,\quad J_n:=\frac{1}{2i\pi }\int_{\Psi\circ(\Gamma _n{\ssstyle\vee}\,\xi_n)}\frac{dx'_{s'}}{x'_{s'}}.
\end{equation}
Clearly the homotopy class of $\Psi\circ \lambda _n$ in $\pi _1({\mb B}^\ast, \Psi(q_0))$ is %equal to
 $J_n\mf m'_{\Psi(c_0)}\,+\,n\mf p'_{\Psi(c_0)}$ and therefore
$$J_n\mf m'_{\Psi(c_0)}\,+\, n\mf p'_{\Psi(c_0)}=\Psi_\ast(I_n\mf m_{c_0}+
n \mf p_{c_0})\,.
$$
Theorem~\ref{invstruperi} provides the equality $I_n=J_n$. Thanks to  (\ref{recspr}),   $\CS(\underline{\F},\br{\mc S}, s)$ and $\CS(\underline{\F}',\br{\mc S}', s')$ have the same real part. Hence, they coincide because on the other hand, their difference is an integer number. Indeed, the exponential of each of them is the linear part of the holonomies of
$\br{S}$ and $\br S'$, which  are analytically conjugated by the biholomorphism $\Psi_{|T}$.

\begin{obs}\label{Ntop}
Notice that the above proof only uses the transverse holomorphy of $\Psi$  on a neighborhood $V$ of the strict transform of each punctured separatrix of $\F$.
Hence we deduce that
any topological {conjugation} which is transversely holomorphic in a neighborhood of each punctured nodal separatrix is a $\mc N$-topological {conjugation}.
\end{obs}

\subsection{Invariance of all the Camacho-Sad indices}\label{invtousCSpreuve}
The proof is based in the Camacho-Sad index formula, which claims that the auto-intersection number of an irreducible component $D$ of the exceptional divisor equals the sum of the Camacho-Sad indices along $D$, at the singular points of the foliation lying on $D$.
We consider  filtrations of the exceptional divisors  $\mc E_S:=E^{-1}_{S}(0)$ and $\mc E_{S'}:=E^{-1}_{S'}(0)$,
$$\mc E_0:=\mc E_S\supset\mc E_1\supset \mc E_2\supset \cdots\quad\hbox{ and }\quad \mc E'_0:=\mc E_S'\supset\mc E'_1\supset\mc E'_2\supset\cdots\,,$$
defined by induction in the following way: $\overline{\mc E_{j-1}\setminus\mc E_j}$
is the union of the components  $D$ of $\mc E_{j-1}$ having  valence 1 in $\mc E_{j-1}$.
Since the dual graphs of these divisors are trees, we eventually obtain the empty set. Clearly $G(\mc E_j)=\mc E'_j$, for all  $j$.
In order to obtain the equalities {(\ref{asts})}
 %$(\ast)_s$
in Theorem~\ref{invariance} for every singular point $s\in \Sing(\underline{\F})$ and each irreducible component $D$ of $\mc D_S$, it suffices to show the following assertion  for all
$j\geq 1$:
\begin{enumerate}
  \item[$(\divideontimes)_j$] \it  the equality {(\ref{asts})}
  % $(\ast)_s$
holds at every point $s\in (\Sing(\mc D_S)\setminus\Sing(\mc E_j))$, for each irreducible component $D$ containing $s$.
\end{enumerate}
We conclude by noting that the index formula provides the implication $(\divideontimes)_j\Rightarrow(\divideontimes)_{j+1}$ and that $(\divideontimes)_0$
express the invariance of the Camacho-Sad indices of the separatrix curve, proved in Section~\ref{invcssepartr}.

\section{Proof of  Classification Theorem~\ref{classification}}\label{8}
We keep the notations (\ref{notast}), (\ref{notationstar}), (\ref{tubemilnorB}), Conventions~\ref{convtubemiln} and we assume the hypotheses of Theorem~\ref{classification}.
We shall construct a global $\mc N$-transversely holomorphic {conjugation} between $\un\F$ and $\un\F'$ inducing a realization of  $(\mf g, {h})$ outside some $\mc N$-separators of $\F$ and $\F'$ in ${\mb B}$, satisfying Assertions  (a), (b) and (c) of Theorem~\ref{classification}.
We proceed by induction, by constructing the desired homeomorphism step by step, over ``elementary pieces'' of an appropriated decomposition  of a neighborhood of the total divisor  $\mc D_S $ in $\mc B_S$ that we will describe in Section~\ref{pieceselem}. These elementary pieces are associated to each singular point $s\in\mr{Sing}(\mc D_{S})$ and to each irreducible component $D$ of $\mc D_{S}$ and they are  denoted by $K_{s}$ and $K_{D}$ respectively.

\subsection{Description of the induction }\label{descripinduc} The construction of the elementary pieces is done in  Section~\ref{pieceselem}.
Extension Lemma~\ref{lemext}  is the key tool which allows us to make the inductive step and  to begin the process.
Given a realization  $(\phi ,\wt\phi ,h)$ of $(\mf g,{h})$ over a fiber $T$ of the Hopf fibration contained in a boundary component of an elementary piece  $K$, this lemma gives a simple topological condition
{(\ref{bullet})}
%$({\sstyle \bullet})$
which allows to extend $(\phi ,\wt\phi ,h)$ to a realization defined over the whole piece. In addition we have that
\begin{itemize}
  \item[1.]
  the restriction of this extension to every Hopf fiber contained in  $\partial K$ also satisfies the condition {(\ref{bullet})}, %$({\sstyle \bullet})$,
  \item[2.] when  $T$ is contained in the intersection of two adjacent elementary pieces, then the realizations given by Lemma~\ref{lemext} over each of these pieces coincide over their common intersection.
\end{itemize}
Thus, if $\mr{Node}(\un{\F})=\emptyset$ and we have a realization over an elementary piece $K^{0}$, we can extend it step by step over a whole neighborhood of $\mc D_{S}$ in $\mc B_{S}$.
In order to achieve the proof in this case, it suffices to be sure that we can apply Lemma~\ref{lemext} in the context where  $D$  is the strict transform of the irreducible component of $S$ meeting $\Sigma$, $T=E^{-1}(\Sigma )$ and letting $K^0$ be the elementary piece associated to $D$. In Section~\ref{deminduction} we will prove the existence of a realization  $(g_1,\wt g_1,h)$ of the {conjugation}  $(\mf g, {h})$ satisfying the condition {(\ref{bullet})} %$(\bullet)_{\wt g_1,\,\wt \psi }$
in Lemma~\ref{lemext} in this context. This will achieve the proof of the theorem in the case $\mr{Node}(\un{\F})=\emptyset$.

If $\mr{Node}(\un{\F})\neq\emptyset$,  we
begin this induction process  in each connected component of $\Sigma$ and we stop it when it would require to make an extension to an elementary piece containing a  singular point $s$ belonging to $\mathrm{Node}(\un\F)\cap\mathrm{Sing}(\mc E_\F)$.
For a  nodal singularity $s$ (resp. $s'$) belonging to the strict transform $\br{\mc S}$ (resp. $\br{\mc S}'$) of a nodal separatrix of $\F$ (resp. $\F'$), Extension Lemma~\ref{lemext} provides a foliated homeomorphism $K_{s}\to K_{s'}$ which can be easily extended to the adjacent elementary pieces $K_{\br{\mc S}}\to K_{\br{\mc S'}}$ by using the product structures of $\un\F$ and $\un\F'$ in $K_{\br{\mc S}}\cong (K_{s}\cap K_{\br{\mc S}})\times[0,1]$ and $K_{\br{\mc S'}}\cong  (K_{s'}\cap K_{\br{\mc S}'})\times[0,1]$.
To complete the process in this case, it suffices to glue the realizations obtained in this way, by constructing in Section~\ref{FinalRemarks} suitable foliated homeomorphisms defined
on the elementary pieces associated to the singular points $s\in\mathrm{Node}(\F)\cap\mathrm{Sing}(\mc E_\F)$.

\subsection{Elementary pieces}\label{pieceselem}
{Recall that we have fixed
%geometric representation $(g,\wt g,h)$ of the
a $\mc N$-conjugation $(\mf g,h)$ in Theorem~\ref{classification}.
By Corollary~\ref{reprfibree} there exists an excellent homeomorphism germ $(\mb B,S)\to (\mb B,S')$ representing the marking determined by $\mf g$.
Let $G:(\mc B,S)\to (\mc B',S')$ be its lifting.}
 We fix two Milnor tubes, cf. (\ref{tubemilnorB}), $\mc T_{\varepsilon }$ for $S$ and  $\mc T'_{\varepsilon '}$ for $S'$,  where $\varepsilon , \varepsilon '>0$ are chosen small enough so that $G(\mc T_{\varepsilon })\subset \mc T'_{\varepsilon '}$ and each real hypersurface $\{|x_s|=1\}$ and $\{|y_s|=1\}$, $s\in \Sing(\underline{\F})$, as well as $\{|x_{s'}|=1\}$ and $\{|y_{s'}|=1\}$, $s'\in \Sing(\underline{\F}')$, separates the tube in two connected components and intersects transversely the boundary in a $2$-torus.
We denote
$$\mc H :=\bigcup_{s\in \Sing(\underline{\F})}\{|x_s|=1\}\cup\{|y_s|=1\}
,\;
\mc H':=\bigcup_{s'\in\Sing(\underline{\F}')}\{|x_{s'}|=1\}\cup\{|y_{s'}|=1\}\,.$$
We call \emph{elementary piece} of  $\mc T_{\varepsilon }$ (resp. $\mc T_{\varepsilon '}'$),
every intersection $K:= \mc K\cap %\inte
{{\mc T}_{\varepsilon }}$ (resp. $K':= \mc K'\cap {%\inte
{\mc T}}{}'_{\varepsilon '}$), where $\mc K$ (resp. $\mc K'$), is the adherence of a connected component of $\mc T_{\varepsilon }\setminus \mc H$ (resp. $\mc T'_{\varepsilon '}\setminus \mc H'$). For each elementary piece $K$ (resp. $K'$), one and only one of the following assertions holds:
\begin{itemize}
  \item $K$ (resp. $K'$), contains a  (unique) point $s$ of   $\Sing(\underline{\F})$ (resp. $s'\in\Sing(\underline{\F}')$), and it is contained in the domain  $\mc W_s$ (resp. $\mc W_{s'}$) of the coordinate chart $(x_s, y_s)$ (resp. $(x_{s'}, y_{s'})$);
  \item $K$ (resp. $K'$), contains a compact set $D^\circ := \overline{D\setminus\bigcup_{s}{\mc W}_s}$, where $D$ denotes an irreducible component of  $\mc D_{S}$ (resp. $\mc D_{S'}$), and  $s$ ranges the set of singular points of $\underline{\F}$ (resp.  $\underline{\F}'$); in addition, if $\varepsilon,\varepsilon'>0$ are small enough then
the restriction of the fibration
   $\pi _D$ to $K\cap\pi _D^{-1}(\partial D^\circ)$ is still a disk fibration;
\end{itemize}
In the first case the elementary piece will be denoted by $K_s$ (resp. $K_{s'}$), and in the second case it will be denoted by $K_D$.
The intersection of two different elementary pieces
is either empty or a solid $3$-torus.

\subsection{Extension of realizations}
\label{subsctext} {By Corollary~\ref{reprfibree} there is}
 a geometric representation  $(g, \wt g,h)$ of the  $\mc N$-{conjugation} $(\mf g, {h})$, with $g$ excellent. We assume that $\varepsilon >0$ is small enough so that  $G:=E_{S'}^{-1}\circ g\circ E_{S|\mc T_{\varepsilon }}$ is defined on $\mc T_{\varepsilon }$ into  $\mc T_{\varepsilon '}'$.
%Up to composing $g$ and $\wt g$ by appropriate isotopies, we can suppose that
Thanks to {properties of excellent maps}  stated in   Definition \ref{excellent}, the restriction of $G$ to the total divisor fulfills the following equalities:
$$
G(K_{{\alpha}}\cap\mc D_S) = K_{G({\alpha})}\cap \mc D_{S'}\,,\quad {\alpha}\in \Comp(\mc D_S)\sqcup\Sing(\underline{\F})\,,
$$
and $G(K_{{\alpha}})$ is an neighborhood of $K_{G(\alpha)}\cap \mc D_{S'}$ in $K_{G({\alpha})}$.\\

We consider an irreducible component $D$ of $\mc D_S$ and a Hopf fiber\linebreak $T:=\pi _D^{-1}(c)\cap \mc T_{\epsilon }$, over a point $c$ in the boundary of $D^\circ$, cf. (\ref{comptrouee}). The connected component $C$ of $\partial D^\circ$ containing $c$, is a circle that bounds a disk  $ \mc W_s\cap D$,  $s\in \Sing(\un\F)$.
The point  $c':=G(c)$ belongs to the boundary of  $D'{}^\circ$, $D':=G(D)$. We denote   $T':=\pi _{D'}^{-1}(c')$ and we assume that  a germ of biholomorphism  $\phi_{S} : (T, c)\to (T', c')$ is given, as well as a pro-germ at infinity $\wt{\phi }_{\infty}  : (\wt T^{*},\infty)\to (\wt T'^{\star},\infty)$
lifting  $\phi _S$, such that  $(\phi_{S} ,\wt\phi_{\infty} ,h)$ is a realization of $(\mf g,{h})$ over $T$ and $T'$.

\begin{convs}\label{convenrevecl}  For $B\subset \wt{\mb B}^* $ and  $B'\subset \wt {\mb B}'^{\star}$, we denote : $$\pi _0(B,\infty) := \limproj_{U\in\mf U_{\F,\Sigma }}\pi _0(B\cap\wt U^\ast)\,,\quad\pi _0(B' ,\infty) := \limproj_{U\in\mf U_{\F',\Sigma '}}\pi _0(B'\cap\wt U^\star)\,.$$
\end{convs}

\begin{lema}[of extension of realizations]\label{lemext} Let
$K$ (resp. $K'$) be one of the elementary pieces   $K_D$ or $K_s$  (resp. $K_{G(D)}$ or $K_{G(s)}$)
and let us denote  $Z:=K\cap \mc D_S$ and $Z':=K'\cap \mc D_{S'}$. We assume that $\wt\phi_{\infty}$ and the  restriction of $\wt g$ to $\wt T^\ast$ induce the same map
%$$(\bullet)_{\wt g,\,\wt \phi_{\infty} }\qquad
\begin{equation}\label{bullet}
\wt \phi _{\infty%\ssstyle \bullet
}=  \wt g%_{\ssstyle \bullet}
 : \pi _0(\wt T^{\ast},\infty) \too
\pi _0(\wt T'{}^{\star},\infty)\,,%$$
\end{equation}
and, in the case that  $K=K_s$, we also assume the following equality concerning the Camacho-Sad indices: $$\CS(\un \F, D, s)=\CS(\un \F', D', G(s))\,.$$
Then there are homeomorphisms $\Phi  : V\to V'$ and $\wt \Phi : \wt V^\ast\to\wt V'{}^\star$ such that for each component
$\br D$ of $\mc D_S$ meeting $K$, the following properties hold:
\begin{enumerate}[(a)]
\item\label{prem} $V$ (resp. $V'$), is an open neighborhood of  $\mc D_S\cap K$  in  $K$ (resp.  $\mc D_{S'}\cap K'$ in $K'$), $\wt\Phi $ lifts   $\Phi $ i.e. $\un q'\circ\wt \Phi =\Phi \circ \un q_{|\wt V^\ast}$ and $\Phi _{|Z\cap\br D^\circ}=G_{|Z\cap\br D^\circ}$;
\item\label{hopfresp}  $\Phi$
preserves the Hopf fibres over $K\cap\br D^\circ$, more precisely, $$\pi _{\br D'}\circ\Phi_{|V\cap\pi _{\br D}^{-1}(V\cap \br D^\circ)} =G\circ {\pi _{\br D}}_{|V\cap \br D^\circ},$$ with  $\br D':=G(\br D)$;
\item\label{extendephi} the germ of the restriction of $\Phi $ to $T$ coincides with  $\phi _S$ and the  pro-germ at infinity of the restriction of $\wt \Phi $ to $\wt T^\ast$ coincides with $\wt \phi _\infty$;
\end{enumerate}
moreover, except in case $K=K_s$, $s\in\mathrm{Node}(\uF)$, we have:
\begin{enumerate}[(a)]
\setcounter{enumi}{3}
\item\label{cestunereal} if we denote by $\Phi _S$ the germ of $\Phi $ along  $K\cap\mc D_S$ and by  $\wt \Phi _\infty$ the pro-germ  at infinity of $\wt \Phi $, then $(\Phi _S, \wt \Phi _\infty, h)$ is a relization of the {geometric} {conjugation}  $(\mf g, {h})$ over $V$ and $V'$, in the sense given in  Conventions~\ref{convtubemiln}; in particular  $\Phi $ is a transversely holomorphic {conjugation} between the restricted foliations  $\un{\F}_{|V}$ and $\un{\F}'_{|V'}$;
\item\label{dercondi} for each  $t\in \br D^\circ\cap K$, the restrictions of
              $\wt g$ and $\wt \Phi $ to  $T_t:=\pi _{\br D}^{-1}(t)$ induce the same map
              $$
              \wt g_{|\wt T_t^\ast%\,\ssstyle \bullet
              }= \wt \Phi _{|\wt T_t^\star%\,\ssstyle \bullet
              } : \pi _0(\wt T^{\ast}_t,\infty) \to
\pi _0(\wt T'{}^{\star}_{t'},\infty)\,, \quad T'_{t'}:=\pi _{\br D'}^{-1}(t')\,,\; t':=G(t)\,.
$$
\end{enumerate}
\noindent In addition, the homeomorphisms obtained above by extension along  $K=K_D$  and $K=K_s$, coincide over the intersection of their domains of definition.
\end{lema}

\subsection{Beginning of the induction}\label{deminduction}
We begin with a realization $(\psi,\wt\psi,h)$ over $\mc N$-collections of transversals $\Sigma$ and $\Sigma'$ of a geometric $\mc N$-{conjugation} $(\mf g,{h})$ between the monodromies of $\F$ and $\F'$ geometrically represented by $(g,\wt g,h)$.
For each connected component $\br\Sigma$ of %the $\mc N$-collection of transversals
$\Sigma$ we consider the irreducible component $\br S$ of $S$ meeting $\br\Sigma$. From Condition ($\sigma$) in the Introduction and Condition (2) in  Theorem~\ref{classification}, we deduce that $\br S$ (resp. $\br S':=g(\br S)$) is a non-nodal separatrix of $\F$ (resp. $\F'$).
Denote by $\br{ \mc S}$ and $s$ (resp. $\br{\mc S}'$ and $s'$) the strict transform of $\br S$ (resp. $\br S'$) and its attaching point in the exceptional divisor $\mc E_{S}$ (resp. $\mc E_{S'}$).
Denote by $\br\Sigma'\supset\psi(\br\Sigma)$ the connected component of $\Sigma'$ meeting $\br S'$.
Thanks to the invariance by $\F$ or $\F'$-isotopies stated in Proposition~\ref{lemmeisotopie}, we can assume that  $\br \Sigma $ (resp. $\br\Sigma '$) is contained in the image $E_S(T_0)$ (resp. $E_{S'}(T'_0)$) of a Hopf fiber  $T_0$ (resp. $T'_0$) in the intersection of the elementary pieces $K_{\br{\mc S}}$ and $K_{s}$ (resp. $K_{\br{\mc S}'}$ and $K_{s'}$).
Analogously, up to composing the geometric representation $(g,\wt g,h)$ of the $\mc N$-{conjugation} $(\mf g,{h})$ with
$S$ or $S'$-isotopic homeomorphisms, cf. Remark~\ref{invparisotopie}, we can assume the equality of germs $(g(\br\Sigma),g(\br\Sigma\cap \br S))=(\br\Sigma',\br\Sigma'\cap \br S')$.
Finally, we also assume that
${x_s=0}$ (resp. ${x_{s'}=0}$) is a reduced local equation of  $\br{\mc S}$ (resp.  $\br{\mc S}'$) and we denote:
$$
\mb T:=E_S(\{|{y_s}|=1, {x_s}=0\})\quad\hbox{and}\quad\mb T':=E_{S'}(\{|{y_{s'}}|=1, {x_{s'}}=0\})\,.
$$
In order to simplify the notations in the rest of this section we will write $\wt\Sigma^{*}$ and $\wt\Sigma'^{\star}$ instead of $\wt{\br\Sigma^{*}}$ and $\wt{\br\Sigma'^{\star}}$.

\indent Now  notice that  $\psi$ and $\wt \psi$ are necessary  holomorphic. This follow from the commutativity of the  diagram {($\star$)} %$(\star)_{\wt \psi}$
in Definition~\ref{defconj}, using Proposition~\ref{monoprogermcan}. Then, in order to begin with the inductive process described in Section~\ref{descripinduc}, we will construct a geometric representation  $(g_1,\wt g_1,h)$ of the $\mc N$-{conjugation} $(\mf g, {h})$ which will satisfy the topological condition
(\ref{bullet})
%$({\sstyle\bullet})$
required by the Extension Lemma~\ref{lemext} in the case $D=\br{\mc S}$,    $T=T_0$ and $(\phi_{S} , \wt \phi_{\infty} , h)=(\psi , \wt \psi , h)$.
This condition is equivalent to the following equality:
\begin{equation}\label{bulletbis}
%$$(\bullet)_{\wt g_1,\,\wt \psi }\qquad
\wt \psi %_{\ssstyle \bullet}
=
\wt g_{1%\ssstyle\, \bullet
} : \pi _0(\wt \Sigma ^{\ast},\infty) \too
\pi _0(\wt \Sigma '{}^{\star},\infty)\,.%$$
\end{equation}
%An easy computation similar to (\ref{calculsimple}) shows that $(\bullet)_{\wt g_1,\,\wt \psi }$
{We can see that (\ref{bulletbis})}
holds if the equality $\wt \psi %_{\ssstyle \bullet}
(\sigma_{0} )=
\wt g_{1%\ssstyle \bullet
}(\sigma_{0} )$ is satisfied by a fixed element $\sigma_{0}$ of $\pi _0(\wt \Sigma ^{\ast},\infty)$.
Thanks to the Lemma~\ref{lemmetop} below, $\wt g%_{\ssstyle\bullet}
(\sigma _0)$ and $\wt\psi %_{\ssstyle\bullet}
(\sigma _0)$ belong to the same fiber of the map $\iota %_{\ssstyle\bullet} 
: \pi _0(\wt \Sigma '{}^{\star},\infty) \to \pi _0(\wt{\mb T}'{}^\star,\infty)$
induced by the inclusion $\iota  : \wt \Sigma '^\star\hookrightarrow \wt{\mb T}'{}^\star$.
{This property will allow us to construct a homeomorphism   $\Theta : {\mb B}'\iso{\mb B}'$ which is $S'$-isotopic to the identity and whose lift $\wt \Theta $ satisfies  $\wt \Theta %_{\ssstyle\bullet}
(\wt g%_{\ssstyle\bullet}
(\sigma _0))=\wt\psi %_{\ssstyle\bullet}
(\sigma _0)$. We achieve the beginning of the inductive process by taking $g_1=\Theta \circ g$ and $\wt g_1:=\wt \Theta \circ\wt g$.
The homeomorphism $\Theta$ is a sort of foliated Dehn twist that we will describe before proving Lemma~\ref{lemmetop} {below,} which is needed to construct it.}\\

{On the disk $\br{\mc S}'\cap\mc W'_{s'}$, we consider the real vector field   $\vartheta$ whose flow is  $(t, y_{s'})\mapsto e^{2i\pi t}y_{s'}$. We fix a smooth function   $u$ with support contained in $\{x_{s'}=0, \varsigma\leq |y_{s'}|\leq 1\}$, $0<\varsigma<1$ taking the value 1 on the circle $C':= \{x_{s'}=0,\ |y_{s'}|=1\}$ %$=K'\cap \br D'^{\circ}$
and we denote by   $Y$ the vector field on $V'$   tangent to   $\un\F'$ and projecting over $u\vartheta$ by $\pi _{\br{\mc S}'}$. The flow
$\Upsilon_t$ of $Y$ is defined on a open neighborhood $U_{I}'$ of $\br{\mc S}'$ in  $\mc W_{s'}$, once we fix an interval $I\subset\mb R$ where we allow the time $t$ to vary.
Hence, we can lift the flow $\Upsilon_{t}$ to a unique map $\wt\Upsilon_t : \wt U_I'^\star\to \wt V'{}^\star$ being the identity on  $|y_{s'}\circ\un q'|\leq \varsigma$, and defining consequently a germ at infinity  $\wt\Upsilon_{t\infty} : (\wt V'{}^\star,\infty)\to(\wt V'{}^\star,\infty)$.
Clearly the germ  $\wt\Upsilon_{n\infty}$ fibers over  $ C'$, i.e. ${\pi _{\br{\mc S}'}\circ\un q'\circ\wt\Upsilon_{n\infty}}_{|\wt{\mb T}'{}^\star}={\pi _{\br{\mc S}'}\circ\un q'_{\infty}}_{|\wt{\mb T}'{}^\star}$, for each $n\in\mb Z$.
It defines a deck transformation $\wt\Upsilon_{n%\ssstyle\bullet
} :\Pi '\iso\Pi '$ of  the natural covering
$$
\rho' :\Pi' \to C'\,,\quad \rho '{}^{-1}(p):=\pi _0(\wt T'_p{}^\star,\infty)\,,\quad T'_p:=\pi _{\br{\mc S}'}^{-1}(p)\,.
$$
We put $\Theta :=\Upsilon_{n_0}$ and $\wt\Theta :=\wt\Upsilon_{n_0}$,
choosing the integer $n_0$ in the following way. First,
 we fix a point $a$ in the circle $C:=\{x_s=0, |y_s|=1\}$ %$=K\cap\br D^\circ$
 and an element  $\nu$ of $\pi _0(\wt T_{a}^{\star},\infty)$, $T_a:=\pi _{\br{\mc S}}^{-1}(a)$. We consider the natural covering  $\rho :\Pi \to C$ of $C$  with fibers $\rho ^{-1}(p):=\pi _0(\wt T^\ast_p,\infty)$ and  the following two covering morphisms over $G_{|C}$
      $$
%(\blacklozenge )_{\Lambda }
%\quad
  \begin{array}{ccc}
    \Pi  & \stackrel{\Lambda }{\too} &  \Pi '\\
    \rho \downarrow\phantom{\rho } & \circlearrowleft & \phantom{\rho '}\downarrow\rho ' \\
    C & \stackrel{G_{|C}}{\too} & C'
  \end{array}\,,\qquad \quad\Lambda = \wt \psi %_{\ssstyle\bullet}
  \quad\hbox{and}\quad \Lambda  =\wt g%_{\ssstyle\bullet}
  \,,
  $$
  defined by $\wt \psi $ and $\wt g$.
Using Lemma~\ref{lemmetop},
$\sigma ':=\wt \psi %_{\ssstyle\bullet}
(\sigma_{0} )$  and $\sigma '':=\wt g%_{\ssstyle\bullet}
(\sigma_{0})$ belong to the same fiber of the map $\iota %_{\ssstyle \bullet} 
: \pi _0(\wt T_{a'}'{}^{\star},\infty)\to\pi _0(\wt{\mb T}'{}^\star,\infty)$, induced by the inclusion map  $\iota  : \wt T'_{a'}{}^{\star}\hookrightarrow \wt{\mb T}'{}^\star$. On the other hand, the action $$
\mb Z\times \pi _0(\wt T'_{a'}{}^{\star},\infty) \too \pi _0(\wt T'_{a'}{}^{\star},\infty)\,\qquad(n,\sigma )\mapsto \wt\Upsilon_{n%\ssstyle \bullet
}(\sigma )\,,
$$
of $\mb Z$ on the  fiber  of $\rho '$ at the point $a':=\psi (a)=G(a)\in C'$ coincide with the action of  $\pi _1(C', a')\simeq \mb Z$ induced by the covering $\rho'$ on this fiber.
The orbits of that action correspond to the fibers of  $\iota %_{\ssstyle \bullet}
$.
We choose $n_0$ to be the unique integer number such that $\Upsilon_{n_0%\ssstyle\bullet
}(\sigma ')=\sigma ''$.}

%{\color{red}By using the direct image  $E_{S'\ast}Y$ of the vector field  $Y$ constructed in the part~A of the fourth step of the proof of Extension Lemma~\ref{lemext} for $K=K_s$, it is easy to obtain a homeomorphism   $\xi : {\mb B}'\iso{\mb B}'$ which is $S'$-isotopic to the identity and whose lift $\wt \xi $ satisfies  $\wt \xi _{\ssstyle\bullet}(\wt g_{\ssstyle\bullet}(\sigma _0))=\wt\psi _{\ssstyle\bullet}(\sigma _0)$.} It suffices to put $g_1=\xi \circ g$ and $\wt g_1:=\wt \xi \circ\wt g$, in order to conclude the proof.

\begin{lema}\label{lemmetop}
The following diagram is commutative:
$$
\begin{array}{ccc}
    \pi_{0}(\wt\Sigma^{*},\infty)  & \stackrel{\wt\psi%_{\bullet}
    }{\too} &  \pi_{0}(\wt\Sigma'^{\star},\infty)\\
    \imath%_{\bullet} 
    \downarrow\phantom{\rho } & \circlearrowleft & \phantom{\rho '}\downarrow\imath%_{\bullet} 
    \\
    \pi_{0}(\wt{\mb T}^*,\infty) & \stackrel{\wt{g%_{\bullet}
    }}{\too} & \pi_{0}(\wt{\mb T}'^\star,\infty)
  \end{array}
$$
 \end{lema}
\begin{dem}
Let $\wt\Sigma_{0}^*$ be a connected component of $\wt \Sigma^*$ and consider the subgroup $\Gamma_{0}\subset\Gamma_{\infty}$ consisting of those elements $\varphi\in\Gamma_{\infty}$ such that $\varphi(\wt\Sigma_{0}^*)=\wt\Sigma_{0}^*$. Let us
denote  $\wt\Sigma_{0}'^{\star}:=\wt g(\wt\Sigma^*_{0})$ and $\wt\Sigma_{1}'^{\star}:=\wt\psi(\wt\Sigma_{0}^*)$ and for $i=0,1$ consider also the subgroups $\Gamma'_{i}$ of $\Gamma'_{\infty}$ consisting of those elements $\varphi$ such that $\varphi(\wt\Sigma'^{\star}_{i})=\wt\Sigma'^{\star}_{i}$. There exists $\gamma\in\Gamma'_{\infty}$ such that $\wt\Sigma'^{\star}_{1}=\gamma(\wt\Sigma'^{\star}_{0})$.
A straightforward computation shows that $\Gamma'_{1}=\gamma\Gamma'_{0}\gamma^{-1}$. On the other hand, if $\varphi\in\Gamma_{0}$ then $\wt g_{*}(\varphi)\in\Gamma_{0}'$. Using
 {($\star\star$)} %$(\star\star)_{\wt \psi }$
we also deduce that $\wt g_{*}(\varphi)\in\Gamma'_{1}$ for all $\varphi\in\Gamma_{0}$.
Consequently\linebreak $\Gamma_{0}'\cap\gamma\Gamma_{0}'\gamma^{-1}\neq\{1\}$.
We fix a point $\wt c\in\wt\Sigma^\star_{0}$ and we identify $\Gamma_{\infty}'$ to $\pi_{1}(\mb B'^{\star},c)$, $c:=q'(\wt c)$, by means of the isomorphism
 $\chi : \varphi  \mapsto q'\circ\dot{\gamma} _{\varphi }$,  $\varphi \in \Gamma_\infty'$, where $\gamma_\varphi $ denotes a path in  $\wt{\mb B}'{}^\star$, having endpoints  $\wt c$ and $\varphi (\wt c)$. Clearly, $\chi(\Gamma_{0}')=\pi_{1}(\Sigma_{0}'^{\ast},c)=\langle\mf m_{c}\rangle\subset\langle \mf m_{c},\mf p_{c}\rangle=\pi_{1}(\mb T'^{\star},c)$, cf. Definition~\ref{periph}.
Since $H_{1}(\mb B'^{\star},\mb Z)$ is torsion free, if $\mf m_{c}^{\alpha}=\chi(\gamma) \mf m_{c}^{\beta}\chi(\gamma)^{-1}$ then $\alpha=\beta$ and consequently $[\mf m_{c}^{\alpha},\chi(\gamma)]=1$. It remains to see that if $\alpha\neq 0$ then $\chi(\gamma)$ belongs to $\pi_{1}(\mb T'^{\star},c)$. Indeed, if this is  the case then $\wt\Sigma_{0}'^{\star}$ and $\wt\Sigma_{1}'^{\star}$ are contained in the same connected component of $\wt{\mb T}'^{\star}$.
That $\chi(\gamma)\in\pi_{1}(\mb T'^{\star},c)$ follows easily from the Sub-Lemma~\ref{amal} below. Indeed, we apply it first to $A=\pi_{1}(W,c)$, where $W\subset\mb B'^{\star}$ is the JSJ block containing $\Sigma'^{\star}$. We deduce that $\chi(\gamma)\in A$ and we pass to the quotient $\bar A=A/\langle \mf p_{c}\rangle$ by its center $\langle \mf p_{c}\rangle$. We express  $\bar A=\mb Z\overline{\mf m_{c}}\ast G$ and we apply again Sub-Lemma~\ref{amal} in order to have that $\overline{\chi(\gamma)}\in\langle \overline{\mf m_{c}}\rangle$. Hence $\chi(\gamma)\in\langle \mf m_{c},\mf p_{c}\rangle$.
\end{dem}

\begin{sublema}\label{amal}
Let $\Gamma=A\ast_{C}B$ be an amalgamated product of groups $A$ and $B$ over a common subgroup $C$. Fix an element $\alpha\in A\setminus\bigcup\limits_{a\in A}a C a^{-1}$. If $\gamma\in\Gamma$ commutes with $\alpha$ then $\gamma\in A$.
\end{sublema}

\begin{dem2}{of the sub-lemma}
For each element $\gamma$ in $\Gamma$, there is a unique natural number $n\ge 1$ such that $\gamma$ can be written in only one of the following ways
\begin{enumerate}[(i)]
\item
$a_{1}b_{1}a_{2}\cdots a_{n}b_{n}$,
\item $b_{1}a_{1}b_{2}\cdots b_{n}a_{n}$,
\item $b_{1}a_{1}b_{2}\cdots a_{n-1} b_{n}$,
\item $a_{1}b_{1}a_{2}\cdots b_{n-1}a_{n}$,
\end{enumerate}
provided that $a_{i}\in A\setminus C$ and $b_{i}\in B\setminus C$. If $[\alpha,\gamma]=1$ then necessarily $\gamma$ is of type (iv). Since $a_{1}^{-1}\alpha a_{1}b_{1}a_{2}\cdots b_{n-1}a_{n}\alpha^{-1}=b_{1}a_{2}\cdots b_{n-1}a_{n}$ and $a_{1}^{-1}\alpha a_{1}\notin C$ it follows that $n=1$.
\end{dem2}

\subsection{Gluing at the nodal singularities}\label{FinalRemarks}

For $s\in \mathrm{Node}(\F)\cap\mathrm{Sing}(\mc E_S)$ both germs of $\uF$ at $s$ and of $\uF'$ at $s':=G(s)$ are holomorphicaly conjugated to the germ, at the origin in $\C^2$, of the same linear complex foliation $\mc L$ having a  multivaluated first integral $xy^\lambda $, where $\lambda \in \mb R_{<0}\setminus\mb Q$ and $1/\lambda$ are the common Camacho-Sad indices at these singularities. This foliation $\mc L$ has a conical structure over the induced (real) foliation by lines $\mc L_{\ssstyle \mb R}$  on the sphere $
\partial\mb P$ with $\mb P=\mb D_{1}\times\mb D_{1}$.
Outside the knots defined by the coordinate axis, the leaves of  $\mc L_{\ssstyle \mb R}$ are contained in the  tori constituted by the level sets of  $|x||y|^{\lambda }$ in that sphere. Restricted to any torus, $\mc L_{\ssstyle \mb R}$  can be seen as the suspension of the rotation of angle $2\pi \lambda $. The following property of $\mc L$ allow us to easily build a topological {conjugation} between $\uF$ and $\uF'$, which is defined on  neighborhoods of the divisors in the elementary pieces associated to $s$ and $s'$ and which coincide with the {conjugations} previously obtained on the adjacent elementary  pieces.
\begin{itemize}\it
  \item let $\phi(x,y):= (\un\phi(x,y),y) $ be a germ along the circle $\mc C:=\{0\}\times\partial\mb D_{1}$ of a biholomorphism defined on a neighborhood of $\mc C$ in the solid torus $\mb D_{1}\times\partial  \mb D_{1}$, which preserves the foliation $\mc L_{\ssstyle \mb R}$. Then there exist positive real numbers $\varepsilon_1<\varepsilon_2<1$  and a global homeomorphism $\Phi : \partial\mb P\iso\partial\mb P$ which coincide with $\phi$ for $|x||y|^\lambda \leq\varepsilon_1$, and is the identity map for  $|x||y|^\lambda \geq\varepsilon_2$.
\end{itemize}
The proof is based on the linearity of  the maps $\un\phi(x, \cdot)$ resulting of theirs commutations with the holonomy map of $\mc L$ along the circle $\mc C$ -which is a non periodic rotation. We can construct $\Phi$ by a suitable interpolation, along the family of torus $(\mb T_r:=\partial \mb D_{r}\times \partial\mb D_{1})_{\varepsilon_1\leq r\leq \varepsilon_2}$ between the restriction of $\phi $ to $\mb T_{\varepsilon_1}$ and the identity map on $\mb T_{\varepsilon_2}$. %We leave  to the reader to complete the details.

\subsection{{Proof of  Extension Lemma~\ref{lemext} for $K=K_D$}}
Thanks to Theorem~\ref{monodimplholo}, $\phi $ conjugates the holonomy representation of  $\un\F$ along
$D^\circ$ and that of $\un\F'$ along $D'{}^\circ$. Therefore, by using the classical lifting path method we have the following:
\begin{enumerate}
  \item [($\blacktriangle$)] \it there are fundamental systems  $(W_{k})_{k\in\mb N}$ and $(W'_{k})_{k\in\mb N}$ of neighborhoods of  $D^\circ$ and $D'{}^\circ$ in  $K$ and $K'$ and there are homeomorphisms $\Phi _k$  from $W_k$ onto $W'_k$, such that:
\begin{enumerate}[1)]
\item for all $k$, $l$, ${\Phi _k}$ and $\Phi _l$
coincide over $W_k\cap W_l$;
\item the intersection of $W_k$ (resp. $W'_k$)  with the fibers of  $\pi _D$ (resp. $\pi _{D'}$), are conformal disks  and hence  $W_k^\ast$ and $K^\ast$ (resp. $W'{}^\star_k$ and $K^\star$) are homotopic;
\item $\Sat_{\un\F}(T\cap W_k,W_k)=W_k$ and  $\Sat_{\un\F'}(T'\cap W'_k,W'_k)=W'_k$;
      \item $\pi _{D'}\circ \Phi  _k = G\circ {\pi _{D}}_{|W_k}$;
     \item the restrictions of $\phi $ and $\Phi _k$ to $W_k\cap T$ coincide;
      \item $\Phi _k$ conjugates $\un\F_{|W_k}$ to $\un\F'_{|W'_k}$.
      \end{enumerate}
\end{enumerate}
\noindent Denote
$
V:=W_0\,, V':=W'_0\,,\Phi :=\Phi _0
$, $\Delta :=V\cap T$ and $\Delta ':=V\cap T'$.
From Property  $(\blacktriangle)$ 2) above follows  that the restriction of $\un q$ to each connected component $\wt V^{\ast\beta }$ of $\wt V^{\ast}$,  ${\beta \in \pi _{0}(\wt V^{\ast})}$ (resp.  $\wt V'{}^{\star\beta' }$ of $\wt V'{}^{\star}$, ${\beta' \in \pi _{0}(\wt V'{}^{\ast})}$) is a universal covering of  $V^{\ast}$ (resp. $V'{}^{\star}$); the same is true for the connected components  $\wt \Delta ^{\ast\alpha }$ of $\wt\Delta ^{\ast}$ over $\Delta ^\ast$, ${\alpha \in \pi _{0}(\wt \Delta ^{\ast}\cap\wt V^{\ast\beta })}$ and for those  $\wt \Delta '{}^{\star\alpha '}$ of $\wt\Delta '{}^{\star}$ over $\Delta '{}^\star$, ${\alpha ' \in \pi _{0}(\wt\Delta '{}^{\star})}$. We fix  $\beta \in \pi _0(\wt V^\ast)$, $\alpha \in \pi _0(\wt\Delta ^\ast\cap\wt V^{\ast\beta })$ and we denote by  $\wt \Phi^\beta :\wt V^{\ast\beta }\to\wt V'{}^{\star\beta '}$ the unique homeomorphism that lifts $\Phi $ over $\wt V^{\ast\beta }$  coinciding with  $\wt\phi $ on $\wt \Delta ^{\ast\alpha }$, where $\beta '$ corresponds to the connected component of $\wt V'{}^{\ast}$ containing $\wt\phi (\wt\Delta ^{\ast\alpha })$. We will prove the following assertions:
\begin{enumerate}[(i)]
  \item\label{i}  $\wt\Phi ^\beta$ does not depend on the choice of  $\alpha \in \pi _0(\wt \Delta ^{\ast}\cap \wt V^{\ast\beta })$;
  \item\label{ii}  the homeomorphism $\wt\Phi : \wt V^\ast\to\wt V'{}^\star$ defined by requiring that its restriction to $\wt V^{\ast\lambda }$, $\lambda \in \pi _0(\wt V^{\ast})$, coincides with $\wt\Phi ^{\lambda}$, satisfies Conditions  $(\ref{prem})$-$(\ref{dercondi})$ of Lemma~\ref{lemext}.
\end{enumerate}
\subsubsection{Proof of  Assertion \hbox{\rm (\ref{i})}}
Using the unicity of the lifting of $\Phi $ over the covering $\wt V^{\ast\beta }$, it suffices
to show that for each ${\kappa \in \pi _{0}(\wt \Delta ^{\ast}\cap\wt V^{\ast\beta })}$ the maps $\wt\Phi ^{\beta }$ and $\wt\phi$ coincide at one particular point of $\wt\Delta ^{\ast\kappa }$. In fact, we will prove the equality of the germs at infinity of their restrictions to $\wt\Delta ^{\ast\kappa }$.

To do that, it suffices to see that $\wt\Phi^{\beta }(\wt\Delta ^{\ast\kappa })$ and $\wt\phi(\wt\Delta ^{\ast\kappa })$ are contained in the same connected component  $\wt\Delta ^{\ast\kappa '}$ of $\wt \Delta '{}^{\ast}\cap\wt V'{}^{\star\beta '}$. Indeed,   $\wt\Phi ^\beta $ conjugates $\wt{\un\F}_{|\wt V^{\ast\beta }}$ to $\wt{\un\F}'_{|\wt V'{}^{\beta '}}$ and it factorizes through the leaf spaces of these foliations. Hence  the following  diagram is commutative:
$$
%(\star)_{\Upsilon,\;\kappa }\hspace{3em}
\begin{array}{ccc}
(\wt\Delta ^{\ast\kappa }, \infty) & \stackrel{\Upsilon}{\longrightarrow} & (\wt\Delta '{} ^{\ast\kappa '}, \infty)\\
{\sstyle\tau _{\ssstyle \wt\Delta ^{\ast\kappa }}}\downarrow\phantom {\tau
_{\ssstyle \wt\Delta ^{\ast\kappa '}}} & {\sstyle\circlearrowright}&\phantom{\scriptstyle\tau _{\wt\Delta '{}^{\ast\kappa '}}}\downarrow
{\scriptstyle\tau _{\wt\Delta '{}^{\ast\kappa '}}} \\
\QF & \stackrel{h}{\longrightarrow} & \QFpr
\end{array}\,,
\phantom{(\star)_{\Upsilon,\;\kappa }\hspace{3em}}
$$
where $\Upsilon$ denotes the germ at infinity $\wt\Phi ^\beta _{|\wt\Delta ^{\ast\kappa }\,\infty}$ of the  restriction of  $\wt\Phi ^\beta$ to $\wt\Delta ^{\ast\kappa }$, the vertical arrows are the canonical pro-germs defined in Section~\ref{progermcan}. On the other hand, this diagram still commutes when we take as  $\Upsilon$ the germ  at infinity $\wt \phi _{|\wt\Delta ^{\ast\kappa }\,\infty}$ of the restriction of $\wt\phi $ to $\wt\Delta ^{\ast\kappa }$, because  $(\phi _S,\wt \phi_\infty ,h)$ is a realization of $(\mf g,{h})$ over $T$. The equality $\wt\Phi ^\beta _{|\wt\Delta ^{\ast\kappa }\,\infty}=\wt \phi _{|\wt\Delta ^{\ast\kappa }\,\infty}$ follows from  Proposition~\ref{monoprogermcan} which asserts that  $\tau _{ \wt\Delta ^{\ast\kappa }}$ et $\tau _{\wt\Delta '{}^{\ast\kappa '}}$ are  $\proan$-monomorphisms.

Consider now the natural connected coverings  $$\chi_{\beta } : \bigsqcup_{t\in D^\circ}\pi_0\left(\wt\Delta ^{\ast}_t\cap \wt V^{\ast\beta },\infty\right)=:  \Pi _0^{\beta }\too D^\circ\,,\quad\Delta _t:=\pi _D^{-1}(t)\cap V\,,$$
$$\hbox{and}\quad \chi'_{\beta '} : \bigsqcup_{t\in D'{}^\circ}\pi_0\left(\wt\Delta'{} ^{\ast}_t\cap \wt V'{}^{\ast\beta '},\infty\right)=:  \Pi _0'{}^{\beta '}\too D'{}^\circ\,,\quad \Delta '_t:=\pi _{D'}^{-1}(t)\cap V'\,.$$
The maps $\wt g$ and $\wt \Phi ^{\beta }$ send each connected component of  $\wt\Delta ^{\ast}_t\cap \wt V^{\ast\beta }$ onto a connected component of $\wt\Delta'{} ^{\ast}_{G(t)}\cap \wt V'{}^{\ast\beta '}$, defining in this way  covering morphisms $\wt g%_{\ssstyle\bullet}
^\beta $ and $\wt \Phi ^{\beta }%_{\ssstyle\bullet}
$ over $G_{|D^{\circ}}$, i.e.
$$
\begin{array}{ccc}
  \Pi _0^{\beta } & \stackrel{\Lambda}{\too} & \Pi _0'{}^{\beta '} \\
 \chi_{\beta } \downarrow\phantom{\chi_{\beta }} & {\sstyle\circlearrowright} & \phantom{\chi'_{\beta '}}\downarrow\chi'_{\beta '} \\
  D^\circ & \stackrel{G_{|D^\circ}}{\too} & D'{}^\circ
\end{array}\,,\quad \Lambda =\wt g%_{\ssstyle\bullet}
^{\beta }\hbox{ and }\Lambda =\wt \Phi ^{\beta }%_{\ssstyle\bullet}
\,,
$$
Notice that over the point  $c$, {the actions of} $\wt \Phi ^{\beta }%_{\ssstyle\bullet}
$ and $\wt g%_{\ssstyle\bullet}
^{\beta }$ coincide at the previously fixed point $\alpha \in \pi _0(\wt \Delta ^{\ast},\infty)\subset \Pi _0^\beta $. Indeed, we have the equalities:
$$\wt \Phi ^{\beta }%_{\ssstyle\bullet}
(\alpha )=\wt\phi %_{\ssstyle\bullet}
(\alpha ) =\wt g%_{\ssstyle\bullet}
^{\beta }(\alpha ),$$
the first one comes from the construction of $\wt\Phi ^{\beta }$ and the second one from the hypothesis  {(\ref{bullet})}.
%$(\bullet)_{\wt g,\,\wt \phi_{\infty} }$.
We deduce the identity  $\wt \Phi ^{\beta }%_{\ssstyle\bullet}
= \wt g%_{\ssstyle\bullet}
^{\beta }$ {on $\Pi_{0}^{\beta}$}. Using again   {(\ref{bullet})} %$(\bullet)_{\wt g,\,\wt \phi_{\infty} }$
we obtain that $\wt \Phi ^{\beta }%_{\ssstyle\bullet}
(\kappa  )=\wt g^{\beta }%_{\ssstyle \bullet}
(\kappa )=\wt\phi %_{\ssstyle\bullet}
(\kappa  )$, for all  ${\kappa \in \pi _{0}(\wt \Delta ^{\ast}\cap\wt V^{\ast\beta })}$, concluding thus the proof of  (\ref{i}).

\subsubsection{Proof of Assertion \hbox{\rm (\ref{ii})}}
Here  $\br D=D$ and  $K\cap\br D^\circ=D^\circ$. Assertions $(\ref{prem})$, $(\ref{hopfresp})$ and  $(\ref{extendephi})$ of Lemma~\ref{lemext} are satisfied by construction. The proof of (\ref{i}) above also shows Assertion  $(\ref{dercondi})$. It remains to show  $(\ref{cestunereal})$, i.e. the pro-germ  at infinity $\wt\Phi _\infty$ defined by $\wt\Phi $, makes commutative the diagrams {$(\star)$} and {$(\star\star)$}
%$(\star)_{\wt\Phi _\infty}$ and $(\star\star)_{\wt\Phi _\infty}$
in  Definition~\ref{defconj}.

\paragraph{\it Proof of the commutativity of {$(\star)$}.} %$(\star)_{\wt\Phi_\infty }$.}
We fix cofinal families  $U_n\in\mf U_{\F, \,\Sigma \cup E_S(T)}$ and $U'_n\in \mf U'_{\F',\Sigma \cup E_{S'}(T') }$, $n\in \mb N$, and we represent $h\in\mr{Hom}_{\proan}(\wt Q^{\F}_{\infty},\wt Q^{\F'}_{\infty})$ by a family of holomorphic maps $h_{n,p}$ from $\tQF{\wt U_n}$ into $\tQFpr{\wt U'_p}$, cf. (\ref{promorphismes}). By definition, the commutativity of  {$(\star)$} %$(\star)_{\wt\Phi_\infty }$
means the commutativity of all the diagrams
 $$   \begin{array}{ccc}
      (\wt V^\ast,\infty) & \stackrel{\wt \Phi_{\infty} }{\too} & (\wt V'{}^\star,\infty) \\
      \tau _n\downarrow \phantom{\tau _n}& &\phantom{\tau '_p}\downarrow \tau '_p\\
      \tQF{ U_n} & \stackrel {h_{n,p}}{\too}& \tQFpr{ U'_p}
    \end{array}\,,\qquad n\geq p\,,
$$
where $\tau _n:=\tau _{\ssstyle \wt V ^\ast\!,\,U_n}$ and $\tau '_p:=\tau _{\ssstyle \wt V'{}^\star\!\!\!,\, U'_p}$ are the canonical pro-germs defined in  Section~\ref{progermcan}. In order to prove this, thanks to Theorem~\ref{incfeuilles}, it suffices to determine open neighborhoods  $V_{n, p}$ of $D^\circ $ in $K_D$, {satisfying} %verifying
\begin{equation}\label{egsursystfond}
V_{n,p}\subset U_n\,,\quad \Phi (V_{n,p})\subset U'_p\quad\hbox{and}\quad \tau '_{n,p}\circ \wt\Phi _{|\wt V_{n,p}} = {h_{n,p}\circ \tau_{n,p} }_{|\wt V_{n,p}}\,,
\end{equation}
where
$$\tau _{n,p} : \wt V_{n,p}^\ast\too \tQF{ U_n} \quad\hbox{and}\quad \tau '_{n,p} : \wt V'{}^\star_{n,p}\to \tQFpr{ U'_p}$$
denote the quotient maps.
Take for   $V_{n,p}$ an open set  $W_{\kappa (n,p)}$ of the fundamental system given by $(\blacktriangle)$, where the index $\kappa (n,p)\in \mb N$ is chosen big enough so that the above inclusions hold.
The equality in (\ref{egsursystfond}) is an equality of maps and it can be checked locally using the open cover $(\wt V_{n,p}^\sigma )_{\sigma }$ de $\wt V_{n,p}^*$.
{Notice that such an argument is not valid if we want to check the equality of pro-germs $\tau '_p\circ\wt\Phi _\infty=h_{n,p}\circ\tau _n$
 without realizing them previously over the open sets $V_{n,p}$, cf. Remark~\ref{attention}.}  Put
 $$
\wt V_{n,p}^\sigma  :=\Sat_{\wt{\un{\F}}}(\wt\Delta^{\sigma} _{n,p}, \wt V_{n,p} ) \,, \quad \sigma \in \pi _0(\wt\Delta _{n,p}^*,\infty)\,,\quad \Delta _{n,p}:=V_{n,p}\cap T\,,
 $$
 where $\wt\Delta^{\sigma} _{n,p}$ is the connected component of $\wt\Delta _{n,p}^*$ corresponding to  $\sigma $.
Property  3) in $(\blacktriangle)$ implies that the open sets  $\wt V_{n,p}^\sigma$ cover $V_{n,p}$.
On the other hand,  Property {\it c)} in the  definition of $(\F,\Sigma \cup E_S(T))$-admissibility of $V_{n,p}$ stated in  Section~\ref{subsec.bons.vois}, implies that each leaf of the restriction of  $\wt{\un\F}$ to $\wt V_{n,p}^\sigma $ meets  (transversely) $\wt\Delta _{n,p}^\sigma $ in exactly one point, defining thus a holomorphic submersion-retraction-first integral
$r^\sigma $ from $\wt V_{n,p}^\sigma $ onto $\wt\Delta _{n,p}^\sigma $. By using  $(\blacktriangle)$ 4), it is clear that $\wt\Phi (\wt\Delta _{n,p}^\sigma )$ is a connected component of
$\wt \Delta '{}^\star_{n,p}\,$, where $\Delta '_{n,p}:=V'_{n,p}\cap T'=\Phi (\Delta_{n,p})$. Analogously, we also have a holomorphic submersion-retraction-first integral
 $r'{}^\sigma $ defined over the saturated set $\wt V'{}^\sigma_{n,p} $ of  $\wt\Phi (\wt\Delta _{n,p}^\sigma )$ in $\wt V'{}^\star_{n,p}$ onto $\wt\Phi (\Delta _{n,p}^\sigma )$.
 Clearly  $\wt V'{}^\sigma_{n,p} =\wt\Phi (\wt V'_{n,p})$ and
 the following diagram commutes
\begin{equation}\label{comretract}
    \begin{array}{ccc}
      \wt V_{n,p}^\sigma  & \stackrel{\wt\Phi _{|\wt V_{n,p}^\sigma}}{\too} & \wt V'{}^\sigma_{n,p}\\
      r^\sigma \downarrow \phantom{r^\sigma }&  \circlearrowright & \phantom{r'{}^\sigma }\downarrow r'{}^\sigma\\
     \wt \Delta _{n,p}^\sigma  & \stackrel{\wt\Phi _{|\wt \Delta _{n,p}^\sigma}}{\too} & \wt\Phi (\wt \Delta _{n,p}^\sigma)
    \end{array}
    \,,
\end{equation}
thanks to  $(\blacktriangle)$  6).
On the other hand, $\wt \Phi $ and $\wt \phi $ coincide on  $\wt \Delta ^\ast$. Therefore, the commutativity of the diagram {$(\star)$ for $\psi=\phi$}
%$(\star)_{\phi }$
in Definition~\ref{defconj} implies the commutativity of the following diagram if the index
$\kappa (n,p)$ is chosen big enough:
\begin{equation}\label{commparhyp}
\begin{array}{ccc}
 \wt \Delta _{n,p}^\sigma  & \stackrel{\wt\Phi _{|\wt \Delta _{n,p}^\sigma}}{\too} & \wt\Phi (\wt \Delta _{n,p}^\sigma)\\
  \tau \downarrow \phantom{\tau }& \circlearrowleft & \downarrow\tau ' \\
  \tQF{U_n} & \stackrel{h_{n,p}}{\too} & \tQFpr{U'_p}
\end{array}\,,
\end{equation}
$\tau $ and $\tau '$ denoting the quotient maps as usual.
It is clear that the commutativity of  the diagrams (\ref{comretract}) and (\ref{commparhyp}), gives us the  relation
$$\tau '_{n,p}\circ \wt\Phi _{|\wt V_{n,p}} = {h_{n,p}\circ \tau_{n,p} }_{|\wt V_{n,p}}$$ of (\ref{egsursystfond}). This finish the proof of the commutativity of {$(\star)$}. %$(\star)_{\wt\Phi _{\infty}}$.
\\

\paragraph{\it Proof of the commutativity of {($\star\star$) for $\psi=\Phi_{\infty}$.}} %$(\star\star)_{\wt\Phi_\infty }$.}
Notice that
$\wt \phi _\infty$ makes commutative the corresponding diagram  {($\star\star$)} %$(\star\star)_{\wt\phi _\infty}$
in Definition~\ref{defconj}:
$$%(\star\star)_{\wt\phi _\infty}\hspace{3em}
\begin{array}{ccc}
\Gamma_\infty & \stackrel{\iota}{\hookrightarrow} & \Gamma _{\wt T^\ast, \,\infty }\\
{\mf g}\downarrow\phantom{\phi_{*}} &
\hspace{-0,2cm}{\sstyle\circlearrowright}
 & \phantom{\wt\phi_{\ast}}\downarrow { \wt\phi_{\infty\ast}}\\
\Gamma_\infty' &
\hspace{-0,2cm}
\stackrel{\iota'}{\hookrightarrow} & \Gamma_{\wt T'{}^\star,\, \infty }
\end{array}\,,
\phantom{(\star\star)_{\wt\phi _\infty}\hspace{3em}}$$
where
$\iota(\varphi ):=\varphi _{|\wt T^\ast}$ (resp. $\iota'(\varphi ):=\varphi _{|\wt T'{}^\star}$).
By applying the following sub-lemma with $W_1:=E_S(T)$, $W_2:=E_S(V)$, $W'_1:=E_{S'}(T')$ and $W'_2:=E_{S'}(V')$,
we obtain directly the commutativity of  {($\star\star$) in our context.} %$(\star\star)_{\wt\Phi_\infty }$.

\begin{sublema}\label{souslemdiastarstar}
Let $W_1\subset W_2$ (resp. $W'_1\subset W'_2$), be submanifolds of ${\mb B} $ (resp. ${\mb B}'$) not contained in  $S$ (resp. $S'$) such that $W_1\cap S\neq\emptyset$ (resp. $W'_1\cap S'\neq\emptyset$). Consider $\wt \Psi _2:\wt W_2\to \wt W'_2$ a lift of a homeomorphism  $\Psi _2 : W_2\to W_2'$ such that $\Psi_{2} (W_1)=W'_1$. We denote   $\wt\Psi _1:=\wt{\Psi}_{2{}_{|\wt W_1^\ast}}$ and we keep the notations introduced in Section~\ref{subsectrealconj}. Then the
following diagram  {($\star\star$) commutes for $k=1$ if and only if it commutes for $k=2$:}
%diagram $(\star\star)_{\wt\Psi _1 }$ commutes if and only if the diagram $ (\star\star)_{\wt\Psi _2} $ commutes, where
      $$
     % {\sstyle (\star\star)}_{\wt\Psi _k }\hspace{1em}
\begin{array}{ccc}
\Gamma_\infty & \stackrel{\iota_k}{\hookrightarrow} & \Gamma _{\wt W_k^\ast,\,\infty }\\
{\scriptstyle\mf g}\downarrow\phantom{\scriptstyle\phi_{*}} & & \phantom{{\sstyle\wt \Psi _{k\ast}}}\downarrow {\sstyle \wt\Psi _{k\ast}}\\
\Gamma_\infty' &
\hspace{-0,2cm}
\stackrel{\iota'_k}{\hookrightarrow} & \Gamma_{\wt W_k'{}^\star,\,\infty }'
\end{array}\,,\quad k=1,2,
      $$  where    $\iota _k(\varphi ):=\varphi _{|\wt W^\ast_k}$ and  $\iota '_k(\varphi ):=\varphi _{|\wt W'{}^\ast_k}$.
      \end{sublema}

\begin{dem2}{of the sub-lemma}
It suffices to observe that the horizontal lines of the diagram below are composition of  monomorphisms:
$$
\begin{array}{ccccc}
\Gamma_\infty & \stackrel{\iota_2}{\hookrightarrow} & \Gamma _{\wt W_2^\ast,\,\infty } & \stackrel{\iota_{12}}{\hookrightarrow} & \Gamma _{\wt W_1^\ast,\, \infty }\\
{\scriptstyle\mf g}\downarrow\phantom{\scriptstyle\phi_{*}} & & \phantom{{\sstyle\wt \Psi _{2\ast}}}\downarrow {\sstyle \wt\Psi _{2\ast}} & \circlearrowleft & \phantom{{\sstyle \wt\Psi_{1\ast}}}\downarrow
 {\sstyle \wt\Psi_{1\ast}}\\
\Gamma_\infty' &
\stackrel{\iota'_2}{\hookrightarrow} & \Gamma'_{\wt W_2'{}^\star, \,\infty }&
\stackrel{\iota'_{12}}{\hookrightarrow} & \Gamma'_{\wt W_1'{}^\star,\, \infty }
\end{array}\,,
$$
with $\iota _{12}(\varphi ):=\varphi _{|\wt W^\ast_{1}}$ and $\iota '_{12}(\varphi ):=\varphi _{|\wt W'{}^\star_{1}}$.
\end{dem2}
This achieves the proof  Assertion $(ii)$ and consequently we have proved Lemma~\ref{lemext} in the case $K=K_{D}$.
%\end{dem2}

\vspace{1em}

\subsection{{Proof of  Extension Lemma~\ref{lemext} for $K=K_s$}}
\label{8.5} We assume that  $y_s=0$ (resp. $y_{s'}=0$) is a reduced local equation of  $D$ (resp. $D'$) supporting the transverse fibers $T=\pi_{D}^{-1}(c)$ and $T'=\pi_{D'}^{-1}(c')$, $c'=G(c)$. We denote by $\br D$ (resp. $\br D'$)
the irreducible component of $\mc D_{S}$ (resp. $\mc D_{S'}$) meeting $K$ (resp. $K'$) whose reduced local equation is $x_{s}=0$ (resp. $x_{s'}=0$).
 We also adopt the following notations:
$$
\mb P_{\lambda,\mu}:=\{|x_s|\leq \lambda , |y_s|\leq \mu \}
\,,\quad
\mb P'_{\lambda,\,\mu}:=\{|x_{s'}|\leq \lambda , |y_{s'}|\leq\mu\}\,,
$$
$$\mb T_{\lambda ,\,\mu  }:=\{|x_s|\leq \lambda , |y_s|=\mu \}\,,\quad
\mb T'_{\lambda ,\,\mu  }:=\{|x_{s'}|\leq \lambda , |y_{s'}|=\mu \}\,,
$$
$$
\mf T_{\lambda ,\mu }:=\{|x_s|=\lambda , \,|y_s|\leq\mu \}\quad\hbox{and}\quad\mf T'_{\lambda , \mu }:=\{|x_{s'}|=\lambda ,\, |y_{s'}|\leq\mu \}\,.
$$
Denote by $X$ (resp.  $X'$) the real vector field tangent to  $\un\F$ (resp. $\un\F'$), whose flow writes as  $(x_{s},y_{s},t)\mapsto (F(x_{s}, y_{s},t) , e^ty_{s})$ (resp.  $(x_{s'},y_{s'},t)\mapsto (F'(x_{s'}, y_{s'},t) , e^ty_{s'})$).\\

\noindent First, we will construct  $\Phi $ and $\wt\Phi $ satisfying Assertions  $(\ref{prem})$-$(\ref{cestunereal})$ of Lemma~\ref{lemext}. Next we will modify these two homeomorphisms, without affecting Properties   $(\ref{prem})$-$(\ref{cestunereal})$, in order to satisfy also Assertion $(\ref{dercondi})$.

\subsubsection{First step: construction of $\Phi$} Thanks to Theorem~\ref{monodimplholo}, the germ  $\phi_S$ conjugates the holonomies of $\un\F$ and $\un\F'$. Therefore it extends to  neighborhoods of the punctured disks $$D_s^\diamond:=\{0<|x_s|\leq 1, y_s=0\}\,,\quad D'_{s'}{}^\diamond:=G(D_s^\diamond)=\{0<|x_{s'}|\leq 1, y_{s'}=0\}\,,$$
defining a unique germ of homeomorphism  $\Phi_{D_s^\diamond}$ of $(\mc W_s,D_s^\diamond)$ onto $(\mc W'_{s'}, D'_s{}^\diamond)$ which conjugates $\un\F$ to $\un\F'$
and commutes with the Hopf fibrations, i.e. $\pi _{D'}\circ \Phi _{D_s^\diamond}= G\circ {\pi _D}$.
Since the germ $\phi_{S}$ is holomorphic and $G$ is excellent, cf. Definition~\ref{excellent}, and consequently  holomorphic on a polydisk $\mb P_{\alpha ,\,\beta }$, we deduce
that the germ $\Phi _{D_s^\diamond}$ is holomorphic for $|x_s|\leq \alpha$ if $\alpha>0$ is  small enough.
In fact, following the construction given in \cite{MarMatMarq} we can assume that $G$ satisfies also the relations
 $$
 x_{s'}\circ G_{|\mb P_{\alpha ,\,\beta }}={x_{s}}_{|\mb P_{\alpha ,\,\beta }}\,,\qquad y_{s'}\circ G_{|\mb P_{\alpha ,\,\beta }}={y_{s}}_{|\mb P_{\alpha ,\,\beta }}\,.
 $$
The equality of the Camacho-Sad indices will allow us to extend the germ $\Phi _{D_s^\diamond}$ to the singular point $s$ and this extension can be represented by a homeomorphism $\Phi _{D_s}$ defined on a polydisk  $\mb P_{1, \,\epsilon }$, $0<\epsilon <\alpha $ into an open set containing
 another polydisk  $\mb P'_{1,\,\epsilon '}$. This homeomorphism $\Phi_{D_{s}}$ conjugates  $\un \F$ to $\un \F'$ and verifies $\pi _{D'}\circ \Phi _{D_s}= {G\circ {\pi _D}}_{|\mb P_{1,\epsilon }}$. In addition,   $\Phi _{D_s}$ is holomorphic on  $\mb P_{\alpha ,\,\epsilon }$.

When $s$ is not a nodal singular point of $\uF$,  we can apply the {conjugation} Theorem stated in  \cite{MatMou} (see also   \cite[\S5.2.1]{Loray}) to construct $\Phi_{D_{s}}$. Otherwise, the nodal singularities are linearizable and  these linearizations can specifically been done without changing the  Hopf fibrations of $D$ and $D'$. Indeed, there exist a coordinate $\wt y_s$ (resp. $\wt y_{s'}$) such that $\uF$ (resp. $\uF'$) is given a linear vector field in the coordinates $(x_s, \wt y_s)$ (resp. $(x_s, \wt y_{s'})$). On the transversals $T$ and $T'$ the holonomy maps around the singularities $s$ and $s'$ are  irrational rotations when expressed in the coordinates $\wt y_{s\,|T}$ and $\wt y_{s'\,|T'}$ respectively. Thus the {conjugation} germ $\Phi_S$ is linear in these coordinates and any linear extension is a {conjugation} between $\uF$ and $\uF'$. In both cases up to restricting  $\alpha $ and $\epsilon >0$, we can also assume that
\begin{itemize}
  \item the vector field  $X$ is transverse to the torus  $\mb T_{\alpha ,\,\epsilon }$ and   $X'$ is transverse to the real analytic hypersurface  $H :=\Phi _{D_s}(\mb T_{\alpha ,\,\epsilon })$,
  \item $H$ does not intersect the torus $\mb T'_{1,\beta }$,
  \item for $\alpha '>0$ small enough, $(\mb P'_{\alpha ',\,1 }\setminus H)$ possesses two connected components and that one not containing $s$, does contain $\mb T_{\alpha ',\beta }$; we denote by  $\mb P'{}^+_{\alpha ',\,1 }$ its adherence.
\end{itemize}
We fix  such a $\alpha '$  and we ``glue''  $\Phi _{D_s}$ with  $\Phi _{\mc C}$ given by applying the following sub-lemma.
We denote by  $D_1$ (resp. $D_1'$) the irreducible component of $\mc D_S$ (resp. $\mc D_{S'}$), containing  $s$ (resp. $s'$), which is different from  $D$ (resp. $D'$).
\begin{sublema} There is
a homeomorphism $\Phi _{\mc C}$ defined on an open neighborhood  $V_{\mc C}$ of the annulus  $\mc C:=\{x_s=0, \epsilon \leq |y_s|\leq 1 \} $  in $\{|x_s|\leq\alpha ', \epsilon \leq |y_s|\leq 1  \}$, into  an open neighborhood $V'_{\mc C'}$ of  $\mc C':=\mb P'{}^+_{\alpha ', 1 }\cap\{x_{s'}=0\}$ in $\mb P'{}^+_{\alpha ', 1 }$,  conjugating  $\un\F_{V_{\mc C}}$ to $\un\F'_{|V'_{\mc C'}}$, which coincides with  $\Phi _{D_s}$ when restricted to  $V_{\mc C}\cap\mb T_{1,\epsilon }$, satisfying the relation  $\pi _{D_1'}\circ{\Phi _{\mc C}}(p)=G\circ{\pi _{D_1}}(p)$, for $|y_s(p)|\geq \beta $ and verifying  $\Sat_{\un\F}(V_{\mc C}\cap \mb T_{1,\epsilon }, V_{\mc C})=V_{\mc C}$.
\end{sublema}
The proof {can be done} %follows easily by
using  the classical lifting path method {and the vector fields $X$ and $X'$}. %We leave the details to the reader.

Thus, we obtain a homeomorphism  $\un\Phi $ defined on  $\mb P_{1,\epsilon }\cup V_{\mc C}$. We choose $\zeta>0$ small enough so that the set $V:=\{|x_s y_s|\leq \zeta\}\subset\mc W_s$ is contained in $\mb P_{1,\epsilon }\cup V_{\mc C}$ and we put $\Phi:=\un\Phi_{|V}:V\to V':=\un\Phi(V)$.

\subsubsection{Second step: construction of $\wt\Phi $}  With Conventions~\ref{convenrevecl}, we have that the natural maps
\begin{equation}\label{isotores}
 \pi _0(\wt{\mf T}^\ast_{1 , 1}, \infty)\too\pi _0(\wt V^\ast,\infty)\,,\quad\pi _0(\wt{\mf T}'{}^\ast_{1 , 1}, \infty)\too\pi _0(\wt V'{}^\ast,\infty)\,,
\end{equation}
induced by the inclusions  $\wt{\mf T}^\ast_{1 , 1}\subset\wt V^\ast$ and $\wt{\mf T}'{}^\ast_{1 , 1}\subset \wt V'{}^\ast$, are bijective.
Hence, each lift $\wt\Phi _{\mf T_{1 , 1}} : \wt{\mf T}^\ast_{\zeta , 1}\to \wt V'{}^\star$  of the restriction of  $\Phi$ to $\mf T^\ast_{\zeta , 1}$,  extends in a unique way to a homeomorphism  $\wt \Phi :\wt V^\ast\to\wt V'{}^\star$ lifting  $\Phi $. In order to choose  $\wt\Phi _{\mf T_{\zeta , 1}}$, we first apply the Lemma~\ref{lemext} in the case  $K=K_D$, for which it is already proved.
We obtain a homeomorphism $\Phi _{D^\circ}$ defined on a neighborhood  $V_{D^\circ}$ of $D^\circ$ and a lift  $\wt \Phi _{D^\circ}$ defined on  $\wt V_{D^\circ}^\ast$
satisfying the assertions of the lemma. Since  Conditions $(\ref{hopfresp})$ and $(\ref{extendephi})$ of  Lemma~\ref{lemext} gives the unicity of these homeomorphisms  and $\Phi$ also satisfies $(\ref{extendephi})$, we deduce that its restriction to $\mf T_{\zeta , 1}$ coincide with
 $\Phi _{D^\circ}$, after taking germs.
Therefore we can define $\wt\Phi _{\mf T_{\zeta , 1}}$ as the extension of $\wt{\Phi}_{D^\circ}$ onto $\wt V^\ast$, up to restricting   $\zeta>0$ so that   $\mf T_{\zeta , 1}$ be contained in  $V_{D^\circ}$.

\subsubsection{Third step:  proof of Assertions $(\ref{prem})$-$(\ref{cestunereal})$} Assertions $(\ref{prem})$-$(\ref{extendephi})$ are trivially satisfied.
In order to show   $(\ref{cestunereal})$ we must prove the commutativity of the two diagrams appearing in  Definition~\ref{defconj}.
The  commutativity of  {($\star\star$)} %$(\star\star)_{\wt\Phi_\infty }$
follows exactly in the same way that in the case $K=K_D$.
It only remains to prove the commutativity of {($\star$)}. %  $(\star)_{\wt\Phi_\infty }$.
Using the previous notations we remark that since the triple $(\Phi _{D^\circ}, \wt\Phi _{D^\circ}, h)$ is a realization of the $\mc N$-{conjugation}  $(\mf g,{h})$ over $V$, its restriction  $(\Phi _{|\mf T_{\zeta , 1}}, \wt\Phi _{|\wt{\mf T}^\ast_{\zeta , 1}}, h)$ is a realization of  $(\mf g,{h})$ over $\mf T_{\zeta , 1}$. Hence, we have the following commutative diagram
\begin{equation}\label{dia}
%{\sstyle (\star)}_{\wt\Phi _{|\wt{\mf T}_{\zeta , 1}^\ast}}\hspace{1em}
 \begin{array}{ccc}
(\wt {\mf T}_{\zeta , 1}^\ast, \infty) & \stackrel{\wt\Phi _{|\mf T_{\zeta , 1}}}{\longrightarrow} & (\wt {\mf T}'{}^\star_{\zeta,1}, \infty)\\
{\tau}\downarrow\phantom {\tau} & {\sstyle\circlearrowright}&\phantom{\scriptstyle\tau }\downarrow
{\tau } \\
\QF & \stackrel{h}{\longrightarrow} & \QFpr
\end{array}\,,
\end{equation}
where
$\tau $ and $\tau '$ denote the canonical pro-germs.
On the other hand, using the techniques appearing in  \cite[\S4.2]{MarMat} in the case of a resonant singularity, and a straightforward computation in the case of a linearizable singularity (excluding the case of a nodal singularity by hypothesis in Lemma~\ref{lemext}), {we can}
%it is easy to
prove the following result.
\begin{sublema}\label{dulaccollier} There are open neighborhoods
$U$ of $\mc D_S\cap K$ in $K$ and  $U'$ of   $\mc D_{S'}\cap K'$ in  $K'$ and there are deformation retractions $R: U^\ast\to U\cap{\mf T}^\ast_{1,1}$ and $R': U'{}^\star\to U'\cap{\mf T}'{}^\star_{1, 1}$  such that if we denote  $U_\lambda :=(\mc D_S\cap K)\cup R^{-1}( {\mf T}^\ast_{1,\lambda  }))$  and  $U'_\lambda :=(\mc D_{S'}\cap K')\cup R'{}^{-1}( {\mf T}'{}^\star_{1,\lambda  }))$, then
\begin{enumerate}
\item the family $(U_\lambda )_{0<\lambda \ll 1}$ (resp. $(U'_\lambda )_{0<\lambda \ll 1}$) is a fundamental system of neighborhoods of $\mc D_S\cap K$ in $K$ (resp. of $\mc D_{S'}\cap K'$ in $K'$);
    \item every  point  $p$ of $U_\lambda $ (resp. $U'_\lambda $) belongs to the same leaf of $\un\F_{|U_\lambda ^\ast}$ (resp. $\un \F'_{|U'{}^\star_\lambda }$) that  $R(p)$, (resp. $R'(p)$).
\end{enumerate}
\end{sublema}
\noindent  Thanks to Properties (1) and (2) in Sub-Lemma~\ref{dulaccollier} above, the (unique) lifts of these retractions define
germs at infinity $R_\infty$ and $R'_\infty$ which  make commutative the following diagrams
$$
\begin{array}{cc}
\xymatrix{(\wt{V}^*,\infty) \ar[rr]^{R_{\infty}}\ar[rdd]_{\tau_{\wt{V}^*}} & & (\tilde{\mathfrak{T}}^*,\infty) \ar[ldd]^{\tau_{\wt{\mathfrak{T}}^{*}}}\\
& & \\
& \wt{\mathcal{Q}}_{\infty}^{\mathcal{F}} &
} &
\xymatrix{(\wt{V}'^\star,\infty) \ar[rr]^{R_{\infty}'}\ar[rdd]_{\tau_{\wt{V}'^\star}} & & (\wt{\mathfrak{T}}'^\star,\infty) \ar[ldd]^{\tau_{\wt{\mathfrak{T}}'^{\star}}}\\
& & \\
& \wt{\mathcal{Q}}_{\infty}^{\mathcal{F}'} &
}
\end{array}
$$ as well as
 $$   \begin{array}{ccc}
      (\wt V^\ast,\infty) & \stackrel{\wt{\Phi} _\infty}{\longrightarrow} & (\wt V'{}^\star,\infty) \\
   {\sstyle  R_\infty}\downarrow \phantom{\sstyle R_\infty}& \circlearrowright &\phantom{\sstyle R'_{\infty}} \downarrow {\sstyle R'_{\infty}} \\
      (\wt {\mf T}^\ast, \infty) & \stackrel{\wt\Phi _{|\wt{\mf T}^\ast\infty}}{\longrightarrow} & (\wt {\mf T}'{}^\star, \infty)
    \end{array}\,.
$$
The commutativity of  {($\star\star$)} %$(\star)_{\wt \Phi }$
follows directly from the commutativity of this three diagrams and that of (\ref{dia}):

$$
\begin{array}{c}
\xymatrix{ (\wt{V}^*,\infty)\ar[rr]^{\wt{\Phi}_{\infty}}\ar[rrd]^{R_{\infty}}\ar[rddd]^{\tau_{\wt{V}'^{\star}}} &  &  (\wt{V}'^\star,\infty)\ar[rrd]^{R_{\infty}'}\ar[rddd]^{\tau_{\wt{V}'^{\star}}} &  & \\
& & (\wt{\mathfrak{T}}^*,\infty)\ar[ldd]^{\tau_{\wt{\mathfrak{T}}^*}}\ar[rr]^{\wt{\Phi}_{|\wt{\mathfrak{T}}^*\infty}} & & (\wt{\mathfrak{T}}'^\star,\infty)\ar[ldd]^{\tau_{\wt{\mathfrak{T}}'^{\star}}} \\
& & & & \\
& \wt{\mathcal{Q}}_{\infty}^{\mathcal{F}}\ar[rr]^{h} & & \wt{\mathcal{Q}}_{\infty}^{\mathcal{F'}} & }
\end{array}
$$

\subsubsection{Fourth step: modification of  $\Phi $ and $\wt \Phi $  to satisfy  (\ref{dercondi})} A priori the homeomorphism $\wt \Phi $ that we have constructed does not have to satisfy   (\ref{dercondi}). However, it induces on $\pi _0(\wt{ \mb T}^\ast_{1,1} ,\infty)$ the same map as $\wt g$,
\begin{equation}\label{egactionD}
\wt \Phi %_{\diamond}
=\wt g%_{\diamond}
 : \pi _0(\wt{ \mb T}^\ast_{1,1} ,\infty)\too\pi _0(\wt{ \mb T}'{}^\star_{1,1} ,\infty)\,.
\end{equation}
Indeed, $\wt \Phi $ and $\wt g$ induce the same map from $\pi _0(\wt{\mf T}^\ast_{1 , 1},\infty)$ onto $\pi _0(\wt{\mf T}'{}^\ast_{1 , 1},\infty)$, because $\wt \Phi _{D^\circ}$ satisfies (\ref{dercondi}). Then
the equality (\ref{egactionD}) follows from the bijections (\ref{isotores}) and
\begin{equation}\label{isotorepr}
\pi _0(\wt{\mb T}^\ast_{1,1}, \infty)\iso \pi _0(\wt V^\ast,\infty)\quad\hbox{and}\quad \pi _0(\wt{\mb T}'{}^{\star}_{1,1},\infty)\iso \pi _0(\wt V^\star,\infty)\,.
\end{equation}
We will compose $\wt\Phi$ to the left with the germ at infinity $\wt \Theta_\infty$
of the lift of a ``Dehn twist'' $\Theta $ along the leaves of $\un\F'$ with support contained in a small neighborhood of $\mb T'_{1,1}$:
$$
\begin{array}{ccc}
  (\wt V'{}^\star,\infty) & \stackrel{\wt \Theta_\infty }{\too} & (\wt V'{}^\star,\infty) \\
  \un q_{\infty}'\downarrow\phantom{\un q_{\infty}'} & \circlearrowleft & \phantom{\un q_\infty'}\downarrow\un q_\infty' \\
  (V',\mc D_S) & \stackrel{\Theta _S}{\too} & (V',\mc D_S)
\end{array}\,.
$$
We will see then that $\wt\Theta _\infty\circ\wt\Phi$ satisfies $(\ref{dercondi})$ as well as $(\ref{prem})$-$(\ref{cestunereal})$, already verified by  $\wt \Phi $.

\vspace{1em}

\paragraph{\textit{A) Construction of $\Theta $ and $\wt\Theta $. }}
Recall that  $\br D'$ is the irreducible component of $\mc D_{S'}$ whose local equation is  $x_{s'}=0$. On the disk $\br D'\cap\mc W'_{s'}$, we consider the real vector field   $\vartheta$ whose flow is  $(t, y_{s'})\mapsto e^{2i\pi t}y_{s'}$. We fix a smooth function   $u$ with support contained in $\{x_{s'}=0, \varsigma\leq |y_{s'}|\leq 1\}$, $0<\varsigma<1$ taking the value 1 on $C':= \{x_{s'}=0,\ |y_{s'}|=1\}=K'\cap \br D'^{\circ}$ and we denote by   $Y$ the vector field on $V'$   tangent to   $\un\F'$ and projecting over $u\vartheta$ by $\pi _{\br D'}$. The flow
$\Upsilon_t$ of $Y$ is defined on a open neighborhood $U_{I}$ of $\br D$ in  $\mc W_{s'}$, once we fix an interval $I\subset\mb R$ where we allow the time $t$ to vary.
Hence, we can lift the flow $\Upsilon_{t}$ to a unique map $\wt\Upsilon_t : \wt U_I^\star\to \wt V'{}^\star$ being the identity on  $|y_{s'}\circ\un q'|\leq \varsigma$, and defining consequently a germ at infinity  $\wt\Upsilon_{t\infty} : (\wt V'{}^\star,\infty)\to(\wt V'{}^\star,\infty)$.
Clearly the germ  $\wt\Upsilon_{n\infty}$ fibers over  $ C'$, i.e. ${\pi _{D'}\circ\un q'\circ\wt\Upsilon_{n\infty}}_{|\wt{\mb T}'{}^\star}={\pi _{D'}\circ\un q'_{\infty}}_{|\wt{\mb T}'{}^\star}$, for each $n\in\mb Z$.
It defines a deck transformation $\wt\Upsilon_{n%\ssstyle\bullet
} :\Pi '\iso\Pi '$ of  the natural covering
$$
\rho' :\Pi' \to C'\,,\quad \rho '{}^{-1}(p):=\pi _0(\wt T'_p{}^\star,\infty)\,,\quad T'_p:=\pi _{\br D'}^{-1}(p)\,.
$$
We put $\Theta :=\Upsilon_{n_0}$ and $\wt\Theta :=\wt\Upsilon_{n_0}$,
choosing the integer $n_0$ in the following way. First,
 we fix a point $a$ in the circle $C:=\{x_s=0, |y_s|=1\}=K\cap\br D^\circ$ and an element  $\nu$ of $\pi _0(\wt T_{a}^{\star},\infty)$, $T_a:=\pi _{\br D}^{-1}(a)$. We consider the natural covering  $\rho :\Pi \to C$ of $C$  with fibers $\rho ^{-1}(p):=\pi _0(\wt T^\ast_p,\infty)$ and  the following two covering morphisms over $G_{|C}$
      $$
(\blacklozenge )_{\Lambda }
\quad
  \begin{array}{ccc}
    \Pi  & \stackrel{\Lambda }{\too} &  \Pi '\\
    \rho \downarrow\phantom{\rho } & \circlearrowleft & \phantom{\rho '}\downarrow\rho ' \\
    C & \stackrel{G_{|C}}{\too} & C'
  \end{array}\,,\qquad \quad\Lambda = \wt \Phi %_{\ssstyle\bullet}
  \quad\hbox{and}\quad \Lambda  =\wt g%_{\ssstyle\bullet}
  \,,
  $$
  defined by $\wt \Phi $ and $\wt g$.
After the equality  (\ref{egactionD}) and the  bijections (\ref{isotores}) and (\ref{isotorepr}), the images $\sigma ':=\wt \Phi %_{\ssstyle\bullet}
(\nu )$  and $\sigma '':=\wt g%_{\ssstyle\bullet}
(\nu)$ belong to the same fiber of the map $\iota %_{\ssstyle \bullet}
 : \pi _0(\wt T_{a'}'{}^{\star},\infty)\to\pi _0(\wt{\mb T}'{}^\star_{1,1},\infty)$, induced by the inclusion map  $\iota  : \wt T'_{a'}{}^{\star}\hookrightarrow \wt{\mb T}'{}^\star_{1,1}$. On the other hand, the action $$
\mb Z\times \pi _0(\wt T'_{a'}{}^{\star},\infty) \too \pi _0(\wt T'_{a'}{}^{\star},\infty)\,\qquad(n,\sigma )\mapsto \wt\Upsilon_{n%\ssstyle \bullet
}(\sigma )\,,
$$
of $\mb Z$ on the  fiber  of $\rho '$ at the point $a':=\Phi (a)=G(a)\in C'$ coincide with the action of  $\pi _1(C', a')\simeq \mb Z$ induced by the covering $\rho'$ on this fiber.
The orbits of that action correspond to the fibers of  $\iota %_{\ssstyle \bullet}
$.
We choose $n_0$ to be the unique integer number such that $\Upsilon_{n_0%\ssstyle\bullet
}(\sigma ')=\sigma ''$.

\vspace{1em}

\paragraph{\textit{B) Proof of Assertions   $(\ref{prem})$-$(\ref{dercondi})$ for   $\wt \Theta \circ\wt\Phi _\infty$.}}
Properties  $(\ref{prem})$-$(\ref{extendephi})$ are trivially satisfied.  Property  $(\ref{cestunereal})$ follows from Proposition~\ref{lemmeisotopie}. It only remains to prove  $(\ref{dercondi})$, that is, to show the equalities
$$(\lozenge)_t\quad(\Theta \circ\wt\Phi _{\wt T^\ast_{t}})%_{\ssstyle\bullet} 
= \wt g_{\wt T^\ast_{t}%\ssstyle\bullet
} : \pi _0(\wt T^\ast_{t},\infty)\too\pi _0(\wt T'{}^\star_{G(t)},\infty)\,,$$
for all $t\in C$
 and for all $t\in D^{\circ}\cap K$.
In the last case, $\wt\Theta$ is the identity near $\wt T_{t}'^{\star}$ and therefore
$\wt\Theta\circ\wt\Phi$ coincides with the homeomorphism $\wt\Phi_{D^{\circ}}$ given by Lemma~\ref{lemext} in the context $K=K_{D}$, which satisfies $(\ref{dercondi})$ as we have already see.
In the case $t\in C$,
 $(\lozenge)_t$ is equivalent to the equality  $(\wt\Theta \circ\wt\Phi) %_{\ssstyle\bullet}
 =\wt g%_{\ssstyle\bullet} 
 $ of the covering morphisms $(\blacklozenge)_{(\wt\Theta \circ\wt\Phi) %_{\ssstyle\bullet}
 }$ and $(\blacklozenge)_{\wt g%_{\ssstyle\bullet}
 }$ defined by $\wt\Theta \circ\wt\Phi _\infty$ and $\wt g_{\infty}$.
Thus, it suffices to show the equality on a single fiber, i.e. $(\lozenge)_a$. This equality is satisfied for the element $\sigma $ of $\pi _0(\wt T^\ast_{t},\infty)$ that we have previously fixed to define $n_{0}$: $(\wt\Theta \circ\wt\Phi) %_{\ssstyle\bullet}
(\sigma )=\wt g%_{\ssstyle\bullet}
(\sigma )$.
Since $\Gamma _\infty$  acts transitively on  $\pi _0(\wt T^\ast_{t},\infty)$, we must show the equality $(\wt\Theta  \circ\wt\Phi)%_{\sstyle\bullet}
(\varphi%_{\ssstyle\bullet}
(\sigma ))=\wt g%_{\sstyle\bullet}
( \varphi %_{\sstyle\bullet}
(\sigma )) $, for all $\varphi \in \Gamma _\infty$, which follows directly from the commutativity of   {($\star\star$) for $\psi=\Phi\circ\Phi_{|T_{a}}$}.
%$(\star\star)_{\wt \Theta \circ \wt\Phi _{|\wt T^\ast_a}}$.
Indeed,
\begin{gather}\label{calculsimple}
\begin{split}
(\wt\Theta  \circ\wt\Phi)%_{\sstyle\bullet}
(\varphi%_{\ssstyle\bullet}
(\sigma ))&=[(\wt\Theta \circ\wt\Phi )_{\ast}(\varphi )]%_{\ssstyle\bullet}
((\wt\Theta \circ\wt\Phi )%_{\ssstyle\bullet}
(\sigma ))=(\wt g_{\ast}(\varphi ))%_{\ssstyle\bullet}
(\wt g%_{\ssstyle\bullet}
(\sigma ))\\& =(\wt g_{\ast}(\varphi )\circ\wt g)%_{\ssstyle\bullet}
(\sigma )= (\wt g\circ \varphi )%_{\ssstyle\bullet}
(\sigma )=\wt g%_{\sstyle\bullet}
( \varphi %_{\sstyle\bullet}
(\sigma ))\,.
\end{split}
\end{gather}
This concludes the proof of Extension Lemma~\ref{lemext} for $K=K_s$.
%\end{dem2}

\end{document}